\documentclass[12pt]{amsart}
\usepackage{longtable}
\usepackage{amsfonts,amssymb,amscd,amsmath,epsfig}
\usepackage{latexsym}
\usepackage{mathrsfs}
\usepackage{tocvsec2}

\usepackage[]{youngtab}

\usepackage[
hyperindex=true,pagebackref=true,bookmarks=true,
colorlinks=true,linkcolor=blue,citecolor=red]{hyperref}

\let\vv\v

\begin{document}

{
\centering 
{\sf\em Dedicated to the memory of Yuri Ivanovich Manin
\medskip\par}
\vskip 0.3cm}

\begin{abstract}
We begin with modular form periods, a focal point of several
Yuri Manin's works. The similarity is discussed
between the corresponding
zeta-polynomials and 
superpolynomials of algebraic links, closely related
to Khovanov-Rozansky polynomials. 
We focus on DAHA superpolynomials and motivic ones, defined
via compactified Jacobians of plane curve singularities
and their counterparts in arbitrary ranks; the 
non-unibranch construction is new. 
They  conjecturally coincide with 
the corresponding generalizations of
{\em L}-functions and satisfy the Riemann
Hypothesis in some sectors of the parameters. 
Presumably, the motivic ones
can be interpreted as certain partition functions of the
Landau-Ginzburg model associated with plane curve singularities; 
RH for them is remarkably similar
to the Lee-Yang circle theorem for Ising models. 
A $q,t$-deformation of the Witten index is obtained as
an application. 
General perspectives of 
the motivic theory of isolated curve and surface singularities are
discussed, including possible implications in number theory.  
Also, we introduce super-analogs of $\rho_{ab}-$ invariants and
discuss super-deformations of the Riemann's zeta. Among other
topics: Verlinde algebras and the topological vertex.
\end{abstract}

\renewcommand{\tilde}{\widetilde}
\renewcommand{\hat}{\widehat}

\newcommand{\BR}{{\mathbb R}}
\newcommand{\BQ}{{\mathbb Q}}
\newcommand{\BC}{{\mathbb C}}
\newcommand{\BP}{{\mathbb P}}
\newcommand{\BZ}{{\mathbb Z}}
\newcommand{\BN}{{\mathbb N}}
\newcommand{\BS}{{\mathbb S}}

\newcommand{\cH}{{\mathcal H}}
\newcommand{\cA}{{\mathcal A}}
\newcommand{\cB}{{\mathcal B}}
\newcommand{\ccF}{{\mathfrak F}}
\newcommand{\cD}{{\mathcal D}}
\newcommand{\cL}{{\mathcal L}}
\newcommand{\cF}{{\mathcal F}}
\newcommand{\cP}{{\mathcal P}}
\newcommand{\cX}{{\mathcal X}}
\newcommand{\cY}{{\mathcal Y}}
\newcommand{\cS}{{\mathcal S}}
\newcommand{\cSol}{\hbox{$\mathcal Sol$}}
\newcommand{\cT}{\hbox{$\mathcal T$}}

\newcommand{\Z}{{\mathbb Z}}
\newcommand{\Q}{{\mathbb Q}}
\newcommand{\N}{{\mathbb N}}
\newcommand{\C}{{\mathbb C}}
\newcommand{\R}{{\mathbb R}}
\newcommand{\X}{{\mathbb X}}
\newcommand{\Y}{{\mathbb Y}}

\newcommand{\CH}{{\mathcal H}}
\newcommand{\CA}{{\mathcal A}}

\def\HH{\mbox{${\mathcal H}$\kern-5.2pt${\mathcal H}$}}

\newcommand{\binomial}[2]{\genfrac{(}{)}{0pt}{}{ #1 }{ #2 }}
\newcommand{\qbinomial}[2]{\genfrac{[}{]}{0pt}{}{ #1 }{ #2 }_q }
\newcommand{\qbinom}[3]{\genfrac{[}{]}{0pt}{}{ #1 }{ #2 }_{ #3 } }


\def\der{\partial}
\def\tensor{\otimes}
\def\gam{\gamma} \def\Gam{\Gamma}
\def\del{\delta} \def\Del{\Delta}
\def\kap{\kappa}
\def\lam{\lambda} \def\Lam{\Lambda}
\def\Comp{{\mathbb C}}
\def\sM{{\mathcal M}}

\newtheorem{theorem}{Theorem}[section]
\newtheorem{maintheorem}[theorem]{Main Theorem}
\newtheorem{proposition}[theorem]{Proposition}
\newtheorem{definition}[theorem]{Definition}
\newtheorem{lemma}[theorem]{Lemma}
\newtheorem{corollary}[theorem]{Corollary}
\newtheorem{notation}[theorem]{Notation}
\newtheorem{remark}[theorem]{Remark}
\newtheorem{example}[theorem]{Example}

\newtheorem{theorem }{Theorem}[section]
\newtheorem{maintheorem }[theorem]{Main Theorem}
\newtheorem{proposition }[theorem]{Proposition}
\newtheorem{definition }[theorem]{Definition}
\newtheorem{lemma }[theorem]{Lemma}
\newtheorem{corollary }[theorem]{Corollary}
\newtheorem{notation }[theorem]{Notation}
\newtheorem{remark }[theorem]{Remark}
\newtheorem{example }[theorem]{Example}

\newtheorem{ maintheorem }[theorem]{Main Theorem}
\newtheorem{ theorem}{Theorem}[section]
\newtheorem{ proposition}[theorem]{Proposition}
\newtheorem{ definition}[theorem]{Definition}
\newtheorem{ lemma}[theorem]{Lemma}
\newtheorem{ corollary}[theorem]{Corollary}
\newtheorem{ notation}[theorem]{Notation}
\newtheorem{ remark}[theorem]{Remark}
\newtheorem{ example}[theorem]{Example}

\newtheorem{thm}{Theorem}[section]
\newtheorem{prop}[thm]{Proposition}
\newtheorem{lem}[thm]{Lemma}
\newtheorem{cor}[thm]{Corollary}
\newtheorem{conj}[thm]{Conjecture}
\newtheorem{con}[thm]{Conjecture}
\newtheorem{dfn}[thm]{Definition}
\newtheorem{df}[thm]{Definition}
 \newcommand{\rem}{{\bf Comment.\ }}
 \newcommand{\rmk}{{\bf Comment.\ }}
 \newcommand{\exmp}{{\bf Example.\ }}
 \newcommand{\ex}{{\bf Example.\ }}
 \newcommand{\prob}{{\bf Problem.\ }}

\newtheorem{note}{Note} 
\renewcommand{\thenote}{}
\newtheorem*{acka}{Acknowledgments}
\newtheorem{ack}{Acknowledgments}
\renewcommand{\theack}{}
\renewcommand{\appendixname}{\bf Appendix}
\renewcommand{\proof}{{\em Proof.\ }}

\hyphenation{
ap-pen-dix as-ymp-tot-ic at-trib-uted at-trib-ut-able
Bry-li-n-sky com-mu-ta-tion de-ge-ne-rate
de-riv-a-tive dis-trib-ute equi-vari-ant ex-tra-or-di-nary  
geo-met-ric griev-ance griev-ous grad-ed ho-lo-no-my ho-mo-thetic
in-fin-ite-ly in-fin-i-tes-i-mal Ha-rish Cha-n-dra mul-ti-plic-able 
non-euclid-ean non-iso-mor-phic non-smooth par-a-digm 
par-a-bol-ic pa-rab-o-loid pa-ram-e-trize phe-nom-e-non 
post-script pseu-do-dif-fer-en-tial pseu-do-fi-nite 
qua-drat-ics quad-ra-ture Han-kel rec-tan-gle semi-def-i-nite 
set-up wide-spread Euler-ian Feb-ru-ary Gauss-ian Grothen-dieck 
Hamil-ton-ian Her-mi-t-ian her-mi-t-ian Jan-u-ary 
Japan-ese Ka-shi-wa-ra Kor-te-weg Le-gendre No-vem-ber Rie-mann-ian 
Sep-tem-ber Za-mo-lo-d-chi-kov Kni-zh-nik quan-tum Op-dam
Mac-do-nald Ca-lo-ge-ro Su-ther-land Mo-ser 
Ol-sha-net-sky  Pe-re-lo-mov in-de-pen-dent ope-ra-tors 
cy-clo-to-mic ra-tio-nal de-gen-er-a-tion 
in-ter-est-ing de-for-ma-tions de-for-ma-tion pro-ce-dure 
fol-lows ope-ra-tors  pre-serve suf-fices ap-proach 
for-mu-las con-sider its com-ple-tion cor-re-spond-ing 
au-to-mor-phism be-cause pro-por-tional fi-nal-ly let-ting 
equi-v-a-lence ge-n-er-al-ized Mac-do-nald iden-ti-ties 
cor-re-s-pond sub-dia-grams par-ti-tion na-t-u-ral-ly 
or-dered stan-dard de-for-ma-tion ar-gu-ment com-bined 
sphe-r-i-cal rep-re-sen-ta-tions tri-go-no-me-t-ric
ge-n-er-al-ly speak-ing pri-m-it-ive ir-re-du-cible 
sum-ma-tion  rep-re-sen-ta-tives pro-por-ti-o-na-li-ty
ultra-sphe-ri-cal Ro-gers}

\def\ffor{\quad\hbox{ for }\quad}
\def\wwhen{\quad\hbox{ when }\quad}
\def\wwhere{\quad\hbox{ where }\quad}
\def\aand{\quad\hbox{ and }\quad}
\def\for{\  \hbox{ for } \ }
\def\iif{ \ \hbox{ if } \ }
\def\when{ \ \hbox{ when } \ }
\def\where{\  \hbox{ where } \ }
\def\and{\  \hbox{ and } \ }
\def\and{\  \hbox{ and } \ }
\def\oor{\  \hbox{ or } \ }
\def\proof{{\em Proof. \  }}

\def\equal{\stackrel{\,\mathbf{def}}{= \kern-3pt =}}

\def\la{\lambda}
\def\La{\Lambda}
\def\om{\omega}
\def\Om{\Omega}
\def\Th{\Theta}
\def\th{\theta}
\def\al{\alpha}
\def\be{\beta}
\def\ga{\gamma}
\def\ep{\epsilon}
\def\up{\upsilon}
\def\Up{\Upsilon}
\def\de{\delta}
\def\De{\Delta}
\def\ka{\kappa}
\def\kapp{\hbox{\bf \ae}}
\def\si{\sigma}
\def\Si{\Sigma}
\def\Ga{\Gamma}
\def\ze{\zeta}
\def\io{\iota}
\def\bio{b^\iota}
\def\aio{a^\iota}
\def\twio{\tilde{w}^\iota}
\def\hwio{\hat{w}^\iota}
\def\gio{\g^\iota}
\def\Bio{B^\iota}

\def\del{\delta}
\def\pa{\partial}
\def\vp{\varphi}
\def\ve{\varepsilon}
\def\inf{\infty}

\def\vph{\varphi}
\def\vps{\varpsi}
\def\vPh{\varPhi}
\def\vep{\varepsilon}
\def\vpi{{\varpi}}
\def\vth{{\vartheta}}
\def\vsi{{\varsigma}}
\def\vrh{{\varrho}}

\def\bph{\bar{\phi}}
\def\bsi{\bar{\si}}
\def\bvp{\bar{\varphi}}

\newcommand{\bS}{{\mathbf S}}
\newcommand{\bH}{{\mathbf H}}
\newcommand{\bF}{{\mathbf F}}
\newcommand{\bE}{{\mathbf E}}

\def\tal{\tilde{\alpha}}
\def\tbe{\tilde{\beta}}
\def\tde{\tilde{\delta}}
\def\tpi{\tilde{\pi}}
\def\txi{\tilde{\xi}}
\def\tPi{\tilde{\Pi}}
\def\tPhi{\tilde{\Phi}}
\def\tV{\tilde{V}}
\def\tJ{\tilde{J}}
\def\tla{\tilde{\lambda}}
\def\tga{\tilde{\gamma}}
\def\tGa{\tilde{\Gamma}}
\def\tvs{\tilde{{\varsigma}}}
\def\tu{\tilde{u}}
\def\tU{\tilde{U}}
\def\tw{\widetilde w}
\def\tW{\widetilde W}
\def\tB{\tilde B}
\def\tv{\tilde v}
\def\tV{\tilde V}
\def\tz{\tilde z}
\def\tb{\tilde b}
\def\ta{\tilde a}
\def\tih{\tilde h}
\def\trh{\tilde {\rho}}
\def\tx{\tilde x}
\def\tf{\tilde f}
\def\tg{\tilde g}
\def\tG{\tilde G}
\def\tk{\tilde k}
\def\tl{\tilde l}
\def\tL{\tilde L}
\def\tD{\tilde D}
\def\tR{\tilde R}
\def\tP{\tilde P}
\def\tH{\tilde H}
\def\tp{\tilde p}

\def\hH{\hat{H}}
\def\hh{\hat{h}}
\def\hR{\hat{R}}
\def\hY{\hat{Y}}
\def\hX{\hat{X}}
\def\hP{\hat{P}}
\def\hT{\hat{T}}
\def\hV{\hat{V}}
\def\hG{\hat{G}}
\def\hF{\hat{F}}
\def\hw{\widehat{w}}
\def\hW{\widehat{W}}
\def\hu{\hat{u}}
\def\hs{\hat{s}}
\def\hv{\hat{v}}
\def\hb{\hat{b}}
\def\hB{\widehat{B}}
\def\hze{\hat{\zeta}}
\def\hsi{\hat{\sigma}}
\def\hrh{\hat{\rho}}
\def\hth{\hat{\theta}}
\def\hy{\hat{y}}
\def\hx{\hat{x}}
\def\hz{\hat{z}}
\def\hg{\hat{g}}
\def\he{\hat{e}}
\def\hE{\widehat{E}}

\def\B{\mathbf{B}}
\def\I{\mathbf{I}}
\def\P{\mathbf{P}}
\def\G{\mathbf{G}}
\def\S{\mathbf{S}}
\def\F{\mathbf{F}}
\def\one{\mathbf{1}}
\def\Sn{\mathbf{S}_n}
\def\0{\mathbf{0}}
\def\H{\mathbf{H}}
\def\V{\mathbf{V}}

\def\f{\mathcal{F}}
\def\çF{\mathcal{F}}
\def\o{\mathcal{O}}
\def\t{\mathcal{T}}
\def\r{\mathcal{R}}
\def\l{\mathcal{L}}
\def\m{\mathcal{M}}
\def\k{\mathcal{K}}
\def\n{\mathcal{N}}
\def\d{\mathcal{D}}
\def\p{\mathcal{P}}
\def\cP{\mathcal{P}}
\def\a{\mathcal{A}}
\def\h{\mathcal{H}}
\def\c{\mathcal{C}}
\def\y{\mathcal{Y}}
\def\e{\mathcal{E}}
\def\v{\mathcal{V}}
\def\z{\mathcal{Z}}
\def\x{\mathcal{X}}
\def\s{\mathcal{S}}
\def\g{\mathcal{G}}
\def\u{\mathcal{U}}
\def\w{\mathcal{W}}
\def\i{\mathcal{I}}
\def\j{\mathcal{J}}
\def\b{\mathcal{B}}

\def\lan{\langle}
\def\llb{(\!(}
\def\ran{\rangle}
\def\rrb{)\!)}
 \def\dim{{\hbox{\rm dim}}_{\mathbb C}\,}
\def\lng{\hbox{\rm{\tiny lng}}}
\def\sht{\hbox{\rm{\tiny sht}}}
\def\sph{\hbox{\rm{\tiny sph}}}
\def\inv{\hbox{\rm{\tiny inv}}}

\def\br#1{\langle #1 \rangle}

\def\rank{\hbox{rank}}
\def\gl{\mathfrak{gl}_N}

\newcommand{\Aut}{\operatorname{Aut}}
\newcommand{\Hom}{\operatorname{Hom}}
\newcommand{\End}{\operatorname{End}}
\newcommand{\Ind}{\operatorname{Ind}}
\newcommand{\ad}{\operatorname{ad}}
\newcommand{\pr}{\operatorname{pr}}
\newcommand{\aweyl}{\tilde{\mathbb S}_n}
\newcommand{\hec}{{\mathcal H}^t_n}
\newcommand{\Func}{{\mathcal F}({\mathbb C}^n,{\mathcal H}^t_n)}
\newcommand{\tr}{\operatorname{tr}}
\newcommand{\Out}{\operatorname{Out}}
\newcommand{\Rad}{\operatorname{Rad}}
\newcommand{\Spec}{\operatorname{Spec}}
\newcommand{\id}{\operatorname{id}}
\newcommand{\Int}{\operatorname{Int}}
\newcommand{\ct} {\operatorname{ct}}

\newcommand{\rat}{{\mathbb Q}}
\newcommand{\real}{{\mathbb R}}
\newcommand{\cplx}{{\mathbb C}}
\newcommand{\zint}{{\mathbb Z}}

\newcommand{\sq}{\phantom{1}\hfill$\qed$}
\newcommand{\Rea}{\Re}
\newcommand{\Ima}{\Im}

\newcommand{\st}{\bowtie}
\newcommand{\modd}{\mbox{\,mod\,}}
\newcommand{\lr}{\langle}
\newcommand{\rr}{\rangle}
\newcommand{\eps}{\varepsilon}
\newcommand{\phk}{\phi^{(k)}}
\newcommand{\psk}{\psi^{(k)}}
\newcommand{\Res}{\mbox{Res}\;}
\newcommand{\sgn}{\mbox{sgn}}
\newcommand{\mn} {\left\{ \begin{array}{c}m\\
n\end{array}\right\}}

\def\sX{\mathscr{X}}
\def\sH{\mathscr{H}}
\def\sY{\mathscr{Y}}
\def\TT{\mathfrak{T}}
\def\JJ{\mathfrak{J}}
\def\HH{\mathfrak{H}}
\def\FF{\mathfrak{F}}
\def\GG{\mathfrak{G}}
\def\CC{\mathfrak{C}}
\def\LL{\mathfrak{L}}

\def\BB{\mathfrak{B}}
\def\AA{\mathfrak{A}}
\def\ZZ{\mathfrak{Z}}
\def\HH{\hbox{${\mathcal H}$\kern-5.2pt${\mathcal H}$}}
\def\HHH{\hbox{${\mathbb H}$\kern-4.2pt${\mathbb H}$}}
\def\tHH{\widetilde{\HH\ }}

\font\smm=msbm10 at 12pt 
\def\symbol#1{\hbox{\smm #1}}
\def\lsmash{{\symbol n}}
\def\rsmash{{\symbol o}}
\def\#{\sharp}

\font\tenbf=cmbx10
\font\tenrm=cmr10
\font\tenit=cmti10
\font\ninebf=cmbx9
\font\ninerm=cmr9
\font\nineit=cmti9
\font\eightbf=cmbx8
\font\eightrm=cmr8
\font\eightit=cmti8
\font\sevenrm=cmr7
\font\sevenbf=cmbx7


\title [Zeta-polynomials, superpolynomials]
{Zeta-polynomials, superpolynomials, \\ DAHA
and plane curve singularities}
\author[Ivan Cherednik]{Ivan Cherednik $^\dag$}

\thanks{$^\dag$ February 18, 2024. 
\ \ \ Partially supported by NSF grant
DMS--1901796}

\address[I. Cherednik]{Department of Mathematics, UNC
Chapel Hill, North Carolina 27599, USA\\
chered@email.unc.edu}

\maketitle

{\small 
{\sf\em Key words:\ } double affine Hecke algebras; 
Macdonald polynomials; Verlinde algebras; triangular groups; 
Mostow groups; linear rigidity; Deligne-Simpson problem; 
knot invariants;
HOMFLY-PT polynomials; iterated torus links; Seifert manifolds;
Lens spaces; 
Alexander polynomials; Khovanov-Rozansky homology; rho-invariants;
compactified Jacobians; affine Springer fibers; 
plane curve singularities; compactified Jacobians; theta-functions;
Hasse-Weil zeta functions; 
Riemann's zeta function; Dirichlet L-functions; 
Riemann hypothesis; modular forms; Iwasawa theory; 
Landau-Ginzburg models; Ising models; Witten index;
topological vertex

}

\renewcommand*\contentsname{\vskip -0.0cm \sc Contents}
{\normalsize
\setcounter{tocdepth}{2}
\tableofcontents
\vfil\eject
}

 \def\q{\mathcal{Q}}
 \def\sht{\raisebox{0.4ex}{\hbox{\rm{\tiny sht}}}}
 \def\bysame{{\bf --- }}
\let\oldt\~   
 \def\~{{\bf --}}
 \def\rr{{\mathsf r}}
 \def\ss{{\mathsf s}}
 \def\mm{{\mathsf m}}
 \def\pp{{\mathsf p}}
 \def\ll{{\mathsf l}}
 \def\aa{{\mathsf a}}
 \def\bb{{\mathsf b}}
 \def\NS{\hbox{\tiny\sf ns}}
 \def\ssum{\hbox{\small$\sum$}}
\newcommand{\comment}[1]{}
\renewcommand{\tilde}{\widetilde}
\renewcommand{\hat}{\widehat}
\renewcommand{\V}{\mathbb{V}}
\renewcommand{\S}{\mathbb{S}}
\renewcommand{\F}{\mathbb{F}}
\newcommand{\dagx}{\hbox{\tiny\mathversion{bold}$\dag$}}
\newcommand{\ddagx}{\hbox{\tiny\mathversion{bold}$\ddag$}}
\newtheorem{conjecture}[theorem]{Conjecture}
\newcommand*{\vect}[1]{\overrightarrow{\mkern0mu#1}}



\section[\sc \hspace{1em}Modular form periods] 
{\sc Modular form periods}  
\subsection{\bf Manin, my teacher}
In 1965-67, Yuri Ivanovich Manin and
Ernest Borisovich Vinberg delivered special courses at
Moscow School no 2 for senior students. This is when I met Yu.I.
With some stretch, I can say that Manin and Vinberg 
were my high school teachers and Vasilii Iskovskikh and
Victor Kac were our tutors (teaching assistants of Yu.I. and
E.B).

Our regular relations with Yu.I. began about 1968, when he
took me as his student at Moscow State University.
My first assignment  
was reading Serre's ``Corps locaux"; I learned Herbrand
theory (the higher ramification), but cannot say the same about
French. 

Let me omit 50 years and go to 2017, the Arbeitstagung devoted
to his 80{\small th} birthday. It was a great meeting! 
Yu.I. and Ksenia Glebovna were terrific hosts, 
many people were around, a perfect view of Rein from 
their apartment etc.

Mostly we discussed anything but mathematics, though
something came up: {\sf\em zeta-polynomials}, certain
combinations of modular form periods satisfying
Riemann Hypothesis, predicted by Manin.

My talk was mostly about DAHA
{\sf\em superpolynomials}
$\h(q,t,a)$ for  {\sf\em double affine Hecke algebras}. 
Superpolynomials  have several
interpretations; the major one is via 
{\sf\em Khovanov-Rozansky
triply graded homology}. 

This direction is in progress.
Conjecturally, topological superpolynomials, 
those from the BPS states (SCFT), DAHA superpolynomials,
motivic ones and $L$-functions of plane curve singularities
coincide (when these theories overlap).
I will focus on the latter three below; this note is introductory, 
with
very few names and 
references.  
\vskip 0.1cm

\subsection{\bf Modular periods}
Let us
begin with Manin's well-known  paper 
{\sf\em ``Periods of parabolic forms and $p$-adic Hecke series"} (1973).
Basically, you consider a parabolic (cusp)
form $\Phi(z)$  of even weight $w$ and calculate its {\sf\em periods}
$r_k(\Phi)=\int_0^{\imath\infty} \Phi(z)z^k dz$ for 
$0\le k \le w-2$. Then the ratios of $r_k$ for even $k$ or
those for odd $k$ are rational numbers, which can be calculated
(Manin's theorem).  

For instance, such $\Phi$ are proportional to 
$\De=e^{2\pi \imath z}\prod_{n=1}^\infty
(1-e^{2\pi\imath n z})^{24} $ for $w=12$. Then $r_2/r_0=
-\frac{2^2 3^4 5}{691}$, $r_3/r_1=-\frac{2^4 3}{5^2}$, etc. 
We note a relation to the Ramanujan's
$\tau(n)\!=\!\si_{11}\!\mod 691$ (1916); Manin reproved it.

The periods are
essentially the values $L_\Phi(s)$ of the corresponding $L$-function
for integer $s$ inside the critical strip. This can be extended to 
$L_\Phi(s,\chi)=\sum_{n=1}^\infty \chi(n)\la_n n^{-s}$
for suitable Dirichlet characters $\chi$ if $\Phi$ is
an eigenfunction of the Hecke operators $T_n$ 
with eigenvalues $\la_n$.

\vskip 0.2cm

The second part of
his paper was on the $p$-adic extrapolations of the ratios of
the periods, which is  closely related to the {\sf\em
Kubota-Leopoldt $p$-adic zeta function} and {\sf\em eigenvarieties}.
Concerning the origins of this direction, let me mention at 
least Barry Mazur and 
Nicholas Katz.

The periods are generally
for any paths $\ga[0,\imath \infty]$ 
for $\ga\in SL(2,\Z)$, but $[0,\imath\infty]$ is sufficient 
due to the modularity of $\Phi$. However, more general paths
do occur in the Manin's paper in process of calculations.
\vskip 0.2cm
 
The $p$-adic extrapolations and the Kubota-Leopoldt zeta (1964)
are closely related to the {\sf\em Iwasawa invariants} of 
$\Ga$-extensions due to Mazur and Wiles in full generality.
The examples of $\Ga$-extensions are some towers
of cyclotomic fields,  where we monitor the class numbers. 
Following Mazur, the Iwasawa invariants are parallel to the
{\sf\em Alexander polynomials}. The covers of the
$S^3\setminus K$ for knots $K$ in their theory are similar
to $\Ga$-extensions; those of $\mathbb P^2$ minus the corresponding
plane curve singularity  are sufficient
for algebraic $K$ due to Libgober
and others. 

The Alexander polynomials are $\h(q\!=\!t,t,a\!=\!-1)$ 
for the DAHA superpolynomials $\h$, and there is a relation 
to $\rho(q,t)$, refined {\sf\em quasi-$\rho$- 
invariants} 
introduced and discussed below (for algebraic knots).

\subsection{\bf Using DAHA}
The modular periods and DAHA are quite different theories,
but there is a clear common denominator:
the action of $SL(2,\Z)$.
The main feature of DAHA is that it  provides  
a universal formalization of
{\sf\em Fourier transforms} and the action of (projective)
$SL(2,\Z)$ in algebra,
harmonic analysis and physics. 

To be more exact, 
DAHA serve the theories with the Fourier
transform and the Gaussian, where the latter
is an eigenfunction of the former. 
The classical Fourier transform, its $q$-counterparts,
the  
Hankel transform and the Verlinde $S,T$-operators are 
basic examples.
DAHA is a universal (flat) deformation
of the Heisenberg and Weyl algebras, so its
role in Fourier analysis is not surprising. 

Moreover, it appeared that DAHA provides 
invariants of iterated torus links. This 
is not very surprising because the {\sf\em Verlinde
algebras} are closely connected with the invariants
of links and  3-folds. In DAHA theory, these algebras become
{\sf\em perfect quotients} of the
polynomial representations. Hopefully,
DAHA and the theory 
of  modular forms and $L$-functions can eventually merge into one, 
but this will require efforts. 
\vskip 0.2cm

Number theory 
already provided some framework for quite a bunch of 
similar directions, which is hardly accidental.  
Let me quote from the Manin's paper:
``... any points of contact with concrete number-theoretical facts, 
whether old or new, 
take on especial significance. They
discipline the imagination, 
and they provide a breathing space and the opportunity to
evaluate the stunning beauty of past discoveries." This
is very much applicable now to the relations between 
number theory and physics (in both directions).

\subsection{\bf Zeta-polynomials}
Next, in his {\sf\em
``Local zeta factors and geometries under Spec $\Z$"} (2014),
Yu.I. conjectured that a certain combination of
$L_\Phi(1),\,\ldots\,, L_\Phi(w-1)$ is a 
{\sf\em zeta-polynomial}\,:
satisfies the functional equation $s\mapsto 1-s$ and the
Riemann Hypothesis.
This was fully confirmed for $w\ge 4$ by Ken Ono, Larry Rolen and
Florian Sprung in their paper
{\sf\em ``Zeta-polynomials for modular form 
periods"} (2016). 

Let
$
M_\Phi(m)\equal\sum_{j=0}^{w-2}
\left(\frac{\sqrt{N}}{2\pi}\right)^{j+1}
\frac{L_\Phi(j+1)}{(w-2-j)!}j^m 
$ for a $\Ga_0(N)$-modular
$\Phi$. Then 
their  zeta-polynomial $Z_\Phi(s)$ 
is a linear combination of 
$M_\Phi(m)$ for $m=0,\ldots, w-2$ with the
coefficients given in terms of 
Stirling polynomials of the 1{\small st} kind
and the Fernando Rodriguez-Villegas transform. 
Manin used the latter too.
There is a relation of zeta-polynomials to the 
Bloch-Kato Conjecture, which is a Galois cohomological
interpretation of the periods, and related advanced
number-theoretical problems.
\vskip 0.2cm

What is important for us is that there is some ``canonical" way
to combine the modular periods 
in a zeta-polynomial, which resembles 
very much the theory of Witten-Reshetikhin-Turaev
invariants of $3$-folds 
invariants and their relations  
to knot invariants. DAHA invariants are ``canonical" combinations
of basic coinvariants in a similar
way.
\vskip 0.2cm

There are various connections of the {\sf\em WRT-invariants} with
modular forms. Let us at least mention {\sf\em 
``Quantum invariants, modular forms, and lattice points II"}
(K.Hikami, 2006). A challenge is to connect 
the inequalities $0\le k \le w-2$ with those in {\sf\em 
Verlinde algebras}, more exactly with the
range of Macdonald polynomials in 
{\sf\em prefect DAHA modules} at roots of
unity, but this is fully open at the moment.
\vskip 0.2cm
 
\section[\sc \hspace{1em}Basic DAHA theory] {\sc Basic
DAHA theory} 
\subsection{\bf Main definitions.}
DAHA, denoted by $\HH$,  were initially introduced to complete the
theory of Knizhnik-Zamolodchikov equations and Quantum Many
Body Problem.  It 
is a universal flat deformation of $\w$,
the Weyl algebra extended by  
the Weyl group $W$; this is for any reducible 
reduced root systems $R$
(in this note).  The
projective $SL(2,\Z)$ due to Steinberg 
acts in $\HH$.  This is actually the braid group $B_3$; the 
 notation will be $\tilde{SL(2,\Z)}.$
\vskip 0.2cm

For $A_1$, $\HH$ is generated by $X^{\pm1},Y^{\pm1},T$ subject to  
group relations $TXTX=1=TY^{-1}TY^{-1}, Y^{-1}X^{-1}YXT^2=q^{-1/2}$
and the quadratic one 
$(T-t^{1/2})(T+t^{-1/2})=1$. The action of
the $\tilde{SL(2,\Z)}$ is:   
\begin{align*}
&\tau_+\!:\!
Y\!\mapsto\! q^{-\frac{1}{4}}XY,\ X\!\mapsto\!X,\ T\!\mapsto\!T,\ 
\tau_-\!:\! X\!\mapsto\! q^{\frac{1}{4}}YX,\  
Y\!\mapsto\!Y,\  T\!\mapsto\!T,
\end{align*}
exactly as for $\w$ (there is no $q$ here). 
The automorphisms  $\tau_{\pm}$,
the generators of  $\tilde{SL(2,\Z)}$,
 are the preimages of {\small
$\begin{pmatrix} 1 & 1\\ 0 & 1\end{pmatrix}$} and 
{\small $\begin{pmatrix} 1 & 0\\ 1 & 1\end{pmatrix}$}, the
standard generators of $SL(2,\Z).$
 The defining relation
of   $\tilde{SL(2,\Z)}$ is simple:
 $\tau_+\tau_-^{-1}\tau_+ =\si=
\tau_-^{-1}\tau_+\tau_-^{-1}$, which formally gives that
 $\si^2$ is central. 
The element $\si^{-1}$ is the (operator) {\sf\em DAHA-Fourier
transform}. 
\vskip 0.2cm

When $t^{1/2}=1$, $T$ becomes the inversion $s$
of $X$ and $Y$, and we arrive at $\w$.
Upon $t=q$, DAHA is closely related to quantum groups.
The case 
$t=q^k$ as $q\to 1$ serves the Harish-Chandra theory
and its $k$-generalization, called Heckman-Opdam theory
(in mathematics), including spherical functions
and Jack polynomials. 
Also, $q\to 0$ is the
$p$-adic limit and  $t\to 0$ is the Kac-Moody limit. 
The action of  $\tilde{SL(2,\Z)}$ generally collapses in the
limits. However it survives in the important limit
to rational DAHA and Hankel transforms, which is
when 
$X=q^x, Y=q^{-y}$, $t=q^k$ and  $q\to 1$ (for $A_1$).

Also, $\tilde{SL(2,\Z)}$ generally
acts in finite-dimensional {\sf\em rigid} $\HH$-modules.
{\sf\em Perfect representations} are such:\, canonical 
irreducible quotients of polynomial representations at
roots of unity $q$ or for any $q$ and {\sf\em singular} $k$. 
The classical {\sf\em Verlinde algebras} are symmetric (spherical)
parts of  perfect representations when 
$t=q$ and $q$ is a root of unity.

\comment{
This commutative
algebra (it is a quotient of a polynomial ring) is $K_0$
of the so-called {\sf\em reduced category} of the quantum group 
of the corresponding type when $q$
is a root of unity. However,
 $\tilde{SL(2,\Z)}$ does not act in the polynomial representation
(which is $K_0(Rep_q)$ for the Lusztig's quantum group).
}

\vskip 0.2cm
We omit the general definition of
$\HH$. For an arbitrary reduced irreducible 
root system $R$ of rank $n$,
it is generated by  pairwise commutative
$X_\la$ for $\la$ in the {\sf\em the weight lattice} $P$, 
pairwise commutative $Y_\la$ and  $T_i$ for
$1\le i\le n$ such that $(T_i-t^{1/2})(T_i+t^{-1/2})=0$,
where $t_{sht}$ and $t_{lng}$ can be generally used when
$\al_i$ is short and long. {\sf\em Spherical DAHA} is
defined as $\p_+ \HH \p_+$ for the $t$-symmetrizer
$\p_+=\sum_{w\in W}t^{l(w)/2}T_w/\sum_{w\in W}t^{l(w)}$.

\subsection{\bf Polynomial representation} 
The key property of DAHA is the {\sf\em PBW theorem}: any element
$H\in \HH$
can be uniquely represented as $H=\sum c_{\la,w,\mu}
X_\la T_w Y_\mu$ for $\la,\mu\in P$ and $w\in W$ (the non-affine
Weyl group).
 Equivalently, there is a faithful 
action of $\HH$ in the {\sf\em polynomial representation}
$\mathscr{X}=\C[X_\la]$: the one induced from
the character $Y_\la\mapsto t^{(\rho,\la)}, T_i\mapsto t_i^{1/2}$.
It is a deformation of the classical {\sf\em Fock
representation} of the corresponding Heisenberg-Weyl algebra.
We will need below the {\sf\em DAHA coinvariant}:
$\{\,H\,\}=\sum c_{\la,w,\mu} t^{(\mu-\la,\rho)+
\text{length}(w)/2}$, which is $H(1)(X_\la\mapsto t^{-(\rho,\la)})$, 
where $H(1)$ is the action of
$H\in \HH$ at $1\in \mathscr{X}.$
\vskip 0.2cm

Technically, the simplest definition of $\HH$ is via
the usage of $T_0$
instead of $Y_\la$. This approach is directly
related to the $T\!\times\! \C^\ast-$
equivariant $K$-theory of affine flag varieties. 
 Then $Y_\la$ are defined as some
products in terms of $T_i$ for 
$i\ge 0$ and their commutativity is some proposition.
The construction of $Y_\la$  
via $\{T_i\}$ is essentially
due to Bernstein-Zelevinsky and Lusztig. 
Finding explicit defining
relations between $X_\la$ and $Y_\mu$ 
is generally involved unless for $A_n$ and in small ranks.

The exact connection with affine flag varieties was 
fully clarified by
Garland-Grojnowski and Ginzburg-Kapranov-Vasserot (1995).
We note that  the action of $\tilde{SL(2,\Z)}$
is far from obvious from the $K$-theoretical viewpoint
in spite of various attempts in this direction.
 The key property of any Fourier transforms
is that they send polynomials to delta-functions, which is
not simple to incorporate geometrically.

\vfil
Generally, there are
3 fundamentally different definitions (interpretations) of DAHA;
more than 3 for $A_n$. The 2{\small nd} definition  is via the 
{\sf\em orbifold}\, fundamental
group of the so-called {\sf\em elliptic configuration space}:\,
$\pi_1^{orb}\Bigl(\bigl\{x\in E^n \mid
\prod_\al (x,\al)\neq 0\}\bigr\}/W\Bigr)$. Here $E$ is an 
elliptic curve and $W$ is the non-affine Weyl group acting
in $E^n$. This action is well-defined because
the roots $\al\in R$ are with integer coefficients
in terms of the fundamental weights, which are considered as 
coordinates of $E^n$.  

Then we take group algebra of this $\pi_1^{orb}$ and 
impose the quadratic relations for $T_i$ as above, which 
elements are half-turns corresponding to simple reflections.
The element $\si$ becomes basically
the transposition of the periods of $E$.  
The elliptic configuration space  
is the ``big cell" in $Bun_G(E)$ for the corresponding
simple Lie group $G$. 

The existence
of the action of $\tilde{SL(2,\Z)}$
is straightforward  from the topological definition. However, 
the polynomial representation is far from immediate in this
interpretation, 
which is almost by construction via the $K$-theory of affine
flag varieties. I constantly use both approaches and the one
via harmonic analysis, the 3{\small rd} major approach. 

\vskip 0.2cm
\vfil
The action of $\tilde{SL(2,\Z)}$ becomes  ``natural"
if DAHA is interpreted via the $q$-version 
of harmonic analysis on symmetric spaces. Namely,
$\si$ becomes the DAHA-Fourier transform, and $\tau_+$ is
associated with 
the multiplication by the Gaussian,  $q^{x^2/2}$. 
Basically, this projective action
is equivalent to the fundamental theorem of Fourier
analysis: the Gaussian 
is an eigenfunction of the corresponding Fourier transform.

{\sf\em Spherical DAHA} for $GL_n$ can be interpreted as 
{\sf\em elliptic Hall algebras} and are related to quantum groups
due to Schiffmann -Vasserot and others. In spite of the usage of
elliptic curves, the action of  $\tilde{SL(2,\Z)}$ 
 is not  ``immediate" in this approach; further work is needed.
Also, it can be interpreted via
certain {\sf\em shuffle algebras}.

\vskip 0.2cm

Given a reduced irreducible root system $R$,
the {\sf\em nonsymmetric Macdonald polynomials} 
$E_\la$ for $\la\in P$ generalize 
the monomials $X_\la$ for the corresponding Weyl algebra and form
a basis of $\mathscr{X}$. They are eigenfunctions of $Y_\mu$
normalized by the conditions $E_\la=X_\la+(\text{lower terms})$.
In this approach, {\sf\em Macdonald polynomials} 
$P_\la$ are the $t$-symmetrizations of $E_\la$ for $\la\in P_+$;
the normalization is $P_\la=X_\la+(\text{lower terms})$.

The action of $\HH$ in $\mathscr{X}$ is sufficiently 
explicit. For $A_1$:
$T\mapsto t^{1/2}s\!+\! 
\frac{t^{1/2}\!- t^{-1/2}}{X^2-1}(s\!-\!1)$, 
$X\mapsto X,\ \, Y\mapsto spT$, where $ s(X)=X^{-1},\,
\, p(X)=q^{1/2}X.$ The divided differences here and for
any root systems $R$ are very standard
in the theory of affine Hecke algebras and are quite common in
related geometry and combinatorics. 

\vskip 0.2cm
For $GL_n$, the corresponding
$\HH$ is generated by $X_i^{\pm 1}, Y_j^{\pm 1}, T_k$,  where 
$1\le i,j\le n$ and
$1\le k \le n-1$ for pairwise commutative $\{X_i\}$ and 
$\{Y_j\}$. One has:   
 $\tau_+(Y_1)\!=\!q^{-1/2} 
X_1 Y_1$,\  $\tau_-(X_1)\!=\!q^{+1/2}
Y_1 X_1$ and so on. The action of $Y_1$ in $\mathscr{X}$ is via
the formula 
$Y_1\!=\!\pi T_{n-1}\ldots T_1$, where 
$\pi: X_1\!\mapsto\! X_2, X_2\!\mapsto\! X_3,
\ldots, X_n\!\mapsto\! q^{-1} X_1$. These formulas are quite
similar for any $Y_i$.
We will use them below when calculating
the DAHA-superpolynomial of trefoil (as an example).
\vskip 0.2cm
\vfil

Going back to the modular periods,
the evaluation map $X\mapsto t^{-\rho}$, {\sf\em
the DAHA coinvariant},
plays
a role of integration $\int_0^{\imath\infty}
\{\cdot\} \Phi dz$,
$E_\la$ replace $z^k$, and the action of $\tga$ corresponds
to the change of the integration path to  $\ga[0,\imath\infty]$.
The main deviation is that the action of
$\ga$ plays a much more significant
role for  DAHA superpolynomials versus that for the periods.
In contrast to $\int_0^{\imath\infty}
\{\cdot\} \Phi dz$, 
the DAHA coinvariant is not stable in any way with respect to
the action of $\tilde{SL(2,Z)}.$ 
\vskip 0.2cm

When
switching to the zeta-polynomials $Z_\Phi(s)$, 
special linear combinations of $z^k$-momenta must
be considered, which resembles our usage of 
$E_\la$. Both constructions are
``canonical" in a sense; the restriction 
$0\le k\le w-2$ seems  somewhat similar to those in Verlinde algebras.
\vskip 0.2cm

The next topic, {\sf\em DAHA-Verlinde algebras},
gives a direct link of DAHA at roots of
unity to number theory. They are some counterparts  
of {\sf\em Tate modules} for the towers of
covers of an elliptic curve ramified at one point; the
{\sf\em absolute Galois group} acts there. 
The Verlinde
algebras are one of the key ingredients of the 
{\sf\em Witten-Reshetikhin-Turaev
invariants}, generalize $K_0$ of the 
{\sf\em reduced category} in 
representation theory of 
{\sf\em quantum groups} at roots
of unity, and that of
integrable modules of {\sf\em Kac-Moody algebras}
(Kazhdan-Lusztig, Finkelberg). 

\subsection{\bf Refined Verlinde algebras}
These algebras in the nonsymmetric variant
are {\sf\em perfect} finite-dimensional
quotients
of 
$\mathscr{X}$, those with the action of $\tilde{SL(2,\Z)}$
and  DAHA-invariant non-degenerate quadratic forms.
They are commutative algebras, but
can be non-semisimple, those related to {\sf\em
logarithmic\,} CFT in examples. 
Technically, we divide $\mathscr{X}$ by the radical of the
evaluation pairing. Such modules exist  either when $q$ is
a root
of unity or for {\sf\em singular} $k$, where $t=q^k$ as above. 

In the case of $A_1$, let $q=\exp(\frac{2\pi i}{N}),
k\in \frac{\Z_+}{2}$ and 
$k\!<\! N/2$. 
The map $X(z)=q^{z}$ can be naturally
extended to an $\HH$-homomorphism
$\C[X^{\pm 1}]\to
V_{2N-4k}$. The latter is the {\sf\em nonsymmetric Verlinde algebra}.
It the space of functions 
{\small
$
f:\{-\frac{N+k+1}{ 2},...,-\frac{k+1}{ 2},-\frac{k}{ 2},\, 
\frac{k+1}{ 2},...,\frac{N-k}{ 2}\}\to \C$
}
with pointwise multiplication and 
the action of
$X,T,Y$ induced from that in $\mathscr{X}$.
These modules are {\sf\em rigid},
which readily gives an action of  
$\tilde{PSL_2(\Z)}$ there,
and in $V_{N-2k+1}^{sym}\!\equal\!
\{v\in V\!\mid\! Tv\!=\!t^{\frac{1}{2}}v\}$.
The low index is the dimension of the corresponding $V$:
 dim$\,V_{2N-4k}\!=\!2N\!-\!4k$, and 
dim$\,(V_{2N-4k})^{sym}$ $=N\!-\!2k\!+\!1$.
The operators 
$X,Y,T$ become unitary in $V_{2N-4k}$ if the 
$q=e^{\frac{2\pi \imath}{N}}$, the
``minimal" primitive $N${\tiny th} root.
\vfil

We classified {\sf\em rigid} modules for $A_1$ 
in {\sf\em ``On Galois action in rigid DAHA modules"} (2017).
They are: ($\boldsymbol{\al}$) $V_{2N-4k}$ as above, 
($\boldsymbol{\be}$) 
non-semisimple $V_{2N+4|k|}$ for
$k\in -\Z_+$ such that $-N/2<k<0$, and ($\boldsymbol{\ga}$) $V_{2|k|}$ for
$k=-\frac{1}{2}-m>-N/2$, where $ m\in \Z_+$. There is
a similar list for the 
{\sf\em little} DAHA 
$\HH'=\lan X^{\pm 2},Y^{\pm 2},T\ran\subset \HH.$ Importantly,
families $(\boldsymbol\al,\boldsymbol\ga)$ have {\sf\em flat} 
$q$-deformations from roots of unity to arbitrary $q$.
The unimodular
$q$ such that  \,arg $q\le \frac{2\pi}{N}$ result in the 
positivity of the invariant form in type $(\boldsymbol\al)$. 
Such a deformation leads to 
some relations between $V,V'$ defined for $N\mid N'$,
generalizing those for
the {\sf\em Tate modules}.
\vskip 0.2cm

The usual {\sf\em Verlinde
algebra} is $V^{sym}_{N-1}$ of type ($\boldsymbol\al$),
 which is for $k\!=\!1$, i.e. for $t=q$.
 Then 
$\tau_+$ becomes the $T$-operator, and
$\si=\tau_+\tau_-^{-1}\tau_+\ $ becomes 
the Verlinde $S$-operator.
The {\sf\em reduced characters} in Verlinde algebras, generally
become the images of the
corresponding Macdonald $E$-polynomials for $V$ and symmetric
ones for $V^{sym}$.

Perfect representations are quotients of the ones obtained
from  $\mathscr{X}$ by fixing the corresponding central
characters, which are of dimension $4N$. For $k=1$, 
the symmetrizations $V^{sym}$ of the latter 
are connected with the category of representations
of {\sf\em small quantum group}. For instance, 
$V^{sym}_{N-1}$  is the 
Grothendieck ring $K_0$
of the so-called {\sf\em reduced category} for $A_1$. The 
perfect representations for  $\Z/2\ni k\neq 1$ are generally beyond 
quantum groups, though the ones of
type $(\boldsymbol\be)$ 
are connected with {\sf\em logarithmic} conformal
field theories and there are other links.

\subsection{\bf The Galois action}
The rigidity provides that the {\sf\em absolute Galois group}
acts in the modules above (including
the usual Verlinde algebras). We use that {\sf\em elliptic
braid group}  $\b_q$
generated by $X,T,Y$ subject to the group relations in the
definition of $\HH$ of type $A_1$ is a renormalization of
the  orbifold fundamental group 
$\pi^{orb}_1(E/\{1,s\})$, where $E$ is an elliptic curve,
$s:x\mapsto -x$.
If $E$ and its origin $o$ 
are defined over some field 
$\Q[q^{1/4}]\subset K\subset \overline{\Q}$, 
then Gal$(\overline{\Q}/K)$ acts projectively
in these modules. 
\vskip 0.2cm

More exactly, 
setting $A=XT,\, B=XTY,\,  C=T^{-1}Y$,
the relations of $\b_q$ and the action of $\tau_{\pm}$ there 
become as follows:
{\small
\begin{align*}
A^2\,=\,1&\,=\,C^2=q^{1/2}B^2,\, \hbox{\, where \,}
ABC=A^2 YT^{-1}Y=
YY^{-1}T=T,\notag\\
&\tau_+:\ A\mapsto A,\ B\mapsto q^{-1/4}C,\ C\mapsto
q^{1/4}C^{-1}BC,
\\
&\tau_-:\  A\mapsto q^{1/4}ABA^{-1}, \
B\mapsto q^{-1/4}A,\  C\mapsto C. \notag
\end{align*}
}
\!\!The classification of $\HH$-modules at roots
of unity $q,t$ becomes the corresponding
multiplicative {\sf\em Deligne-Simpson problem} with specific 
quadratic relations for $A,B,C,D=T^{-1}$. They can be arbitrary quadratic
for DAHA of type $C^\vee C_1$ (Sahi, Noumi-Stokman);
let me also mention Oblomkov-Stoica (2009).
This algebra is generated
by $A,B,C,D$ such that  $ABCD=1$ satisfying any quadratic relations:  
those for  the monodromy of 
the {\sf\em Heun equation.} There are links to SCFT.

\vskip 0.2cm
The images of $\b_q$ in type $(\boldsymbol\al)$ rigid modules with
positive-definite invariant forms  are
finite and we obtain finite covers of $\mathbb P^1$
ramified at $0,1,\infty$ and  $o\in E(K)$,
where $A,B,C,D$ are the corresponding
monodromies. 
When $t=1$, we arrive at 
unramified covers of $E$.

\vskip 0.2cm

The case of the Hermitian 
invariant forms {\sf\em with one minus}
is interesting. Then the images
of $\b_q$ are discrete groups. The smallest
nontrivial such $V$ is for {\sf\em little} $\HH$;
dim$V=3$. In particular, we obtain then {\sf\em all}  
{\sf\em Livn\'e lattices} in $PU(2,1)$,
which are examples of
the {\sf\em Mostow groups}. Livn\'e used
a branched 2-cover of degree $2$ of
the {\sf\em universal elliptic curve}.

More generally (for the same $V$),
there  is a direct connection with 
the theory of {\sf\em equilateral triangle groups} in
$PU(2,1)$; for instance, see 
{\sf\em ``Complex hyperbolic triangle groups"} (R.E. Schwartz, 2002)
and {\sf\em ``Cone metrics on the sphere and Livn\'e's lattices"}
(Parker, 2006). 
\vskip 0.2cm

We mention here that the (regular) {\sf\em Inverse Galois Problem}
is based on {\sf\em rigid triples}, which are $\{a,b,c\}$
generating a group $G$ and
satisfying $abc=1$. They are assumed from
given conjugacy classes in $G$ and 
the rigidity means essentially the
uniqueness of such $\{a,b,c\}$  up to (simultaneous)
conjugations in $G$. 
We need $\{a,b,c,d\}$ here, $4$ points in $\mathbb P^1$ and 
the {\sf\em linear rigidity\,} (in matrices) based on
Katz' theory of rigid systems 
(in the variant due to M. Dettweiler and others). 

Such covers extend  
the Belyi's theorem and
Grothendieck's program of {\sf\em dessins d'enfants\,} to $E$;
let us mention  Beilinson-Levin (1991). We 
deal only with very ``small" covers: those from DAHA
modules. Our towers are  similar to those from Tate modules,
though our (ramified) version of
 $T_p(E)=\varprojlim (E/E_{p^n})$ 
is {\sf\em not} a module over $p$-adic numbers. Here and above
see {\sf\em ``On Galois action in rigid DAHA modules"}.

\section[\sc \hspace{1em}Knot invariants via DAHA]
{\sc Knot invariants via DAHA}
\subsection{\bf DAHA-Jones polynomials}
 They can be defined for any reduced root system
$R$; in the non-reduced case, see author's
{\sf\em ``Jones polynomials of torus knots via DAHA"} and
{\sf\em ``DAHA-Jones polynomials of torus knots"}. The system
$C^\vee C_1$ is discussed in the latter.
Algebraic torus knots
are for $r,s>0$ such that gcd$(r,s)=1$. They can be
represented as
$T(r,s)=\{x^r=y^s\}\cap S_\ep^3$ for a small sphere $S_\ep^3$
centered at $0$. 
The formula for the corresponding {\sf\em DAHA-Jones invariant}
is:
$$
J^\la_{r,s}(q,t)=\bigl\{\tilde{\ga}\bigl(
\frac{E_\la}{E_\la(t^{-\rho})}\bigr)(1)\bigr\},\, 
\bigl\{F(X)\bigr\}=F(X\mapsto t^{-\rho}),
\, \la\in P_+.
$$  
where  $(r,s)^{tr}$
is the 1{\small st} column of $\ga\in SL(2,\Z)$, 
$\tilde{\ga}$
is its action in $\HH$, and the Laurent polynomial
$\tga(E_\la)(1)$ is 
$\tga(E_\la)\in \HH$ applied to $1\in \mathscr{X}$. 
It is not necessary
to assume that $r,s>0$ in this definition; 
the corresponding torus
knots will be non-algebraic for $rs<0$.

A remarkable theorem is that it is always a $q,t$-polynomial 
(up to some fractional power $q^\bullet t^\bullet$)
in spite of the $q,t$-singularities
of $E/E_\la(t^{-\rho})$, which exist only for some 
$\la$  at roots of unity $q$  and
for singular $t$.  The product formula for 
 $E_\la(t^{-\rho})$ is the Macdonald evaluation conjecture
in the nonsymmetric variant (now a relatively simple
theorem).

The formulas for  $J^\la_{r,s}(q,t)$ and their generalizations
below  for
iterated torus knots will give exactly the same invariants
if $E_\la$ is replaced by the 
{\sf\em symmetric} Macdonald polynomials $P_\la$ for $\la\in P_+$. 
Note that 
$P_\la(t^{-\rho})=P_\la(t^{\rho})$. We 
employ
the $t$-symmetrization inside the coinvariant $\{\cdot\}$. 

The $P$-polynomials become the Schur-Weyl characters for $R$
when $t=q$; they do
not depend on $q,t$ in this case; then $P_\la(t^{\rho})$ is
the corresponding $q$-dimension. Colored {\sf\em Jones polynomials} 
for $T(r,s)$ are $J^{\la}$ for $A_1$, $\la=m\om_1$
and $q=t$, which is 
up to some power $t^\bullet$.

\vskip 0.1cm

The usage of the $E$-polynomials is important. 
For instance, 
$E_{\om_i}=X_{\om_i}$ in the $A_n$-case, and 
$J^{\om_i}_{r,s}(q,t)$ are closely related to the left-hand side
of the {\sf\em Shuffle Conjecture}, $\nabla e_n[X]$, 
 proved by Carlsson-Mellit in their 
{\sf\em ``A proof of the shuffle conjecture"} (JAMS, 2018).
The right-hand side of this conjecture (in the most general
setting) is the corresponding {\sf\em
motivic superpolynomial} (below).
In the uncolored case $\la=\om_1$; {\sf\em we will 
omit $\la$ in $J^{\la}$ and other formulas in this case.} 
\vfil

\vskip 0.2cm

{\sf\em $A_n$-stabilization.} 
For any $G=1+q\C[[q,t]]+t\C[[q,t]]\in C[q,t]$ and rational $u,v$, 
let 
 $ \bigl(q^u t^v G\bigr)^\circ\equal G.$ 
If $G$ also depends on $a^{\pm 1}$ (below), then
$ \bigl(q^\bullet t^\bullet a^\bullet\, G\bigr)^\circ \in 
1+q\C[q,t]+t\C[q,t]+a\C[q^{\pm 1},t^{\pm 1}]$.
\vskip 0.2cm

Given 
a Young diagram $\la= (\la_1\ge \la_2\ge \cdots \ge \la_m>0)$, 
$\la=\sum_{i=1}^m c_i \om_i$ is considered a weight for 
$A_{n}$, where $n\!\ge\! m\!-\!1$ and $c_i$ is the number
of columns of size  $i$ in $\la$.

The claim is that given $r,s$ as above,
there exists a unique polynomial
$\h_{r,s}(q,t,a)\in \C[q,t^{\pm1},a]$     
such that $J_{r,s}^\la(q,t)^\circ=\h_{r,s}^\la(q,t,a=-t^{n+1})$ for 
$J_{r,s}^\la$ of type $A_n$ for $n\!\ge\! m\!-\!1$    
Automatically, $\h_{r,s}(q,t,a)^\circ=\h_{r,s}(q,t,a).$

\vskip 0.2cm

{\sf\em Comments.} The starting point of this
theory was due to Aganagic-Shakirov (2011) and the author (2011). 
Concerning the related physics, let me mention
at least the paper by Gukov, Iqbal, Kozcaz and Vafa (2010).
The stabilization of $J^\la$ for $A_n$
is based on a DAHA theorem due to 
Schiffmann-Vasserot (2012). Let me mention here 
the proof of the DAHA-superduality for torus knots
by Gorsky-Negut (2013).

It was conjectured by the author  that the $a$-stabilization
holds for $B,C,D$ too; topologically, we generalize
{\sf\em Kauffman polynomials}. Moreover, 
superpolynomials were calculated for a couple of knots 
by the author and R.Elliot (2016) for the
exceptional series from (Deligne-Gross, 2002):
$
A_0\subset A_1\subset A_2\subset D_2\subset F_4\subset D_4\subset
F_4\subset E_6\subset E_7\subset E_8.$
Our paper was mostly on the superpolynomials in the case of the
{\sf\em annulus}. 
\vskip 0.2cm

Let us state a version of the stabilization conjecture for $C_n$
(from my paper).
One has $R_+=\{\ep_i\pm \ep_j, 2\ep_i\}$ for
$1\le i,j\le n$ such that $i\le j$, 
$P_+=\{\la=\sum_{i=1}^n\la_i \ep_i, \la_1\ge \la_2\ge \cdots\ge 
\la_n\ge 0\}$.  Assuming that
$\la_m>0$, they become dominant weights for $C_n$ for any $n\ge m$. 
We will treat $\la$ as Young diagrams.
DAHA invariants depend now on $q$, 
$t=t_{\sht}=q^{k_{sht}}$ and the additional
parameter $u=t_{lng}=q^{k_{lng}}$. We set 
$\rho_k=\frac{1}{2}\sum_{\al>0}k_\al \al$ and replace
$t^\rho\mapsto q^{\rho_k}$ in all formulas. 
\vskip 0.2cm
 
Given a torus knot $T_{r,s}$ and $\la=(\la_i)$ with $\la_m>0$,
conjecturally there exists a polynomial
$\h^{C}_{r,s}(\la\,;\,q,t,u,a)$ with integral
coefficients  in terms of positive powers of 
$a,q,u$ and $t^{\pm 1}$ such that for any $n\ge m$
\begin{align*}
&\h^{C}_{r,s}(\la\,;\,q,t,u,a=-t^{n-1})=
J^{C_n}_{r,s}(\la\,;\,q,t,u)^\circ.
\end{align*}
The uncolored invariants are for $\om_1=\ep_1$, which is
minuscule for $C_n$.

\subsection{\bf The case of trefoil}\label{sec:tref}
Let us calculate $\h_{3,2}$ for uncolored trefoil; we
begin with $A_1$. As above,
$\{H\}\!\equal\! H(1)(X\!\mapsto\! t^{-\rho})$,
where $t^{-\rho}=t^{-\frac{1}{2}}$ for $A_1$.
By $\,\sim$\,, 
we mean ``\,up to $q^{\bullet}t^{\bullet}$\,". One has:
\begin{align*}
&J_{3,2}\!=\!
\{\tau_+\tau_-^2(X)\}\!\sim\!\{(XY)(XY)X(1)\}\!\sim\!
\{Y(X^2)\}\\
&=t^{-\frac{1}{2}}q^{-1}X^2-
t^{\frac{1}{2}}+t^{-\frac{1}{2}}|_{X^2\mapsto t^{-1}}\sim
1+qt-qt^2.
\end{align*}
When $q=t$, we obtain the {\sf\em Jones polynomial}:
$J_{3,2}(q\mapsto t)^\circ= 1\!+\!t^2\!-\!t^3.$

We use that $E_1\!=\!X$:\,
$Y(X)\!=\!(qt)^{-\frac{1}{2}}X$. 
Using the formula for $Y_1$ for $A_n$ above
and the
action of $\tau_{\pm}$ on $X_1,Y_1$, 
we obtain that  $J_{3,2}^\circ
=1+qt-qt^{n+1}$, 
which gives that $\h_{3,2}=1+qt+aq$ for $a=-t^{n+1}$.
The relations
$\h_{3,2}(a\!\mapsto\! -t)=1$ and    
$\h_{3,2}(a\!\mapsto\! -t^{2})=1+qt-qt^2$
are sufficient to fix it uniquely if it is known
that $deg_a\h=1$. Generally, 
$deg_a\h^{\la}_{r,s}\!=\!
|\la|\bigl(\text{Min}(r,s)\!-\!1\bigr)$. 
A remarkable simplicity of $\h_{3,2}$ is fully clarified 
in the approach via motivic superpolynomials (below).
\vskip 0.2cm

The case of $T(2p+1,2)$ is quite similar.  For $p=1,2,\ldots\,$,
one has:
$\h_{2p+1,2}=
1 + qt + q^2 t^2 +\cdots+ q^p t^p+  aq(1 + qt + \cdots+(qt)^{p-1}).$
Similarly, the Khovanov-Rozansky polynomials
are the simplest for these knots. 

\Yboxdim6pt
The formula  becomes significantly
more involved with colors. Let $\la=m\om_1=$ 
$\yng(2)\cdots\yng(1)$ ($m$ boxes) for
$m=1,2,\ldots\ $. Then:
\begin{align*}
&\h^\la_{2p+1,2} =
\frac{(q;q)_m}{(-a;q)_m(1-t)}\!\sum_{k=0}^m (-1)^{m-k}\\
&(qt)^{\frac{m-k}{2}}
\bigl((q^{\frac{m(m+1)}{2}}-q^{\frac{k(k+1)}{2}})
(t/q)^{\frac{m-k}{2}}\bigr)^{2p+1}\\
&\frac{ (t;q)_k(-a;q)_{m+k}(-a/t;q)_{m-k}(1-q^{2k}t) }
{(q;q)_k(qt;q)_{m+k}(q;q)_{m-k}},
\end{align*}
where $(a;q)_n=(1-a)\cdots (1-a q^{n-1})$.\
This formula was proposed by {
Dunin-Barkowski, Mironov, Morozov, Sleptsov, Smirnov} 
(2011-12), and, independently, by 
{Fuji, Gukov, Sulkowsky} (2012). 
A somewhat different formula is $\h^\la_{3,2}\! =\!$
$\sum_{k=0}^m q^{mk}t^k ${\Large $\frac{(q;q)_m(-a/t;q)_k}
{(q;q)_k(q;q)_{m-k}}$} (only for trefoil). 
The justifications were obtained via DAHA, i.e. for
the DAHA superpolynomials.
The Habiro's formula (2000)   
is for $p=1, a=-t^2, t=q$. 
\vskip 0.2cm

Let us mention here  
{\sf\em ``Torus knots and quantum modular forms"} devoted
to color Jones polynomials for $T(2p+1,2)$ 
(K.Hikami-Lovejoy, 2014), and the {\sf\em Kontsevich-Zagier series}
from
{\sf\em 
``Vassiliev invariants and a strange identity related to the 
Dedekind eta-function"} (Zagier, 2001). 
Presumably, our refined formulas above can be used
in a similar way. See also Example 5 from 
{\sf\em ``Quantum modular forms"} (Zagier, 2010).

In the case of uncolored trefoil for the system $C_n$, 
$\h^{C}_{3,2}(\om_1\,;\,q,t,u,a)=$
$1+q t+a (q t-q u)+a^2 \bigl(-q u+q^2 u-q^2 t u\bigr)
+a^3 \bigl(-q^2 t u+q^2 u^2\bigr).$

The superpolynomial for the {\sf\em exceptional series}
for $T(3,2)$ and for the adjoint representation is
due to Cherednik-Elliot (2016): 

\renewcommand{\baselinestretch}{0.5} 
{\small
\noindent
\( 
\h^{E}_{3,2}(adj;\,q,t,a) \; = \;  1 + q(t - ta + a^2 - a^4 + 
t^{-1}a^5 - t^{-1}a^6) + q^2(t^2a^2 - ta^3 + a^4 + ta^5 + 
t^{-1}a^6 - 3a^6 + t^{-1}a^7 + a^7 - t^{-1}a^8 - t^{-1}a^9 + 
t^{-1}a^{10} - t^{-2}a^{11}) 
 + q^3(t^{-1}a^6 - a^7 + ta^7 + t^{-1}a^8 - a^8 - ta^8 
+ a^9 - 2t^{-1}a^{10} + a^{10}
+ t^{-2}a^{11} - t^{-1}a^{11} - a^{11} - t^{-2}a^{12} + 
2t^{-1}a^{12} - t^{-2}a^{13} + t^{-2}a^{15}) 
 + q^4(t^{-2}a^{12} - t^{-1}a^{12} + t^{-1}a^{13} - a^{13} 
- t^{-2}a^{14} + t^{-1}a^{14} + t^{-2}a^{16}
 - t^{-1}a^{16} - t^{-3}a^{17} + t^{-2}a^{17})
 + q^5(-t^{-3}a^{18} + t^{-2}a^{18}).
\)
}
\renewcommand{\baselinestretch}{1.2} 

\noindent
Here $a=-t^{\frac{h}{6}}$ for the Coxeter number $h$ for the 
$A,D,E$ there.

\subsection{\bf Iterated links} We will begin
with the DAHA construction for iterated torus knots. 
For any sequence 
$\ga_1,\ga_2,\ldots,\ga_\ell \in SL(2,\Z)$,
we set 
$J^\la=\biggl(\cdots \tga_{\ell-1}
\Bigl(\tga_\ell\bigl(\frac{E_\la}{E_\la(t^{-\rho})}\bigr)(1)
\Bigr)(1)\cdots 
\biggr)(t^{-\rho})$. This
is due to Cherednik-Danilenko; the invariant $J^\la$ 
depends only on the isotopy type of the 
corresponding {\sf\em 
iterated torus knot}; see an example below. 

This is somewhat similar to Manin's work
{\sf\em ``Iterated integrals of modular forms and noncommutative 
modular symbols"} (2005). Basically, $\int_{0}^\infty$
is replaced by $\int_p^{q}$ for rational $p,q$ in this
paper and
multiple zeta values occur. 
When the coinvariant is replaced by the
corresponding integral formula (a DAHA theorem), the similarity
becomes less speculative. 
\vskip 0.2cm

For instance,  one obtains for  $K= 
C\!ab(53,2)C\!ab(13,2)C\!ab(2,3)$:
{
$$
J_{\l}^\la\!=\!\bigl\{\p_\la \bigr\},\   \p_\la\!=\,   
\Downarrow\!
\begin{pmatrix} 3 & \ast\\ 2 &\ast\\ \end{pmatrix}
\Downarrow\!
\begin{pmatrix} 2 & \ast\\ 1 &\ast\\ \end{pmatrix}
\Downarrow\!
\begin{pmatrix} 2 & \ast\\ 1 &\ast\\\end{pmatrix}
\!\bigl(\frac{E_\la(X)}{E_\la(t^{-\rho})}\bigr),
$$
}
\noindent
\!\!where the $\ga$-matrices act via their lifts to Aut$(\HH)$, 
$\Downarrow\!\!H\!\equal\!H(1)$,  
$\bigl\{H\bigr\}${\!\footnotesize $\equal$\!}
$H(1)(t^{-\rho})$ is the {\sf\em coinvariant}, and
$E_\la$ is the 
$E$-polynomial for dominant $\la$.
Here
$C\!ab(a,b)K$ is $T(b,a)$ plotted at the boundary of the solid torus
around a given knot $K$  for the Seifert zero-framing.

Generally, given a sequence {\footnotesize
$\begin{pmatrix} r_1 & r_2 & r_3 & \cdots
\\ s_1 & s_2 & s_3 & \cdots \end{pmatrix}
$} of the 1{\small st} columns of $\ga_1,\ga_2,\ga_3, \ldots$,
which is {\footnotesize
$\begin{pmatrix} 3 & 2 & 2\\ 2 & 1 & 1 \end{pmatrix}
$} in the example above, the corresponding  cable is
 $\cdots C\!ab(a_3,r_3)Cab(a_2,r_2)C\!ab(a_1,r_1)$ 
for $a_1\!=\!s_1$, $a_2\!=\!r_1s_1r_2+s_2,$\, 
$a_3\!=\!a_2 r_2 r_3+s_3\,$, and so on.
Note that 
$C\!ab(s_1,r_1)=T(r_1,s_1)=C\!ab(r_1,s_1)$;
the transposition of
 $r_i,a_i$ for $i>1$ changes the knot. 

The $a$-stabilization theorem  is 
the same as for torus knots. 
Generally, deg$_a(\h^{\la})=|\la|(mult-1)$
for  
 $mult=\text{Min}(r_1,s_1)\cdot r_2\cdot r_3\cdots$\,,
which is the multiplicity  of singularity.
When $a\!=\!0$, deg$_{\,t}(\h^\la)=\de\sum_{i} m_i^2$ for 
$\la=(m_1\!\ge\! m_2\!\ge \cdots)$, and deg$_{\,q}(\h^\la)=
\de\sum_{i} (m'_i)^2$ for the transposition 
$\la'$ of $\la$; $\de$ is the arithmetic genus of the
singularity (below). The justification
is under minor assumptions; the superduality and the specialization 
$q\!=\!1$ provide the reduction to pure columns. 
\vskip 0.2cm

{\sf\em Topological invariance.} 
For torus knots, the isotopy invariance means that 
$T(r,s)$,  $T(s,r)$ and $T(-s,-r)$ must have coinciding 
 $\h^\la$, and
that $\h^\la=1$ for any $T(1,s)$. Also, the superpolynomial of
$T(-s,r)$, which knot is the mirror image of
$T(s,r)$, must be the ``conjugation" 
$q,t,a\mapsto q^{-1},t^{-1},a^{-1}$ of $\h^\la$ for $T(s,r)$ 
up to a factor  $q^\bullet t^\bullet a^\bullet$. This 
requires the usage of the DAHA-automorphism $\eta$, 
which we omit.

We note that the symmetry 
 $\h^\la_{r,s}=\h^\la_{s,r}$  can be 
a challenge for (other) algebraic and
algebraic-geometric approaches even in the case of 
HOMFLY-PT polynomials. It is a simple DAHA lemma for us. 
 
Generally, the theorem is that DAHA superpolynomials depend 
only the isotopy type of the corresponding 
iterated torus link. It is interesting, all fundamental properties of the
{\sf\em coinvariant} are needed, but not a difficult one. 
The case of iterated
{\sf\em links} is more involved, but this is
mostly because {\sf\em splice diagrams} are used,
which are not too simple. The 
additional DAHA fact is 
the integral formula for action of $\tau^{-1}$, which
we will provide below when discussing {\sf\em DAHA vertex}.
\vskip 0.2cm

{\sf\em From knots to links}.
The construction becomes more ramified for iterated torus
{\sf\em
links}. We will provide here the procedure from
in {\sf\em ``DAHA approach to 
iterated torus links"} (Cherednik-Danilenko, 2015).

\vskip 0.2cm

First, we switch from $P_\la$ to the so-called 
$J$-polynomials. For $A_n$: 
\begin{align*}\label{P-arms-legs}
\tilde{P}_\la\equal h_\la P_\la \for 
h_\la=\prod_{\Box\in\la}
(1-q^{arm(\Box)}t^{leg(\Box)+1}),
\end{align*}
for the Macdonald polynomials $P_\la$. 
Here $arm(\Box)$ is the number of boxes in the same row as $\Box$
strictly after it; $leg(\Box)$ is
the number of boxes in the column of $\Box$
strictly below it. The $J$-polynomials and the construction
below are for any root systems $R$.

We now set
$\tilde{\p}^\la(X)=\biggl(\cdots \tga_{\ell-1}
\Bigl(\tga_\ell\bigl(\tilde{P}_\la\bigr)(1)
\Bigr)(1)\cdots 
\biggr)$, called {\sf\em basic prepolynomials}.
General {\sf\em prepolynomials} are defined inductively as follows.
Given two prepolynomials $\tilde{\p}_1,
\tilde{\p}_2$ and $\ga\in PSL(2,\Z)$ (can be $id$), we
define a new {\sf\em prepolynomial} 
$\biggl(\tga\bigl(\tilde{\p}_1\tilde{\p}_2\bigr)\biggr)(1)$.

Combinatorially, we obtain a union of trees with 
marked {\sf\em last
vertices} and the corresponding ends ``colored" by the diagrams
$\{\la\}=\{\la^1,\cdots,
\la^{\kappa}\}$, where  $\kappa$ is the number of connected
components of a link.
The prepolynomials for the (maximal) roads in this union
are those for the connected components. A union of several trees
is for a disconnected union of the 
corresponding links. 
\vskip 0.2cm

This is not the end. Given two prepolynomials,
$\tilde{\p}(X)$ and $\tilde{\q}(X)$,
 for the components colored by the sequences $\{\la^i\}$
and  $\{\mu^j\}$ of diagrams,
$$
J^{\{\la\},\{\mu\}}\equal
\bigl\{\tilde{\q}(Y) \tilde{\p}(X)\bigr\}/
LCM\bigr(\tilde{P}_{\la^i}(t^\rho), \tilde{\q}_{\mu^j}(t^\rho)
\text{ for all $i,j$}\bigl).
$$
The division by the LCM provides that $J^{\{\la\},\{\mu\}}$
are polynomials in terms of $q,t$ up to some (possibly
fractional) powers of $q$ and $t$. This can be extended to any 
root systems $R$. The $a$-stabilization theorem in type $A$
is as for knots:
we arrive at $\h^{\{\la\},\{\mu\}}(q,t,a)$. 

Topologically, we consider  the $1${\small st}  link in
the horizontal solid torus and the $2${\small nd} in the
complementary vertical one. Applying the coinvariant $\{\cdot\}$, 
we obtain a topological invariant of the resulting link in
$S^3$. Changing $Y$ by $Y^{-1}$ in this definition
corresponds to changing the orientation
of the $2${\small nd} link versus the $1${\small st}.  
Algebraic links are when $\tga$ are products of 
positive powers of $\tau_{\pm}$ and
the linking numbers between the components of
the link for  $\tilde{\p}$ and that for $\tilde{\q}$ are
positive. 

For instance,
$J=\bigl\{\tau_-(\tilde{P}_{\lambda^1}
\tilde{P}_{\lambda^2})\bigr\}/LCM\bigr(\tilde{P}_{\lambda^1}(t^\rho), 
\tilde{P}_{\lambda^2}(t^\rho)\bigl)$  for the 
Hopf 2-link $L$
with the linking number $+1$ (algebraic), colored by the
Young diagrams $\lambda^1$ and $\lambda^2$.
When $\lambda^1\!=\!\yng(1)\,\!=\!\lambda^2$,
 the corresponding $\mathcal{H}$ becomes
$1-t+qt + aq$. Using the presentation of $L$ with $2$ solid tori,
we obtain the alternative formula:
$J=\bigl\{\tilde{P}_{\lambda^1}(Y)
\tilde{P}_{\lambda^2}(X)\bigr\}/LCM\bigr(\tilde{P}_{\lambda^1}(t^\rho), 
\tilde{P}_{\lambda^2}(t^\rho)\bigl)$. Their coincidence
is one of the key DAHA identities in the proof
of the topological invariance of $J$ for links. 


\subsection{\bf Superduality and RH} 
The DAHA-superduality is based on the  
$q\leftrightarrow t$-symmetry of type-$A$ {\sf\em stable} 
Macdonald polynomials and some properties of the action
of projective $PSL_2(\Z)$. 

In physics, the {\sf\em superduality} is related to 
the {\sf\em $S$-duality} in SCFT for the
BPS states, and to the {\sf\em CPT symmetry}. 
We note that various formulas 
and properties of superpolynomials were obtained and/or
conjectured by physicists.  Their works are mostly experimental, 
though the BPS states can be defined rigorously.

In topology, the superduality for 
the Khovanov-Rozansky polynomials,
{\sf\em  KhR-polynomials}, is a difficult matter.
As far as we know, it was justified only for 
positive iterated torus links when one of their components
is colored by a row or a column (uncolored otherwise).
A related problem is the
definition of the 
{\sf\em reduced} KhR-polynomials, which was resolved only
partially; see {\sf\em ``Khovanov-Rozansky homology of 
two-bridge knots and links"} (Rasmussen, 2005). Our superpolynomials are
counterparts of {\sf\em reduced} KhR-polynomials. Also, 
considering {\sf\em links} is generally involved in the KhR-theory.
Let us mention here {\sf\em Soergel modules}, an important tool
in this theory.

The DAHA construction is for {\sf\em any} iterated torus links 
and {\sf\em arbitrary} colors.
Moreover, $J^{\{\la\},\{\mu\}}$ 
can be defined for any reduced root systems and $C^\vee C_1$.  
When the theories overlap, it is expected that topological
superpolynomials,
the DAHA ones, motivic superpolynomials and those
from physics coincide up to
renormalizations. Also, there are combinatorial conjectures
and connections with the {\sf\em Heegard-Floer cohomology}.
The Alexander polynomials and the $\rho_{ab}$-invariants, discussed
below, are related to the latter.

\vskip 0.2cm

The DAHA superduality, conjectured by the author, was justified 
by Gorsky-Negut for torus knots and Cherednik-Danilenko for
iterated torus links. In terms of the standard DAHA parameters:
$\h^\la(q,t,a)=
q^\bullet t^\bullet \h^{\la'}(\frac{1}{t},\frac{1}{q},a)$, 
where by $q^\bullet t^\bullet$, we mean ``up 
to some power of $q,t$";\,  $\la'$ is the transposition of 
the diagram $\la$.

\vskip 0.2cm

The conjectural coincidence of the 
DAHA  superpolynomials $\HH$  with motivic ones 
will be stated below
for algebraic links and ``rows" (weights $m\om_1$). The
motivic superpolynomials
are defined by now only in this generality.
They are conjectured to coincide 
with the corresponding {\sf\em flagged $L$-functions} (below). The 
superduality for $L$-functions is the {\sf\em functional equation},
not very difficult to check (for $m=1$).

The $a$-stabilization 
and superduality 
are expected to hold for $B,C,D$.
For the $C$-{\sf\em hyperpolynomials} (above)
and
the transposition $\la\mapsto \la'$, 
the conjecture in
{\sf\em ``Jones polynomials of torus knots via DAHA"} was:
$$
\h^{C}(\la\,;\,q,t,u,a)=q^\bullet t^\bullet u^\bullet 
\h^{C}(\la'\,;\,t^{-1},q^{-1},u^{-1}t/q, -a q u).
$$
We conjecture in the case of $u=t$ that $\h^C$ are in terms
of $q,t^{\pm 1},\aa\equal qta^2$ (odd powers of $a$ vanish) 
for iterated torus knots; then the superduality is that for $A$:
it fixes
$\aa$ and sends $q\leftrightarrow t^{-1}$. 
Let us provide the hyperpolynomial for $T(4,3)$ for $u=t$. One has:\,
$\h^C_{4,3}(\yng(1)\,;q,t,t,\aa)$=
{\small
$$
1+q t+q^2 t+q^2 t^2+q^3 t^3
+\aa (-1+q^3-2 q t- 2q^2 t^2-q^3 t^3)
+\aa^2(q t-q^2 t+q^2 t^2).
$$
}

\noindent
We can take here $\aa=a^2 t^2$, 
which is super-invariant too. Then the restriction to $C_n$
will be $\aa=t^{2n}=t^h$ for the Coxeter number $h$.  
 
\vskip 0.2cm

{\sf\em HOMFLY-PT polynomials.}
The  $a$-stabilization of our $J$ for $q\!=\!t$ corresponds to
the relation between the 
HOMFLY-PT polynomials, $H\!O\!M(t,\aa;\la)$, and the quantum group
 invariants for $A_n$ (the WRT invariants). 
Namely, the latter are
essentially $H\!O\!M(t,\aa=t^{n+1};\la)$. The stabilization 
of the QG invariants is  connected 
with the {\sf\em Deligne category} $Rep(GL(\up))$. 
See {\sf\em
``New realizations of deformed double current algebras 
and Deligne categories"} (Etingof, Kalinov, Rains, 2020).

We note the symmetry $k\mapsto 1/k, \up\mapsto \up k$ 
is discussed in their work in Section 4.3. With $\up\mapsto -\up k$, it 
occurs in Section \ref{sec:2ndfig} below in the context of Riemann's zeta,
but because of different reasons (not because of deformed current
algebras). 
\vskip 0.2cm 

The definition of $H\!O\!M(t,\aa;\la)$
is especially simple in the uncolored case, which is for 
$\la=\yng(1)\,$ (i.e. for $\la=\om_1$ for $A_n$).
The following {\sf\em skein  relation} is sufficient to define
them (the reduced ones):
{\small
$$
\aa^{1/2}
H\!O\!M(\nwarrow\kern-10pt\nearrow\kern-10.5pt
\nearrow\kern-11pt\nearrow
\kern-11.5pt\nearrow)
\!-\!\aa^{-1/2}
H\!O\!M(\nearrow\kern-10.2pt\nwarrow\kern-10.7pt
\nwarrow\kern-11.2pt\nwarrow
\kern-11.7pt\nwarrow)\!=\!
(t^{1/2}\!-\!t^{-1/2})
H\!O\!M(\uparrow\uparrow),\ H\!O\!M(\bigcirc)=1.
$$
}

Given $\la$ (type $A$),
$H\!O\!M^\circ (t,\aa;\la)\!=\!\h^{\la}(q=t,t,\aa=-a)$
for {\sf\em iterated torus knots}, 
where $\circ$ is as above (up to $t^\bullet \aa^\bullet$) 
and $H\!O\!M$ is reduced: $1$ for the unknot.
The coincidence is due to 
Cherednik (torus knots), Morton-Samuelson (iterated torus knots), and
Cherednik-Danilenko (iterated torus links). 

We note
that the LCM-normalization
of our $\h^{\{\la\},\{\mu\}}$ must be replaced
by the division by {\sf\em one} 
$\tilde{P}_{\la^i}(t^\rho)$ (or that for one $\mu^j$)
to match {\sf\em reduced} ``HOMFLY-PT" for links. 
The reduced $H\!O\!M$ is 
defined with respect to one ``distinguished" components of
a link. Both normalizations coincide
for links with uncolored non-distinguished components
or if all colors coincide.  
Even in these cases, the passage from non-reduced KhR
polynomials of links to reduced ones is known only partially. 
In full generality,
 the LCM-normalization is  ``non-topological".

\vskip 0.2cm
For the uncolored trefoil, i.e. for $T(3,2)$ and when 
$\la=\om_1=\yng(1)\, $:\,
$H\!O\!M\!=\!\aa(t+t^{-1}-\aa),\, H\!O\!M^\circ\!=\! 1+t^2-t\aa$; 
recall that  
$\h\!=\!1+qt+qa$. 
The Alexander polynomials $Al(t)$ are generally
$H\!O\!M(t,\aa\!=\!1)/(1-t)^{\kappa-1}$  for links 
with $\kappa$ components; in particular,
$Al\!=\!t^{-1}-1+t,\, Al^\circ\!=\!1-t+t^2$ for trefoil. 

The simplest link is the Hopf $2$-plus-link, $2$
unknots with the linking number $+1$. 
Then:
$H\!O\!M=\aa^{1/2}$ 
{\large $\,\frac{1+\aa-t-t^{-1}}{t^{1/2}-t^{-1/2}}$}
and $Al^\circ=1.$  


The superduality becomes $t^{\frac{1}{2}}\to -t^{-\frac{1}{2}}, 
\aa^{\frac{1}{2}}\to \aa^{-\frac{1}{2}}$ for  $H\!O\!M$;
it is obviously compatible with 
the skein relation above (in the uncolored case).
Generally, the Young diagram  $\la$ goes
to its transpose.
The symmetry $t^{\frac{1}{2}}\to -t^{-\frac{1}{2}}$
holds for $Al$. However, it does not hold for
the Jones polynomials and for 
the  (quantum group) $A_n$-invariants. The latter are basically
$H\!O\!M(t,\aa=t^{n+1};\la)$, where the substitution $\aa=t^{n+1}$
is obviously incompatible with the superduality. 

\vskip 0.2cm

\comment{
We note that superduality is not related to the DAHA
Fourier transform, which is basically 
the (projective) action of $\si\in PSL_2(\Z)$. As we will
demonstrate, it can  
considered as an algebraic counterpart of 
the Hasse-Weil functional equation for curves over $\F_q$.
It is expected to be important in the theory of $q$-zeta functions.
}

\vskip 0.2cm

{\sf\em RH for superpolynomials.}
After our talks with Yu.I. in 2017, I focused on RH
for DAHA superpolynomials. We need to adjust the parameters:
$\H(q,t,a)\equal\h(qt,t,a)$, i.e. we switch to $q_{new}=q/t$.
Importantly, $q_{new}$ is fixed under the superduality. 
Then 
$\H(q,1/(qt),a)=q^\bullet t^\bullet 
\H(q,t,a)$ and the ``weak" (qualitative) 
RH is the claim that   
$|\xi|=1/\sqrt{q}$ for the $t$-zeros $\xi$ of  $\H(q,t,a)$ for 
{\sf\em sufficiently small $q$}. We consider RH only
for uncolored
algebraic knots; otherwise the transposition $\la\mapsto \la'$
is necessary. There is a variant with colors: for
{\sf\em rectangle} Young diagrams. 
\vskip 0.2cm

Such {\sf\em weak RH} can be
justified for (uncolored)  {\sf\em motivic} superpolynomials,
conjecturally coinciding with the DAHA ones. {\sf\em Strong RH}
for $a=0$ states:\, $|\xi|=1/\sqrt{q}$ 
 holds for  $0<q\le 1/2$ for {\sf\em any} algebraic knot.
This is the exact  bound conjecturally. 

Let us mention that strong RH  holds for any 
 $q\!>\!0$ in the case of the family of uncolored
 $T(2p+1,2)$. In this case
$\H(q,t,a\!=\!0)=\frac{1-(qt^2)^{p-1}}
{1-qt^2}$, where $q$ is``new", i.e. after the
substitution $q\mapsto qt$.  
Experimentally, RH holds for any $q<1$ only for this family. 

\vfil
We note that the value $q=1$ is special for torus knots.
Then $\H(q\!=\!1,t,a\!=\!0)$ becomes then a product of cyclotomic
polynomials due to the Shuffle Conjecture (now
a theorem). However, we are looking for the
minimal $q_0$ such that  RH holds
for {\sf\em any} $q<q_0$; this bound $q_0$  
is smaller than $1$ generally, including
sufficiently large torus knots.

\vskip 0.2cm

Numerically, the bound $q_0$ tends to $\frac{1}{2}$ for 
 $C\!ab(13+2m,2)C\!ab(2,3)$ as $m\to \infty$, which is
not proven rigorously, but probably a very difficult to check.
This is the only such family we found. These cables 
correspond to the 
singularity rings  $\r=\C[[z^4,z^6\!+\!z^{7+2m}]]$
(see below). Interestingly, the bound $q_0$
frequently become greater (better!) for multiple
cables or if the cables begin with
torus knots different from $T(3,2)$.
For instance, it is somewhat  better for 
$C\!ab(53,2)C\!ab(13,2)C\!ab(2,3)$ versus 
$C\!ab(13,2)C\!ab(2,3)$; numerically, $0.6816$ versus $0.6686$
for $a=0$.
\vskip 0.2cm 
\vfil

  This is from my paper
{\sf\em ``Riemann hypothesis for DAHA
superpolynomials and plane curve
singularities"} (2018). There are many examples
of superpolynomials there, including colored ones and links.
Weak RH can
be stated for algebraic links too; namely, the conjectural
claim is that for $a=0$ sufficiently small the number
of pairs of exceptional (non-RH) zeros is 
$\kappa-1$, where $\kappa$ is the number of
components of an algebraic uncolored link. See the paper
concerning rectangle diagrams.

Generally,
RH totally fails for {\sf\em non-algebraic}
knots/links and beyond rectangle diagrams taken as colors
for algebraic ones. It seems a really algebraic phenomenon.
Another special feature of algebraic knots is the positivity of
the coefficients of $\h$ for algebraic knots colored by rectangles.
The positivity conjecture for rectangles is the last unresolved
problem from my initial paper. There is a version for the
algebraic links (our papers with Danilenko). 
\vskip 0.2cm
 
 The substitution $q\mapsto qt$ in the passage from $\h$
to $\H$ occurs above as  
a technicality:\, the DAHA superduality
$q\leftrightarrow t^{-1}$ then becomes $q\mapsto q, t\mapsto 1/(qt)$.
However, the latter is exactly 
the Hasse-Weil symmetry from the functional equation for curves
over finite fields. There is some connection with $q$-deformations
of Riemann's zeta 
and the Dirichlet $L$-functions (below); it is based on my ``RH paper". 
Generally, it is the passage from the superpolynomials
of links to those for 
Seifert 3-folds and their special infinite sums.

\section[\sc \hspace{1em}Plane curve singularities]
{\sc Plane curve singularities}
This section provides a conjectural formula for 
superpolynomials of algebraic links colored by
``rows" in terms of the
corresponding plane curve singularities. It corresponds
to the most general case of affine Springer fibers of type $A$
and matches well the DAHA formulas. 

\subsection{\bf Basic facts}\label{sec:basic}
Algebraic links are 
intersections of plane curve singularities at $(0,0)\in \C^2$
with small $S^3\subset \C^2$ centered at $(0,0)$;
they are {\sf\em knots} for irreducible (unibranch) 
singularities. For such knots,  the corresponding
(local) singularity rings can be considered
inside $\C[[z]]$, where $z$ is the uniformizing parameter.
They are any local 
rings $\r\subset \C[[z]]$ with $2$ generators in $(z)=z\C[[z]]$ 
and the localization $\C((z))$. Such rings
are always {\sf\em Gorenstein}. 

The simplest 
topological invariants of a singularity are its
multiplicity  dim$\,\C[[z]]/\C[[z]]\mathfrak{m}$\, for
the maximal ideal $\mathfrak{m}\subset \r$, and the arithmetic
genus \,$\de=$\,dim$\,\C[[z]]/\r$, the {\sf\em Serre number}.

\vskip 0.2cm

The rings
$\r=\C[[x\!=\!z^r,y\!=\!z^s]]$ for $r,s\in \N$
such that gcd$(r,s)\!=\!1$ correspond to unibranch
{\sf\em
quasi-homogeneous
singularities} $x^s\!=\!y^r$ and torus knots $T(r,s)$.
The multiplicity is  Min$(r,s)$ and  
$\de =\frac{(r-1)(s-1)}{2}$, which is actually 
due to Sylvester (the Frobenius coin problem). 
The simplest ``non-torus" family 
is $\r\!=\!\C[[z^4,z^6\!+\!z^{7+2m}]]$ for $m\!\in\! \Z_+$,
which are 
of multiplicity $4$ and with $\de_m\!=\!8\!+\!m$.

\vskip 0.2cm
{\sf\em From families to towers.}
The simplest {\sf\em family} is
 $\r=\C[[z^r,z^{s+mr}]]$ 
for $m\in \Z_+$. Our families can be naturally interpreted as
towers of extensions of $\C[[x,y]]$ via 
the Puiseux theory. This is related to 
the theory of  
{\sf\em Drinfeld-Vl\'eduts bound} (1983) and 
the paper by Manin-Vl\'eduts {\sf\em
``Linear codes and modular curves"}
(1985). It is for growth of the arithmetic genus in some towers 
of curves $X$, which can be singular. The 
{\sf\em Artin-Schreier towers} provide some important examples
here.  

We do plane curve singularities, when the ``curve" is $1$ point,
and the 
problems becomes about finding some formulas-bounds for 
$|\j_m(\F_q)|$ for the corresponding compactified Jacobians
(below). These numbers are the values 
of motivic superpolynomials as $t=1, a=0$.  
For smooth projective curves $X$ 
over $\F_q$, the Hasse-Weil-Deligne formula can be used in terms
the eigenvalues of ``Frobenius". We use different tools,
but the functional equation and even some form of Riemann
Hypothesis work for plane curve singularities.

\vskip 0.1cm
The formulas for
$|\j_m(\F_q)|$ and our superpolynomials 
can be viewed as counterparts of 
{\sf\em Iwasawa polynomials} for class numbers 
in $\Gamma$-extensions. 
According to Barry Mazur: \,
$\Ga$-extensions can be considered as
counterparts  of abelian coverings of
$S^3$ ramified at a given link, where the Iwasawa polynomials
can be seen as counterparts of  Alexander
polynomials ($q=t, a=-1$ for us). This is for any links.
For algebraic links,
cyclic (algebraic) coverings of $P^2$ branched over  
certain singular curves (can be assumed rational)
are sufficient to consider; Libgober (1980) and others.
This is similar to our towers.

\vskip 0.2cm
\vfil
{\sf\em Valuation semigroup.}
It is one of the key in the theory of curve singularities.
The definition of
this semigroup is as follows:
$\Ga\equal\bigl\{\,\nu_z(f), 0\neq f\in \r\subset \o\equal\C[[z]]
\,\bigr\}$, where
$\nu_z$ is the valuation, the order of $z$. We readily obtain that
$\de=|\Z_+\setminus \Ga|$. Importantly, $\Ga$ 
gives the topological type  
of the corresponding algebraic knot (considered up to isotopy), 
which is due to Zariski and others. Thus, topological invariants
of rings $\r$ are exactly those expressed in terms of $\Ga$.

\vskip 0.2cm

For instance, the Alexander polynomial is immediate
via $\Ga$. Namely,
$Al^\circ$ is  
$(1-t)\sum_{\nu\in \Ga}t^\nu$ for any $\r\subset \o$
(for its $\circ$-normalization). 
For instance, it is
$(1-t)(\frac{1}{1-t}-t)=1-t+t^2$ for $T(3,2)$.
The theory of topological equivalence  
of algebraic {\sf\em links} is significantly
more ramified; {\sf\em splice diagrams} are very helpful. 
Generally, the coincidence of
semigroups for the components and the corresponding 
pairwise linking  numbers (all must be positive) is sufficient
due to the Reeve theorem.
The pairwise linking numbers can be algebraically calculated
via the ring of singularity, which is not too involved. 
\vskip 0.2cm

For $C\!ab(53,2)C\!ab(13,2)C\!ab(3,2)$ 
above (note the change $(2,3)\mapsto (3,2)$), 
the ring is $\r\!=\!$
$
\!\C[[x\!\!=\!z^8,y\!=\!z^{12}\!+\!z^{14}\!+\!z^{15}]].
$ 
The Newton's pairs are generally
$\{r_1,s_1\},\{r_2,s_2\}, \cdots\,$, and 
the Puiseux-type equation is
$y=x^{\frac{s_1}{r_1}}\Bigl(1+ c_1 
x^{\frac{s_2}{r_1 r_2}}
\bigl(1+c_2 x^{\frac{s_3}{r_1 r_2 r_3}}(\cdots)\bigr)\Bigr)$
for generic $c_i$. {\sf\em We will assume that $r_1<s_1$, 
which can be always imposed}. 

The arithmetic genus is $\de\!=\!42$, and the 
valuation semigroup $\Gamma=\lan 8,12,26,53\ran$. 
Generally, $\Ga=\lan r_1r_2r_3,\,
a_1r_2 r_3,\, a_2 r_3,\, a_3\ran$ for the cable parameters 
$(a_i,r_i)$ above (here $r_1<s_1$ is used). Recall that
$a_1\!=\!s_1$, $a_2\!=\!r_1s_1r_2+s_2,$\, 
$a_3\!=\!a_2 r_2 r_3+s_3\,$, and so on for any number 
of $\{r_i,s_i\}$. In this example, 
the Newton's pairs are $\{(2,3),(2,1),(2,1)\}$.

\comment{
The Newton's pairs are $\{(3,2),(2,1),(2,1)\}$. They result
in the Puiseux-type equation of this  
singularity, which  is { 
$x=y^{\frac{2}{3}}\bigl(1+ c_1 
y^{\frac{1}{3\cdot 2}}
(1+c_2 y^{\frac{1}{3\cdot 2\cdot 2}})\bigr)$}. 
}
\vskip 0.2cm

The passage from the base field $\C$ to finite 
fields $\F_q$ 
for $q=p^k$ and prime $p$ is sufficiently straightforward;
it will be needed below. We begin with 
$\r$ over $\C$,  define it over $\Z$, which is always 
doable {\sf\rm within a given isotopy type}, and then
consider $\r\otimes _{\Z} \F_p$. 
A prime number $p$ is called a {\sf\em prime of good
reduction} if $\Gamma$ remains unchanged 
over $\F_p$ upon this procedure. This definition is
adjusted to the topological invariance.
\vskip 0.2cm

All primes $p$ are good for the
rings $\C[[x\!=\!z^r,y\!=\!t^s]]$ as above.
Presumably,
there are no prime $p$  of bad reduction in this sense
within a given topological
type for {\sf\em any} algebraic knots: given any $p$, there exists
$\r$ representing a given knot where this $p$ is good.

To give an example, let
$\r\!=\!\Z[[x\!=\!t^4,y\!=\!t^6\!+\!t^7]]$.
Then  $\Ga\!=\!$ {\small $\{0,4,6,8,10,13,14,16,17,18,\ldots\}$} 
and $\de=8$. This $\r$ has bad reduction only at $p=2$. 
Indeed, 
$\nu_z(y^2-x^3)=14$ in $\F_2$,  which is $13$ for $p\neq 2$. However,
this singularity is equivalent over $\C$
(analytically, not only topologically) to the one for
$\Z[[t^4+t^5,t^6]]$, where bad $p$ is $3$. We obtain that 
the corresponding cable
has no primes of bad reduction. 

\subsection{\bf Compactified Jacobians}\label{sec:stand}
Let $\r\subset \o\equal\F[[z]]$ be the ring of an irreducible 
plane curve singularity over any field $\F$. The corresponding 
{\sf\em flagged compactified Jacobian}
 $\j_\ell$, considered as a set of $\F$-points by now, 
is formed by {\sf\em standard flags} 
$\vect{M}=
M_0\!\!\subset\!\! M_1\!\!\subset\!\! \cdots\!\!\subset\!\!
 M_\ell\!\subset\!
\o\!=\!\F[[z]]$ of 
$\r$-submodules $M_i$ of $\o$ such that\  \,
(a)\, $M_0\ni \phi=1\!+\!z(\cdot)$ (where $(\cdot)\in \o$), 
\, (b)\,
dim$\,M_i/M_{i-1}\!=\!1 \hbox{ and } M_i=M_{i-1}\oplus 
 \C\, z^{g_i}(1+z(\cdot))$, \,  and (c)\, (important)\,
$g_i<g_{i+1},\, \text{ where } i\ge 1$. We will call them
$\ell$-flags.

When $\ell=0$ ($0$-flags),
there is  only one condition: $\o\supset M\ni \phi=1+z(\cdot)$.
Equivalently, $\De(M)\ni 0$, where 
$\De(M)\equal \{\nu_z(v)\mid 0\!\neq\!v\!\in\! M\}$.

Generally, $\De(M)$ are $\Ga$-modules for any $\r$-modules
$M$,\, i.e. $\Ga+\De\subset \De$.
{\sf\em Standard}  $\De$ are those in $\Z_+$
containing $0$ and,
therefore, containing the whole $\Ga$. Thus, standard $M$ are
those with  standard $\De(M)$.

For quasi-homogeneous
singularities $\r=\F[[x=z^r,y=z^s]]$, where gcd$(r,s)=1,\, r,s>1$,
all standard $\Ga$-modules $\De$ 
come from some 
standard $M$. There are several ways to see this. Piontkowski used the
method of syzygies, which also gives that the corresponding
cells are affine spaces and result in combinatorial formulas
for their dimensions. Also, the $\C^*$-action and the 
{\sf\em Bialyncki-Birula theorem} can be used.

This is a special feature of quasi-homogeneous (plane curve,
unibranch) singularities;
generally, not all  $\De$ are present in the 
decomposition of $\j_0$. 
For instance, for $\F[[z^4,z^6+z^7]]$, two from $25$ such $\De$
are not in the form $\De(M)$ for any standard $M$, 
which phenomenon is due to Piontkowski. It seems that this is
always the case unless for quasi-homogeneous singularities. Also,
generally, not all Piontkowski cells are affine spaces.  
\vskip 0.2cm
\vfil

Let us supply $\j_0$ with a structure of a {\sf\em projective}
variety. We will describe the corresponding reduced scheme. 
By construction, this set is naturally 
a {\sf\em disjoint} unions of
quasi-projective varieties, those for 
different values of the deviations of $M$ (below). Importantly, they
can be combined in one projective variety. 
The main steps are as follows. 

First, any standard $M$ contains the
ideal $(z^{2\de})=z^{2\de}\o$. Indeed, the latter is
the {\sf\em conductor} of $\r$ for any Gorenstein $\r$, the greatest
ideal in $\o$ that belongs to $\r$. Using this,
 $\phi=1+z(\cdot)\in M$ (it is standard) implies that 
 $\phi\cdot (z^{2\de})=
(z^{2\de})\subset M$. 

Second, let $dev(M)\equal   
\de-\text{dim}(\o/M)$, its {\sf\em deviation} from $\r$; this is for
any $\r$-modules in $\o$.
 Then, $dev(M)\ge 0$
for standard $M$ and it is $0$ if and only if $M=\phi\r$
for some $\phi$ as above. The latter modules are called
{\sf\em invertible}. They form the {\sf\em generalized
Jacobian variety} of this singularity, which is an algebraic
group.  The third step (the key) is based on the fact that 
$z^{dev(M)}M\supset 
(z^{2\de})$ for standard $M$ due to 
Pfister-Steenbrink. Equivalently, 
$dev(M)+\De(M)\supset 2\de+\Z_+$. 
\vskip 0.2cm

Finally, let $M\mapsto M'\equal z^{dev(M)}M$.
Then $dev(M')=dev(M)-dev(M)=0$. It establishes an 
identification of standard $M$ with $\r$-modules  
$(z^{2\de})\subset M'\subset \o$ such that  $dev(M')=0$. 
The inverse map is 
$M'\mapsto z^{-d} M'$ for $d=\text{Min}\{
\nu_z(m)\!\mid\! 0\!\neq\!m\in M'\}$. Then $\{M'\}$,
the {\sf\em compactified Jacobian}, becomes a  projective
subvariety of 
 the Grassmannian of the
subspaces of the middle dimension in
 $\o/(z^{2\de})$. It is irreducible (Rego), which holds only for
plane curve singularities among all Gorenstein ones. 
 Then  $\j_\ell$ become
natural fiber spaces over $\j_0;$ the fibers are not too difficult
to describe, which will be used below.

\vfil
\vskip 0.2cm
{\sf\em Affine Springer fibers.} 
The definition requires the equation $F(x,y)=0$ for the
generators $x,y$ of $\r$. 
 Let $n$ and $m$ be the top $x$-degree and 
$y$-degree of this equation, which we assume irreducible. Then
our $\j_0$ can be interpreted as a (parahoric)
affine Springer fiber $\x_\ga$ 
defined either for  $GL_n$ or  for $GL_m$;
the equation connecting $x$ and $y$ becomes
the corresponding {\sf\em characteristic equation}.

The case of {\sf\em arbitrary} $F(x,y)$,
not irreducible and not square-free, 
will be addressed  below; topologically, this is the case
of algebraic links colored by any rows. 
\vfil

 Generally, ASF 
are due to Kazhdan-Lusztig (1988). Their description 
entirely in terms of $\r$ is
a remarkable feature of type $A$. The definition is via
$GL_n$ or via $GL_m$; but the corresponding ASF are 
isomorphic. This is not immediate from their definition (below).
The standard modules $M$ and the definition of
$\j_0$ given above do {\sf\em not} require 
the equation $F(x,y)=0$; only $\r\subset \o$ is needed. 
\vskip 0.2cm
\vfil

For semisimple Lie algebra
$\mathfrak{g}$ and any field $\F$, let $\mathfrak{g}[[x]]=
\mathfrak{g}\otimes_\F \F[[x]]$ and 
 $\mathfrak{g}((x))=
\mathfrak{g}\otimes_\F \F((x))$. Accordingly,
we define $G[[x]]$ and $G((x))$
for simply-connected $G$ with Lie$(G)=\mathfrak{g}$.

\vfil

Given $\ga\in  \mathfrak{g}[[x]]$, $\x_\ga \equal \{g\in  
G((x))/G[[x]]
\mid g^{-1}\ga g\in  \mathfrak{g}[[x]]\}$, where we assume that
the centralizer of $\ga$ in $G((x))$ is {\sf\em anisotropic}
(the nil-elliptic case).  Then $\x_\ga\cong \j_0$
in type $A$, where
the singularity is $P(x,y)=0$ for the characteristic polynomial
$P(x,y)=$\,det\,$(\mathbf 1 y-\ga)$. The corresponding 
{\sf\em orbital integral}
will become  $\h^{mot}(q,t\!=\!1,a\!=\!0)$ for
the motivic superpolynomials defined below, where $\F=\F_q$. 
It is conjectured to be a topological 
invariant, which implies that so is the orbital integral.
For instance, 
$\h^{mot}(q\!=\!1,t\!=\!1,a\!=\!0)$ is 
conjecturally the Euler characteristic $e(\x_\ga)$. 

We note that our compactified Jacobians occur
as {\sf\em Jacobian
factors} if projective {\sf\em rational} singular curves are considered;
they are basically {\sf\em Hitchin fibers} over $\mathbb P^1$.
However, {\sf\em factorizable} Lie groups and algebras (below)
are, generally,  beyond Hitchin fibers. Given a
{\sf\em factorizable Lie algebras} 
$\mathfrak{G}$, 
we considered the families of
subtori $\t\subset \g$  with fixed characteristic polynomials
in the corresponding factorizable Lie group
(Cherednik, 1983). The definitions are  as follows.

The factorizable Lie algebras  $\mathfrak{G}$ are vector bundles 
over a smooth projective curve $E$, with the
structure of {\sf\em relative}
Lie algebra over $E$. The generic fiber must 
 be $\mathfrak{g}$, but some fibers can be  
non-semisimple Lie algebras (all are of the same dimension).
 The factorization
conditions are $H^0(E,\mathfrak{G})=\{0\}=
H^1(E,\mathfrak{G})$ for \vv{C}ech cohomology,
which readily implies that genus$(E)\le 1$. If 
$E$ is singular then $\mathfrak{G}$ must be assumed
torsion free. Main applications are for $E=\mathbb P^1$. 

Such $\mathfrak{G}$
are in 1-1 correspondence with {\sf\em not necessarily unitary
classical $r$-matrices} \,$r(u,v)\in \mathfrak{g}^{\otimes 2}$: 
those satisfying the 
identity $[r^{12},r^{13}+r^{23}]=
[r^{13},r^{32}]$, where $r^{ij}$ is $r^{ij}$ for $u=u_i, v=u_j$
considered with values in $\mathfrak{g}^{\otimes 2}$ 
embedded in the components $i,j$ of   
$\mathfrak{g}^{\otimes 3}$. The parameters $u_i$ are {\sf\em local}:
near $0$. Additionally, 
we assume that
$r-\Om/(u-v)$ is regular at $0$ for the ``permutation matrix" 
$\Om\in 
\mathfrak{g}^{\otimes 2}$, the Casimir element.  
\vskip 0.2cm

The link to ASF is basically 
as follows. Let $\g$ be the group scheme over $E$ with the
Lie algebra $\mathfrak{G}.$ We obtain that 
$H^0(E,\g)$ and $H^1(E,\g)$ are trivial. The starting point 
is a subscheme $\t\subset \g$, which is assumed a maximal
subtorus at the generic point of $E$. 
Since $H^1(E,\g)=\{0\}$, any cocycle 
$\phi$ in the {\sf\em generalized Jacobian}, which is
 $H^1(E,\t)$, becomes
the boundary $\{\phi_i\phi_j^{-1}\}$ for an open cover
$E=\cup_i U_i$ and $\phi_i\in H^0(U_i,\g)$. Then  
$\t_\phi=\phi_i^{-1}\t \phi_i\subset \g$ is another
toric subscheme with 
the characteristic polynomial coinciding with that of  $\t$. 

Generally, any $\t$ can be represented
as $\t=\mathbb G_m (C)$ for a projective curve 
$C$ covering $E$, {\sf\em possibly singular}. 
Let $\overline{Jac}(C)$ be
the compactification of the generalized Jacobian  
of $C$. Then 
the {\sf\em Jacobian factors} will be the contributions 
of singular points of $C$ to $\overline{Jac}(C)$.
One can take here $E=\mathbb P^1$ and consider
rational curves  $C$ with only one singularity.
Then it will give our $\j_0$ for the corresponding $\r$. 

\comment{
Our $\j_0$ can be interpreted as a (parahoric)
{\sf\em affine Springer fiber} 
for $GL_n$ or  $GL_m$,
where $n$ and $m$ are the top $x$-degree and $y$-degree in the equation
of the corresponding  plane curve  singularity.
 Generally, ASF
are due to Kazhdan-Lusztig (1988). Their description 
entirely in terms of $\r$ is
a remarkable special feature of type $A$. One can use here
$GL_n$ or $GL_m$; the corresponding ASF coincide,  which 
is not immediate from their definition (see below).
The standard modules and $\j_0$ depend only on the {\sf\em
spectral curve}, which is given by $\r$ in our case.
\vskip 0.2cm
}

\comment{
Also, let us mention here the Oblomkov-Shende-Rasmussen conjecture
(2012), 
which is about {\sf\em nested Hilbert schemes} of pairs of
ideals in $\r$; see below.
 It is an extension of the prior Oblomkov-Shende 
conjecture proved by Maulik for HOMFLY-PT polynomials; see below.}
 
\subsection{\bf Motivic superpolynomials}
The rings $\r$ and $\r\subset \o$ will be now over 
 $\F=\F_q$. 
Following 
{\sf\em ``DAHA and plane curve singularities"}
(Cherednik-Philipp, 2017), the
{\sf\em motivic superpolynomial} of $\r$  is:
\vskip 0.1cm
\noindent
$\h^{mot}\!\equal\!
\sum_{\{M_0\!\subset\cdots\subset\!
 M_\ell\}\in\j_\ell(\F)}t^{\hbox{\tiny dim}(\o/M_\ell)}a^\ell
$
\ for all flags $\vect{M}\subset \o$, where $\ell\ge 0$. 

The flags are actually not necessary in this definition
due to the following theorem. Let 
$r\!k_q(M)\!\equal$\,dim${}_{\F_q}
M/\mathfrak{m}M$ for the maximal ideal $\mathfrak{m}$
of $\r$. Then 
$\h^{mot}=\sum_M t^{\hbox{\tiny dim}(\o/M)}(1+aq)\cdots
(1+aq^{r\!k_q(M)-1})$, where the summation is over all 
{\sf\em standard} $M\subset \o$. The justification uses
Proposition 2.3 from the paper mentioned above.
\vfil

We conjectured there 
that $\h^{mot}=\h$, i.e. the motivic one for $\r$
coincides with the uncolored DAHA superpolynomial
$\h(q,t,a)$ associated with 
the link of the singularity associated with  $\r$. The definition of
$\h^{mot}$ and this
conjecture  were extended later (with Philipp) to
torsion free sheaves of any rank$=m$ over irreducible plane curve
singularities, corresponding to the DAHA superpolynomials
$\h^\la$ for $\la=m\om_1$. The latest development is the 
generalization to 
non-unibranch singularities to be considered below. As we mentioned
above, this corresponds to affine Springer fibers of type $A$
with the most general characteristic polynomials. So it is
a natural setting here.

The DAHA superpolynomials depend on $q$ polynomially by construction
and are topological invariants (a theorem). Thus,
this conjecture includes the claims that
$\h^{mot}$ polynomially depend on $q$ and that 
these polynomials  are {\sf\em topological} invariants of the
corresponding plane curve singularities. This was justified by 
the author for some families:\, 
when $\Ga$ has  $2$ generators (the case of torus knots),
or $3$ generators; the latter was with restrictions. 
Generally, algebraic/analytic types of plane curve singularities 
depend on ``continuous" parameters;
 the classification  is essentially known and we use its elements
when considering $\Ga$ with $2-3$ generators mentioned above.

\vskip 0.2cm
 
We note  that counterparts
of the motivic superpolynomials can be defined in characteristic
$0$:\,  for any
$p$-adic integral domains $\o$ instead of $\F_q[[z]]$
and its {\sf\em orders} $\r$,
subrings with the same localization field. 
They count {\sf\em standard} $\r$-modules  $M\subset \o$,
those containing a unit in $\o$,  
with the weights $t^{deg} a^{r\!k}$. Here $|\o/M|=q^{\deg}$
for $\r/\mathfrak{m}_{\r}=\F_q$ and $r\!k=\text{dim}_{\F_q}
M/\mathfrak{m}_{\r}M$ for the maximal ideal 
$\mathfrak{m}_{\r}\subset \r$. 

There will be no quasi-projective varieties and 
{\sf\em Witt vectors} will be used, but the 
procedure is similar. 

For instance, let $\o=\Z_p[[\pi]]$ for the
$p$-adic $\Z_p$, 
$\pi^s=p$ and $\r=\Z_p[[x\!=\!p,y\!=\!\pi^r]]\subset \o$, where 
 $r,s>0$ and gcd$(r,s)=1=$gcd$(p,s$ (the tamely ramified case). 
 The corresponding  
superpolynomial will be then the same as the one for
$\r=\F_p[[z^s,z^r]]\subset \o=\F_p[[z]]$ in relatively simple examples. 
Generally, there are many possible
domains $\o$ in the $p$-adic case. Counterparts of
plane curve singularities are complete subrings in $\o$ with $1$ and
one generator (and the same field of rationals). Generally, the action of 
 \text{Gal}$(\overline{\Q}_p/\Q_p)$ becomes significantly
more involved in the $p$-adic theory.

\vfil
\vskip 0.2cm
{\sf\em Piontkowski cells.}
We set $\De(\vect{M})=\{\De(M_i)\}$. It is standard for standard
$M$ in 
the following sense.  An abstract  sequence
of $\Ga$-modules 
$\vect{\De}\!=\!
\{\De_0\!\subset\!\cdots \!\subset\!\De_\ell\subset \Z_+\}
$ is called {\sf\em standard} if $\De_0$ contains $\Ga$, 
$\De_i=\De_{i-1}\cup\{g_i\}$,
and $g_i<g_{i+1}$ for $1\le i\le \ell$. Given a standard 
$\vect{\De}$, the corresponding  {\sf\em Piontkowski cells}
is 
$\j_\ell(\vect{\De})
\equal\bigl\{ \vect{M}\in \j_\ell \mid \De(\vect{M})=
\vect{\De}\bigr\}$. 

These cells are subsets in
$\j_\ell$ and  $\j_\ell=\cup\, \j_\ell(\vect{\De})$, 
where
the union is disjoint. These cells are not always affine
spaces ${\mathbb A}^m$ and some can be empty. Empty cells
always occur (in examples) 
unless  for quasi-homogeneous singularities $x^s=y^r$;
all cells are non-empty affine spaces for them. Beyond them,
all cells can be affine spaces for  some ``non-torus" 
exceptional ``small" families. Then the lists of
empty cells and 
dim$\,\j_\ell(\vect{\De})$  for the other cells are sufficient
to know. 

The varieties $\j_\ell(\vect{\De})$   are  
conjectured {\sf\em configurations of affine
spaces}, i.e. 
unions and differences of affine spaces $\mathbb{A}^m$ 
in a bigger $\mathbb{A}^N$, not always equidimensional and
connected. This includes the            
non-unibranch generalization considered below. 
Then $\h^{mot}$ becomes
with the coefficients in
terms of $[\mathbb A^1]\in K_0(V\!ar/\mathbb{F})$ instead of $q$,
i.e. motivic indeed; $K_0(V\!ar/\mathbb{F}_q)\supset [X]\mapsto
|X(\mathbb F_q)|$ is an important
{\sf\em motivic measure}.  


Let us mention that the connection between the dimensions
of these cells in $\j_0$ and the {\sf deviations}
in the case of $\F_q[[z^r,z^s]]$ was observed by Lusztig-Smelt 
(1995). This is a special case of superduality. The deviations
are readily given in terms of $\De(M)$. 
\vskip 0.2cm
\vfil

The coincidence of $\h^{mot}$ with the
DAHA superpolynomials $\h^{daha}$
is checked in many examples, including the cases when some  
Piontkowski cells are {\sf\em not}  affine spaces. If all of them
are such, the method of {\sf\em syzygies}  
provides formulas
for their dimensions, and the calculation of $\h^{mot}$
is mostly reduced to the combinatorics of $\De(M)$; 
Dyck paths occur for $\h^{mot}_{r,s}$,
etc. Thus, $\h^{daha}=\h^{mot}$ is actually an advanced  
version of the {\sf\em Shuffle Conjecture}. 

We note that motivic superpolynomials 
are always {\sf\em significantly} faster
to calculate than {\sf\em flagged $L$-functions} defined below.
This is especially true when explicit formulas 
for  dim$\,\j_\ell(\vect{\De})$ are known (for $T(r,s)$ and for
several exceptional families of cables). Otherwise,  
DAHA calculations are, generally, faster than motivic ones,
especially with colors. 

\vskip 0.2cm
If {\sf\em some} covering of $\j_0$ by affine cells exists,  
then the coefficient of $q^i$
in  $\h^{mot}$ for $t\!=\!1,a\!=\!0$ is the {\sf\em Betti number}
$b_{2i}\!=$\,rk\,$H_{2i}(\j_{0};\R)$ and $b_{2i+1}\!=\!0.$
In particular, $\h^{mot}(q\!=\!1,t\!=\!1,a\!=\!0)$ 
is the Euler number $e(\j_0)$.
The latter is the rational Catalan number 
$\frac{1}{r+s}\binom{r+s}{r}$ for $\r=\F[[z^r,z^s]]$ 
(Beauville), where gcd$(r,s)=1$ as above. This is 
the number of all standard $\De$ for such 
 $\r$, which are 1-1 with Dyck paths 
in the rectangles ``$r\times s$"; 
the approach to $e(\j_0)$ via the count of standard
$\Ga$-modules $\De$ is due to Piontkowski. However,
this number is bigger than $e(\j_0)$ unless for
torus knots (quasi-homogeneous singularities).  
\vfil

We conjectured with Ivan Danilenko that the relation
to Betti numbers of $\j_0$ always holds for the 
corresponding DAHA 
superpolynomials.
More generally, the conjecture is
that the {\sf\em geometric superpolynomials} defined in terms
of {\sf\em Borel-Moore homology} of $\j_{\ell}$ 
coincide with the DAHA superpolynomials for any algebraic knots. 
The geometric superpolynomials coincide with 
motivic ones if $\j_{\ell}$ can be covered by affine spaces,
which is by the definition of the Borel-Moore homology. 
\vskip 0.2cm

{\sf\em From knots to links.} The consideration of
non-unibranch plane singularities 
colored by $m\om_1$ (pure rows) is necessary for the theory
of ASF   
of type $A$ with arbitrary (not only irreducible)
characteristic polynomials.
Also, they occur 
in the inductive formulas for the superpolynomials, 
similar to the Rosso-Jones formula in topology, even
if we begin with unibranch uncolored plane singularities. 

The ring will be now $\r\subset \o\equal 
\oplus_{i=1}^\kappa e_i\o_i$,
where $\o_i=\F[[z_i]]$ and $e_i e_j=\de_{ij} e_i$. 
We set $z\equal \sum_{i=1}^\kappa z_i$ and 
$e\equal\sum_{i=1}^\kappa e_i$;
then $z_i=z e_i$ and $e$ is the unit element $1$ in the ring $\o$. 
Generally, $f_i$ will be the projection $f e_i$ for any $f\in \o$. 
Here, $\r$ must contains $1=e$ and have $2$ generators: 
$x=\sum_{i=1}^\kappa x_i$ and 
$y=\sum_{i=1}^\kappa y_i$ in $\mathfrak{m}_{\o}=z\o$. 
Also, the localizations of the projection 
$\r_i$ of $\r$ onto $\o_i$ must be full $\F((z_i))$.

By construction, 
$\prod_{i=1}^\kappa
F_i(x,y)=0$, where $F_i(x_i,y_i)=0$ for the corresponding 
irreducible polynomials $F_i$ for $\r_i$, {\sf\em which will be
assumed all non-proportional}. The assumption that $F(u,v)=
\prod_{i=1}^\kappa F_i(u,v)$ 
is {\sf\em square-free} is standard for curve singularities.

\vfil
For $\F=\C$, the equation $F(u,v)=0$
gives the corresponding singularity (with $\kappa$ branches).
The corresponding link is 
$\{F(u,v)=0\}\cap S^3_\ep$ in $\C^2$
with the coordinates $u,v$; it has $\kappa$ components.
Its isotopy type gives the topological type of the singularity.

The passage from $\C$ to $\F_q$ is the same as in the 
unibranch case. Namely, we pick $x,y\in \Z[[z]]$ within
a given topological type and then switch to 
$\F_q$ for $q=p^m$ provided that $p$ is a prime of {\sf\em
good reduction}. By definition, good $p$ are such that 
the corresponding $F_i$ remain irreducible and  
pairwise non-proportional over $\F_q$. The  semigroups $\Ga_i$ for $\r_i$
and the pairwise linking numbers must remain unchanged.
The latter conditions are entirely algebraic:\, 
the linking numbers are the corresponding intersection numbers,
which can be defined via $\r$.

The notion of good reduction is necessary for the
conjectural coincidence of motivic superpolynomials
with the DAHA superpolynomials and topology; the corresponding
$p$ must be good. 
The coincidence conjecture can be extended to 
$F(u,v)=\prod_{i=1}^\kappa F_i(u,v)^{c_i}$, i.e. to
arbitrary $F$, not only square-free. The algebraic links colored
by the sequence of weights 
$\si=\{c_i\om_1, 1\le i\le \kappa\}$ occur 
on the DAHA side in this case.
The sequence $\si=\{c_i\}$ will be assumed ordered: 
$c_1\ge c_2\ge \cdots c_\kappa>0$. These inequalities
can be achieved by permuting $\{F_i\}$. 
\vskip 0.2cm

We extend the sequence of
{\sf\em uniformizing parameters} $\{z_i, 1\le i\le \kappa\}$.
It will be now $\{\ze_i, 1\le i\le \tau\}$, where 
$\tau\equal \sum_{i=1}^\kappa c_i$. The connection is as follows:
$z_1=\ze_1+\cdots+\ze_{c_1}, 
z_2=\ze_{c_1+1}+\cdots+\ze_{c_1+c_2},\,\ldots,\,
z_\kappa=\ze_{\tau-c_{\kappa}+1}+\cdots+
\ze_{\tau}.$ 

Accordingly, $\{\ep_i, 1\le i\le \tau\}$ will be the
extended sequence of idempotents:\, 
$e_1=\ep_1+\cdots+\ep_{c_1}$, and so on. We set
$\Om=\F[[\ze_1,\ldots, \ze_{\tau}]]=\sum_{i=1}^{\tau} \Om_i$
for $\Om_i=\ep_i\Om$, where $1\le i\le \tau$, and
$\o_i=e_i \Om$ for $1\le i\le \kappa$.
 The prior $\r$ can be naturally embedded
in $\Om$:\, 
$\r\subset \o=\sum_{i=1}^\kappa \o_i \subset \Om$,
i.e. we embed diagonally for the segments of $\{\ze_i\}$
associated with the multiplicities $c_i$.  
\vskip 0.2cm

{\sf\em Standard modules} $\m$ are by definition
$\r$-invariant $\F$-subspaces $\m\subset \Om$
such that $\Om \m=\Om$. 
As above, $r\!k_q(\m)\!\equal$\,dim${}_{\F_q}
\m/\mathfrak{m}_{\r}\,\m$, where 
$\mathfrak{m}_{\r}=\r\cap z\o$. 

The minimal $q$-rank is then $r\!k_{min}=c_1$.
The maximal $q$-rank, $r\!k_{max}$, is that for $\m=\Om$.
It equals $\sum_{i=1}^\kappa m_i c_i$, where $m_i$ are mutiplicities
of the singularities associated with
$F_i(x,y)=0$, which are $m_i=\min\{\Ga_i\setminus\{0\}\}$
for the corresponding $\Ga_i$. 

The 
{\sf\em motivic superpolynomial} of $\r,\si$ is defined as follows:
$$
\h^{mot}_{\si}=\sum_{\m} t^{\hbox{\tiny dim}(\Om/\m)}
\prod_{j=r\!k_{min}}^{r\!k_q(\m)-1}(1+aq^j)
\text{ summed over standard } \m.
$$  
The first product $\Pi$ here (for $r\!k_q(\m)=r\!k_{min}=c_1$) is $1$.
The conjecture is that $\h^{mot}$  depend polynomially on $q,t,a$
and coincide with the DAHA superpolynomials  
$\h^{{\la}}(q,t,a)$ for the corresponding links colored
by the sequences $\{\la\}=\{c_1\om_1,c_2\om_1,\ldots, 
c_\kappa \om_1\}$ (only pure rows). In particular, they are topological
invariants. The conjecture is well-checked.

\subsection{\bf Some examples}
Let us begin with  $\j_\0$ for
$\r=\F[[z^4,z^6+z^7]]$ with $K=C\!ab(13,2)C\!ab(2,3)$ 
discussed above. All cells are affine spaces 
and we show only dim$\,=$dim\,$J_0(\De)$ in the table
below for the corresponding
sets of gaps $D\equal \De\setminus \Ga$ for standard $\De$. 
One has $dev(\De)=|D|$ and dim\,$\o/M=\de-|D|$, which gives the
power of $t$. 
Two standard $\De$
from $25$ have no standard $M$, namely for
$D=[2,15]$ and $D=[2,11,15]$. The table of $D$ and the
corresponding dimensions of the cells $\j_0(\De)$ is:
{\footnotesize
\begin{table*}[ht!]
\[
\centering
\begin{tabular}{|l|l|}
 \hline 
\hbox{$D$-sets} & $dim$\\
\hline
$\varnothing$ & 8\\
15 & 7\\
11,15 & 6\\
7,11,15 & 6\\
9,15 & 7\\
9,11,15 & 5\\
7,9,11,15 & 4\\
3,7,9,11,15 & 4\\
5,9,11,15 & 5\\
5,7,9,11,15 & 3\\
3,5,7,9,11,15 & 2\\
1,5,7,9,11,15 & 4\\
\hline
\end{tabular}
\hspace{0.1cm}
\begin{tabular}{|l|l|}
\hline
\hbox{$D$-sets} & $dim$\\
\hline
1,3,5,7,9,11,15 & 2\\
2,7,11,15 & 6\\
2,9,15 & 7\\
2,9,11,15 & 6\\
2,7,9,11,15 & 5\\
2,3,7,9,11,15 & 4\\
2,5,9,11,15 & 5\\
2,5,7,9,11,15 & 3\\
2,3,5,7,9,11,15 & 1\\
1,2,5,7,9,11,15 & 3\\
1,2,3,5,7,9,11,15 & 0\\
2,15 and 2,11,15 & $\emptyset$\\   
\hline   
\end{tabular}
\]
\end{table*}
}

\comment{
We provide the top ranks of $M$ corresponding to $D$-sets.
{\footnotesize
\begin{table*}[ht!]
\[
\centering
\begin{tabular}{|l|l|l|}
 \hline 
\hbox{$D$-sets} & $dim$ & $rk$\\
\hline
$\varnothing$ & 8 & 1\\
15 & 7 & 2\\
11,15 & 6 & 2\\
7,11,15 & 6 & 2\\
9,15 & 7 & 2\\
9,11,15 & 5 & 3\\
7,9,11,15 & 4 & 3\\
3,7,9,11,15 & 4 & 2\\
5,9,11,15 & 5 & 2\\
5,7,9,11,15 & 3 & 3\\
3,5,7,9,11,15 & 2 & 3\\
1,5,7,9,11,15 & 4 & 2\\
\hline
\end{tabular}
\hspace{0.1cm}
\begin{tabular}{|l|l|l|}
\hline
\hbox{$D$-sets} & $dim$ & $rk$\\
\hline
1,3,5,7,9,11,15 & 2 & 3\\
2,7,11,15 & 6 & 2\\
2,9,15 & 7 & 2\\
2,9,11,15 & 6 & 3\\
2,7,9,11,15 & 5 & 3\\
2,3,7,9,11,15 & 4 & 3\\
2,5,9,11,15 & 5 & 3\\
2,5,7,9,11,15 & 3 & 4\\
2,3,5,7,9,11,15 & 1 & 4\\
1,2,5,7,9,11,15 & 3 & 3\\
1,2,3,5,7,9,11,15 & 0 & 4\\
\,2,15 and 2,9,15 & $\emptyset$ & 0\\
\hline   
\end{tabular}
\]
\end{table*}
}
}

The whole (uncolored) superpolynomial is: 
$\h(q,t,a)=$
\renewcommand{\baselinestretch}{1.2} 
{\small
\(
1 + q t + q^8 t^8 + q^2 \bigl(t + t^2\bigr) 
 + 
 q^3 \bigl(t + t^2 + t^3\bigr) + q^4 \bigl(2 t^2 + t^3 + t^4\bigr) 
+  q^5 \bigl(2 t^3 + t^4 + t^5\bigr) + q^6 \bigl(2 t^4 + t^5 
+ t^6\bigr) + 
 q^7 \bigl(t^5 + t^6 + t^7\bigr) + 
a \bigl(q + q^2 \bigl(1 + t\bigr) + q^3 \bigl(1 + 2 t + t^2\bigr) 
+ q^4 \bigl(3 t + 2 t^2 + t^3\bigr) + 
    q^5 \bigl(t + 4 t^2 + 2 t^3 + t^4\bigr) + q^6 \bigl(t^2 
+ 4 t^3 + 2 t^4 + t^5\bigr) + 
    q^7 \bigl(t^3 + 3 t^4 + 2 t^5 + t^6\bigr) + q^8 \bigl(t^5 
+ t^6 + t^7\bigr)\bigr)+  
a^2 \bigl(q^3 + q^4 \bigl(1 + t\bigr) + q^5 \bigl(1 + 2 t + t^2\bigr) 
+ q^6 \bigl(2 t + 2 t^2 + t^3\bigr) + q^7 \bigl(2 t^2 + 2 t^3 
+ t^4\bigr) + q^8 \bigl(t^3 + t^4 + t^5\bigr)\bigr) +
a^3 \bigl(q^6 + q^7 t + q^8 t^2\bigr).
\)
}
\renewcommand{\baselinestretch}{1.2} 

For instance, there are $3$ cells of dimensions $7$ in $\j_0$
(for $a=0$). Namely, those  with $D=[15],[9,15],[2,9,15]$
and  $t^7, t^6, t^5$. Generally, the number of cells of
dim\,$=\de-1$ is the multiplicity of singularity;
it equals the coefficient of $t$ for $q\!=\!1,a\!=\!0$  due to the
superduality, which is for $\Z_+\!\setminus \De=\{1\},\{2\},\{3\}$ in this
example. Only $\{1\}$ results in dim\,$=1$ and $qt$ in $\h^{mot}$. 
We note 
the {\sf\em reciprocity involution} of standard modules
$\De\mapsto \De^\vee- \min\{\De^\vee\}$ for
$\De^\vee= \Ga\setminus\{(2\de-1)-D\}$, which preserves
their {\sf\em dim}. For instance, $[15]\mapsto [2,9,15]$
$[2,9]\mapsto [2,9]$.


\vskip 0.2cm

In the case of  $C\!ab(53,2)C\!ab(13,2)C\!ab(3,2)$ discussed above
and the corresponding ring $\r=\F[[z^8,z^{12}\!+\!z^{14}\!+\!z^{15}]]$,
one has:  $\h(q,t\!=\!1,a\!=\!0)=$
\renewcommand{\baselinestretch}{1.2} 
{\small
\(
q^{42}+7 q^{41}+24 q^{40}+56 q^{39}+104 q^{38}+166 q^{37}
+236 q^{36}+306 q^{35}+370 q^{34}+424 q^{33}+465 q^{32}
+492 q^{31}+507 q^{30}+510 q^{29}+504 q^{28}+488 q^{27}
+466 q^{26}+437 q^{25}+406 q^{24}+370 q^{23}+335 q^{22}
+298 q^{21}+264 q^{20}+230 q^{19}+199 q^{18}+168 q^{17}
+143 q^{16}+118 q^{15}+97 q^{14}+78 q^{13}+63 q^{12}+48 q^{11}
+38 q^{10}+28 q^9
+21 q^8+15 q^7+11 q^6+7 q^5+5 q^4+3 q^3+2 q^2+q+1.
\)
}
\renewcommand{\baselinestretch}{1.2}

\comment{
{\small
\begin{align*}
&q^{-42}\h(q,t\!=\!1,a\!=\!0)=
1\!+\!7 q^{-1}\!+\!24 q^{-2}\!+\!56 q^{-3}\!+\!104 q^{-4}\!+\!166 
q^{-5}\\
&+236 q^{-6}+306 q^{-7}
+370 q^{-8}+424 q^{-9}+465 q^{-10}+492 q^{-11}\\
&+507 q^{-12}
+510 q^{-13}+504 q^{-14}+488 q^{-15}+466 q^{-16}+437 q^{-17}
\\
&+406 q^{-18}+370 q^{-19}+335 q^{-20}+298 q^{-21}+264 q^{-22}
+230 q^{-23}\\
&+199 q^{-24}+168 q^{-25}+143 q^{-26}+118 q^{-27}
+97 q^{-28}+78 q^{-29}\\
&+63 q^{-30}+48 q^{-31}+38 q^{-32}+28 q^{-33}+21 q^{-34}+
15 q^{-35}+11 q^{-36}\\
&+7 q^{-37}+5 q^{-38}+3 q^{-39}+2 q^{-40}
+q^{-41}+q^{-42}.
\end{align*}
}
}

\vskip 0.2cm
Here $\de=42$, which corresponds to $q^{42}$ (invertible modules).
The coefficients of $q^{i}$ are the Betty number
$b_{2i}$ (the odd ones vanish), and the 
 Euler number
$e(\j_0)$ is $8512$ (which is for  $a\!=\!0,t\!=\!1,q\!=\!1$). 

\vskip 0.2cm

 The simplest $\h^{mot}$  is for trefoil $T(3,2)$. Its
singularity ring is $\r\!=\!\F_q[[z^2,z^3]]$ with 
$\Ga=\mathbb Z_+ \setminus \{1\}$. There are no primes of bad
reduction for this and any torus knots. 
The standard modules are $M_\la=(1+\la z)$ 
(invertible ones) of  dim$\,\o/M=1$ and $M=\o$, where there are 
$2$ generators (dim$=\!0$).
The standard flags for $\ell=1$  are 
$\{M_0=M_\la\subset M_1=\o\}$; the dimension is dim$\,\o/M_1=0$
for them.  Thus 
$\h^{mot}=1\, (\hbox{for\,} \o)+qt\, (\hbox{counting invertible modules})
+aq\,(\text{counting flags})$. This calculation is almost equally
simple for $T(2p\!+\!1,2).$
\vskip 0.2cm
\vfil

{\sf\em Links.} The simplest example is the  uncolored
Hopf (algebraic) 2-link, the one for $F(u,v)=uv$.
We take 
$\o=\F_q[[e_1,e_2,z_1,z_2]]$ for the idempotents $e_i$ such that
$z_ie_j=\de_{ij}z_i$, and $\r=\F_q[[e_1+e_2, x=z_1, y=z_2]]$,
where $e_1+e_2$ is the unit $1$ in $\o$. 
The standard modules are $\m_\circ=\o$ of $q$-rank $2$ and 
dim$(\o/\m)=0$,
and 
$\m_\al=(e_1+ \al e_2)\r$ for $0\neq \al\in \F_q$ with $r\!k\!=\!1$
(invertibles) and dim$(\o/\m)=1$.
 The latter contain $z_1, z_2$:\,
 $(e_1+ \al e_2)x=\al z_1, 
(e_1+ \al e_2)y=\al z_2$. Thus, $\h^{mot}=(q-1)t+(1+aq).$

Hopf links are quite interesting. Let us
provide $\h^{mot}_{2,1,1}$ 
 for the Hopf $3$-link with the linking numbers $\{+1,+1,+1\}$ and the
1{\small st} component colored
by $2\om_1=\yng(2)\,$. We take $\r=\F_q[[1,x,y]]$
for $ x=\ze_1+\ze_2+\ze_3+\ze_4, y=\ze_1+\ze_2-\ze_3+c\ze_4$, where
$c\neq \pm 1$ and $p\neq 2$. 
The corresponding embeddings are:
$
\r\subset \o=\F_q[[\ep_1+\ep_2, \ep_3,\ep_4,
\ze_1+\ze_2,\ze_3,\ze_4]]\subset \Om=\F_q[[\ep_i,\ze_i, 1\le i\le 4]].
$
One has: 
\vfil

\noindent
\renewcommand{\baselinestretch}{1.2} 
{\small
\( \text{\normalsize $\h^{mot}_{2,1,1}=$}
1+a^2 q^5-2 t+q^2 t+q^3 t+t^2-2 q^2 t^2+q^4 t^2+q^2 t^3
-q^3 t^3-q^4 t^3+q^5 t^3 +a \bigl(q^2+q^3-2 q^2 t+q^4 t+
q^5 t+q^2 t^2-q^3 t^2-q^4 t^2+q^5 t^2\bigr)\\
=(q-1)^2 (1+q) qt^2 (1+q t)\, +\, 
\bigl((q-1) t (2+2 q+q^2+(q^2-1)t)\bigr)
(1+aq^2)\, +\,  (1+aq^2)(1+aq^3).
\)
}
\renewcommand{\baselinestretch}{1.2}

\section[\sc \hspace{1em}Zetas for singularities]
{\sc Zetas for singularities}
\subsection{\bf Generalizing Galkin's zeta} 
V.M. Galkin studied in 1973 zeta- and $L$-functions
for Gorenstein rings in dimension one. 
Plane curve singularities are an important particular case.
Let us consider the unibranch case: when 
$\r\subset \o=\F_q[[z]]$ with $2$ generators and
the localization $\F_q((z))$. 
Recall that $\de=$dim\,${}_{\F_q}\o/\r$
$=|\,\Z_+\setminus \Ga\,|$.
He and St\"ohr (in 1998) considered links; adding $a$ and the 
case of arbitrary ranks closely follow what we did for the 
superpolynomials.                                          
\vfil

The {\sf\em admissible  flags} of ideals in $\r$ are 
$\vect{M}\!=\!\{M_0\!\subset\! M_1\!\subset\!
\cdots\!\subset\!  M_\ell\!\subset\! \r\}$ such that 
$\{z^{-m_0}M_i\subset \o\}$ 
for $m_0=\text{Min}(\De_0)$ are {\sf\em standard flags} in
$\o$ as in
Section \ref{sec:stand}.
The {\sf\em flagged zeta function} is:
$$Z(q,t,a)\!\!\equal\!
\sum_{\vect{M}} 
a^\ell t^{\hbox{\tiny \,dim}(\r /M_{\ell})}=
\!\sum_M t^{\hbox{\tiny \,dim}(\r/M)}(1\!+\!aq)\cdots
(1\!+\!aq^{r\!k_q(M)-1}),
$$
where the summation is over all admissible flags
$\vect{M}\subset \r$
and over all ideals $M\subset\r$ in the 2{\small nd}
formula; dim$=$dim$_{\F_q}$,
$r\!k_q(M)=$\,dim\,$M/\mathfrak{m}M.$ 

The 2{\small nd} formula becomes
with $\prod_{i=r\!k_{min}}^{r\!k_q(M)-1}(1\!+\!aq^i)$ 
in the non-unibranch case, similar to that 
for the motivic superpolynomials; $r\!k_{min}$ is the
maximal multiplicity of irreducible factors in the 
factorization of the corresponding $F(u,v)$. 
 The interpretation
in terms of the standard flags, which is the 1{\small st} formula, 
becomes somewhat technical in the non-unibranch case.

The 
{\sf\em flagged $L$-function} is then
$L(q,t,a)\!\!\equal\!\! (1-t)Z(q,t,a)$; it is a polynomial
in terms of $t,a$, and $t^{-\de}L(q,t,a)$ is invariant
under $t\mapsto 1/(qt)$, which is the {\sf\em 
functional equation}; we note that $Z(q,t,a)$ does not satisfy
the latter. 
In contrast to the smooth case, the Riemann Hypothesis 
holds only for sufficiently small $q$ (if we know that the
dependence on $q$ is polynomial).

The definition of the Galkin zeta, which is $Z(q,t,a)$ for $a=0$
(no flags), is sufficiently standard: a Dirichlet series. 
The functional equation for $L$ is actually
surprising because there is no Poincar\'e duality for 
singular varieties (unless intersection cohomology is used or so).
St\"ohr found a short entirely combinatorial proof of this fact, 
a significant
simplification of that from the John Tate's 
thesis. Tate's $p$-adic proof works well for curve singularities.
\vskip 0.2cm

The key here is the following property of $\Ga$, which is 
actually the defining property of {\sf\em Gorenstein rings}\,: 
the map $g\mapsto g'=2\de\!-\!1\!-\!g$
identifies $\{g\in \Z_+\setminus \Ga\}$ (the set of ``gaps") with
$\{g'\in \Ga\setminus \{2\de+\Z_+\}\}$. For instance, the last
gap, which is $2\de-1$, maps to $g'=0$.
\vskip 0.1cm

Let
$\H^{mot}(q,t,a)\!\equal\!\h^{mot}(qt,t,a)$
for motivic $\h^{mot}$ above:  we switch
to $q_{new}=q/t$ as we did for the DAHA superpolynomials
when defining $\H(q,t,a)$.
 The $\H^{mot}\leftrightarrow L$ {\sf\em coincidence conjecture} 
reads:
\begin{align*}
&\H^{mot}(q,t,a)\!=\!L(q,t,a) \text{ for 
any plane curve  rings $\r\subset \o$},\\ 
&\H^{mot}(q,t,a\!\!=-1/q)\!=\! L_{\hbox{\tiny prncpl}}
(q,t)\, \text{ for Gorenstein $\r\subset \o$}.
\end{align*}. 
 The latter 
is the {\sf\em Z\'u{\oldt{n}}iga 
zeta function}: for $a=0$ and when
the summation is only over {\sf\em principle} 
$M\subset\r$. It coincides with $L(q,t,a\!=-1/q)$, indeed;
see the alternative formula for flagged $L$ provided above.
The conjecture can be extended to the non-unibranch case and
to any ranks (in the DAHA setting, any colors $m\om_1$).  

When $q\to 1$ (for the ``field" with $1$ element):\,
$Z_{\hbox{\tiny prncpl}}=\sum_{\nu\in \Ga}t^\nu$ and
$L_{\hbox{\tiny prncpl}}=(1-t)Z_{\hbox{\tiny prncpl}}$
is the Alexander polynomial
$Al^\circ(t).$ Namely,
$$
\lim_{q\to 1}L_{\hbox{\tiny prncpl}}=(1-t)\bigl(\sum_{i=1}^{\de}
t^{g_i}\bigr)+t^{2\de} \text{ for } \{g_i\}=\Ga \setminus 
(2\de + \Z_+).
$$
We obtain that $\bigl(L_{\hbox{\tiny prncpl}}-t^{2\de}\bigr)/(1-t)$
becomes $\de$ when $q\to 1$ and $t=1$. 
\vskip 0.2cm

The conjectural coincidence of $\H(q,t,a)$ (DAHA), 
$\H^{mot}(q,t,a)$ and $L(q,t,a)$
identifies the superduality for the former with the
functional equation for the later. 
Recall that the conjectural upper RH-bound for $\H(q,t,a\!=\!0)$
is $q\le 1/2$ in the unibranch uncolored case, which is
 far from ``arithmetic"
$q=p^m$ for $L(q,t,0)$. 
\vskip 0.2cm

The coincidence $\H^{mot}(q,t,a)$ and $L(q,t,a)$ is simple
for $t=1$ (for any rings $\r\subset \o$, not only Gorenstein ones).
Indeed, any admissible flag of ideals $\vect{M}'\subset \r$
is $\,z^{m'}\, \vect{M}\,$ for standard $\vect{M}\subset \o$, where
$m',\vect{M}$ are uniquely determined  by $\vect{M}'$. 
Vice versa, given a standard flag $\vect{M}$, let:
$$\{m \mid z^{m}\vect{M}\subset \r\}
=\{0\le m_1<m_2<\cdots<m_k<2\de\}\cup \{2\de+\Z_+\}$$ 
for some $k=k_M$ and $\{m_i\}$.  
The contribution of $z^m \vect{M}\subset \r$ for such $m$ to 
$L(q,t,a)=(1-t)Z(q,t,a)$ is 
$(1-t)t^{m-dev(M_\ell)}$. Recall:
$dev(M)=\de -\text{dim}\,\o/M$.
Thus, all such $z^m \vect{M}$ contribute
$(1-t)t^{-dev}
\bigl(\sum_{i=1}^{k_M} t^{m_i}+t^{2\de}/(1-t)\bigr)$. This will be
$1$ for $t\to 1$ (and any $q)$. We obtain that $L(q,t=1,a)=
\sum_{\ell=0}^{2\de-1} |\j_{\ell}(\F_q)|a^\ell$, which is
$\H^{mot}(q,t=1,a)$.
\vskip 0.2cm

{\sf\em Using Hilbert schemes.}
For a {\sf\em rational} projective curve $C\subset \mathbb P^2$,
 the following identity
is a natural object for physicists and mathematicians
(Gopakumar-Vafa and Pandharipande-Thomas):
$
 \sum_{n\!\ge\! 0} q^{n\!+\!1\!-\!\de} e(C^{[n]})\!=\!
\sum_{0\le i\le \de} n_{C}{(i)}\bigl(\!\frac{q}
{(1\!-\!q)^2}\!\bigr)
^{i+1\!-\!\de}  
$
for the Euler numbers of {\sf\em Hilbert schemes} 
 $C^{[n]}$. The points of the latter are
zero-cycles of $C$, collections of  ideals at any points, of the
(total) colength $n$. Here $\de$ is the arithmetic genus of $C$,
$n_{C}(i)$ are some numbers. The passage from a series 
to a polynomial is far from obvious
even in this relatively simple case.
It is significantly more subtle 
to prove that $n_{C}{(i)}\!\in\! \Z_+$ 
(G\"ottsche, and then Shende for all $i$); 
versal deformations of 
singularities appeared necessary in the Shende's proof.


Switching to local rings $\r$ of singularities, the 
following conjecture is 
for {\sf\em nested\,}  Hilbert schemes
$H\!ilb^{\,[l\le l+m]}$ formed by pairs  of ideals $\{I,I'\}$ 
of colengths $l,l+m$ such that
$\mathfrak{m}I\!\subset\!
I'\!\subset\! I\!\subset\!\r$
 for the maximal ideal $\mathfrak{m}\subset \r$.
One needs the {\sf\em weight $t$-polynomial} 
$\, \mathfrak{w}(H\!ilb^{\,[l\le l+m]})$ defined for the weight
filtration of $H\!ilb^{\,[l\le l+m]}$
due to  Serre and Deligne.
 The
{\sf\em Oblomkov-Rasmussen-Shende conjecture} (2012) states that 
$$
\sum_{l,m\ge 0}q^{2l}\,a^{2m}\,
t^{m^2}\,\mathfrak{w}(H\!ilb^{\,[l\le l+m]})$$
is proportional to the Poincar\'e  series of 
the HOMFLY-PT triply graded
homology  of the
corresponding link. The connection with the perverse
filtration of $\j_0$  is due to Maulik-Yun
and Migliorini-Shende. The ORS series is a geometric
variant of $Z(q,t,a)$.

The ORS-conjecture adds $t$ to 
the {\sf\em Oblomkov-Shende conjecture}, which was
extended by adding  colors $\la$ and then 
proved by Maulik.

The passage from series to polynomials required superpolynomials.
The corresponding topological ones are {\sf\em reduced} KhR-polynomials.
The  Cherednik-Danilenko conjecture was that they coincide with
the uncolored DAHA  
$\h(q,t,a)$. The passage  
to polynomials is a nontrivial  step. This  is manifest
for the DAHA and motivic superpolynomials. Though it is a conjecture
that the motivic ones for plane curve singularities
depend on $q$ polynomially.   


\subsection{\bf The passage to links} 
Let us provide a definition of $L$-functions for 
{\sf\em non-unibranch} 
plane curve singularities colored by rows. Without the colors 
and for $a=0$, such $L$ are due to Galkin (1973)
and St\"ohr (1998). The following definition is
preliminary for {\sf\em colored} links.  

Recall that $\mathcal{H}^{mot}_\sigma$ was defined for 
$\sigma=(c_1\ge c_2\ge \cdots \ge c_{\kappa}>0)$ and 
$\mathcal{R}=\mathbb{F}_q[[x,y]]\subset \o\!=\!
\mathbb{F}_q[[z,e_1,\ldots,e_\kappa]]\subset 
\Omega\!=\!\mathbb{F}_q[[z,\epsilon_1,\ldots,\epsilon_\tau]]$ for
$\tau=\sum_{i=1}^\kappa c_i$,
$\epsilon_i\epsilon_j=\epsilon_i\delta_{ij}$, 
$e_1\!=\!\epsilon_1+\cdots
+\epsilon_{c_1}, 
e_2\!=\!\epsilon_{c_1+1}+\cdots+\epsilon_{c_1+c_2}$, etc. 
Accordingly, 
$z_i=z e_i (1\le i\le \kappa)$ and  
$\zeta_i=z\epsilon_i (1\le i\le \tau)$. 
Let $\varpi_{\mathsf{m}}=\sum_{i=1}^\kappa\epsilon_{m_i}$
for all sequences $\mathsf{m}=(m_1,m_2,\ldots,m_\kappa)$
such that $0\le m_1<c_1,\, 0\le m_2<c_2$, 
$ \tilde{\mathcal{R}}=
\mathbb{F}_q[[x\varpi_{\mathsf{m}},y\varpi_{\mathsf{m}}]]$,
$ \tilde{\mathfrak{m}}=
(x\varpi_{\mathsf{m}},y\varpi_{\mathsf{m}})$, and
$\tilde{\Omega}=\oplus_{\mathsf{m}}\mathcal{R}\varpi_{\mathsf{m}}=
\tilde{\mathfrak{m}}\oplus \sum_{\mathsf{m}}
\mathbb{F}_q\varpi_{\mathsf{m}}.$
Note that
$\tilde{\Omega}=\tilde{\mathcal{R}}=
\mathcal{R}$ in the uncolored case.
For $r\!k_q(M)=\text{dim\,}(M/\tilde{\mathfrak{m}}M)$, 
dim$\,=\text{dim\,}_{\mathbb{F}_q}$, we set: 
$$
Z_{\sigma}(q,t,a)\!=\!\sum_{M} 
t^{\text{dim}(\tilde{\Omega}/M)}\!\!\prod_{j=c_1}^{r\!k_q(M)-1}
\!(1+aq^j),\ 
L_\sigma(q,t,a)\!=\!(1\!-\!t)^\tau Z(q,t,a),
$$
where 
$M\subset \tilde{\Omega}$  are $\tilde{\mathcal{R}}$-invariant 
such that 
$r\!k_q(M_i)\ge c_i$ for  $M_i=M\,e_i$. 
\vskip 0.2cm

Setting $\mathbf{H}_\sigma(q,t,a)=
\mathcal{H}^{mot}_\sigma (qt, t, a)$ (as above), we conjecture that
$\mathbf{H}_\sigma(q,t,a)=
L_\sigma(q,t,a)$ {\sf\em at least in the uncolored case}, when 
$\sigma\!=\!(1,\ldots,1)$.
\vskip 0.2cm

Let $\sigma=(2,1)$ for the Hopf $2$-link. Then, 
$\mathcal{R}=\mathbb{F}_q[[x=\zeta_1+\zeta_2,y=\zeta_3]]$,
$\tilde{\mathcal{R}}=\mathbb{F}_q[[\zeta_j]]$  and
$\tilde{\Omega}=\mathbb{F}_q(\epsilon_1\!+\!\epsilon_3)\oplus
\mathbb{F}_q(\epsilon_2\!+\!\epsilon_3)\oplus
(\zeta_j).$
It is easy to see that $\mathcal{H}^{mot}_\sigma=
(1+aq^2)+(q^2-1)t$; the corresponding $L$ is as follows.

The modules $M\subset \tilde{\Om}$
of $r\!k_q(M)=2$ are those generated over $\tilde{\mathcal{R}}$
by $\{\epsilon_1\!+\!\epsilon_3,\epsilon_2\!+\!\epsilon_3\}$ (containing
$\tilde{\mathfrak{m}}$),
$\{\epsilon_1+\epsilon_3,\zeta_2^u\}$, 
$\{\epsilon_2+\epsilon_3,\zeta_1^u\}$, 
and $\{\zeta_1^u+\alpha\zeta_3^w,
\zeta_2^v+\beta\zeta_3^w\}$,
where  $u,v,w\ge 1$,\, $\alpha,\beta\in \mathbb{F}_q$ 
such that $\alpha\neq 0$ or $\beta\neq 0$. 
They contribute $1+\frac{2t}{1-t}+\frac{(q^2-1)t^3}{(1-t)^3}$ to
$Z_\sigma(q,t,a)$.

The submodules of $q$-rank $3$ are
those generated by 
$\{\zeta_1^u, \zeta_2^v,\zeta_3^w\}$,  which contribute
$(1+a q^2)\frac{t^2}{(1-t)^3}$. Combining and switching to $L$, we obtain:
 $L_\sigma(q,t,a)=(1-t)^3 Z_\sigma(q,t,a)=
(q^2-1)t^3+(1+t)(1-t)^2+t^2(1+aq^2)=q^2t^3+1-t+t^2 a q^2=(1+aq^2t^2)
+(q^2t^2-1)t.$ 

More generally, let $\sigma\!=\!(m,1)$. Then   
$\mathcal{H}^{mot}_\sigma\!\!=\!
(1\!+\!q^ma)\!+\!(q^m\!-\!1)t$ and $\mathbf{H}^{mot}_\sigma=
(1\!+\!q^m t^ma)+(q^m t^m\!-\!1)t$, which coincides with
$L_\sigma$.

\vskip 0.2cm

{\sf\em Double trefoil.} Let us provide the
motivic uncolored superpolynomial for $T(6,4)$
and comment on it and its $L$-counterpart.
One has:
$\mathcal{H}^{mot}_{6,4}(q,t,a)=$
\renewcommand{\baselinestretch}{1.2} 
{\small
\(
q^8 t^8-q^7 t^8
+q^3 t^2 \bigl(q+t-2 q t+q^2 t-2 q^2 t^2+2 q^3 t^2
-q^2 t^3+q^4 t^3-q^2 t^4+q^4 t^4-q^3 t^5+q^4 t^5\bigr) (1+a q)
+q^2 t \bigl(1+q-t+q^2 t-t^2-q t^2+q^2 t^2+q^3 t^2-q t^3+q^3 t^3-q^2 t^4
+q^3 t^4\bigr) (1+a q) (1+a q^2)
+\bigl(1-t+q t-q t^2+q^2 t^2\bigr) (1+a q) (1+a q^2) (1+a q^3).
\)
}
\renewcommand{\baselinestretch}{1.2}

It coincides with
$\mathcal{H}^{daha}_{6,4}$ in formula (6.4) of Section 6.2
of {\sf\em ``DAHA approach to iterated torus links"}
(Ch, Danilenko, 2015); 
the notation there was $\hat{\mathcal{H}}^{min}$.
Also, it coincides with the reduced (uncolored)
Khovanov-Rozansky polynomial defined and obtained
via {\sf\em Soergel modules}, which is from Example 1.3
in {\sf\em ``Torus link homology"} (Hogancamp-Mellit, 2019).
Their expression 
is $\mathcal{H}^{daha}_{6,4}(q\mapsto 1/t, t\mapsto q, a)$ 
multiplied by $\frac{1+a}{(1-q)^2}$. 

Here $\mathcal{R}=\mathbb{F}_q[[x=z_1^2-z_2^2,y=z_1^3+z_2^3]]
\subset \mathcal{O}=\mathbb{F}_q[[e_1,e_2,z_1,z_2]]$; recall that
all our rings contain $1=\sum_i e_i$. 
A natural basis in $\mathcal{O}/\mathcal{R}$ is formed by
(the images of) 
$\{e_1,z_1,z_2,z_1^2,z_1^3,z_1^4,z_1^5,z_1^7\}.$ We use different
signs in $x$ and $y$ to have $z_{1,2}^6$ in
$\mathcal{R}$; so are all $z_{1,2}^m$ for $m\ge 8$. 
Thus, dim$_{\mathbb{F}_q}\mathcal{O}/\mathcal{R}=8$
and 
the contribution of any {\sf\em invertible modules} $M$ (those
with $r\!k_q=1$) to
$\mathcal{H}^{mot}$ is $(q^8-q^7)t^8$. The 
difference $(q^8-q^7)$ here is because the cyclic generator
of $M$ must contain $e_1+\alpha e_2+\ldots$ for $\alpha\neq 0$
(it must be standard). 

Concerning other modules, 
$r\!k_q(\mathcal{O})=4$. 
There are $2$ more families of standard modules with $r\!k_q=4$:
$M'_{\alpha}$ generated by $e_1+\alpha e_2, z_1,z_2, z_1^2$
and $M''_{\alpha,\beta}$ generated by $e_1+\alpha e_2+\beta z_1, 
z_1-\alpha z_2, z_1^2, z_2^2$, where $\alpha\neq 0$.
They give
{\small $\bigl(1+(q-1)t+q(q-1)t^2\bigr)
(1+a q)(1+a q^2) (1+a q^3)$} in $\mathcal{H}^{mot}$.
\vskip 0.2cm

Let us discuss a bit  the corresponding contributions to the
$L$-function.
The {\sf\em conductor} of $\mathcal{R}$
is $\mathcal{C}=z^8\,\mathcal{O}$, 
the greatest $\mathcal{O}$-submodule in $\mathcal{R}$.
It has $4$ generators, $z_{1,2}^8$ and $z_{1,2}^9$,
 over $\mathcal{R}$ and dim$_{\mathbb{F}_q}\mathcal{R}/\mathcal{C}=6$.
Thus, the contribution of ideals $M_{i,j}=z_1^{i}\mathcal{C}+
z_2^{j}\mathcal{C}$ for $i,j\ge 0$  to
$Z$ is $t^6(1+a q)(1+a q^2) (1+a q^3)/(1-t)^2$.  The other
two families are based on 
$z^8 M'$ and $z^8M''$ of $q$-rank $4$.
They contribute 
$t^6(1+t (q-1)+t^2(q^2-q)(1+a q)(1+a q^2) (1+a q^3)$ to $L$. 
The expected general formula is:
$L_{6,4}(q,t,a)=\mathcal{H}^{mot}_{6,4}(qt,t,a)=$
\renewcommand{\baselinestretch}{1.2} 
{\small
\(
1-2 t+t^2+q^3 t^4-2 q^3 t^5+q^3 t^6+q^4 t^6-2 q^4 t^7+q^4 t^8
+q^5 t^8-2 q^5 t^9+q^5 t^{10}+q^6 t^{10}-2 q^6 t^{11}+q^7 t^{12}
+q^6 t^{14}-2 q^7 t^{15}+q^8 t^{16}+t \bigl(1-t+q t-t^2-q t^2+q^2 t^2
+t^3-q t^3-2 q^2 t^4+2 q^3 t^4+2 q t^5-q^3 t^5
+q^4 t^5-q t^6-3 q^3 t^6+q^4 t^6+2 q^2 t^7+q^3 t^7-q^4 t^7
+q^5 t^7-q^2 t^8-3 q^4 t^8+q^5 t^8+2 q^3 t^9-q^5 t^9
+q^6 t^9-2 q^5 t^{10}+2 q^6 t^{10}+q^4 t^{11}-q^5 t^{11}
-q^5 t^{12}-q^6 t^{12}
+q^7 t^{12}-q^6 t^{13}+q^7 t^{13}+q^7 t^{14}\bigr) (1+a q)
+t^3 \bigl(1-t+q t
+q^2 t^2-t^3-2 q t^3+t^4+2 q^3 t^4-a q^4 t^4+q t^5-3 q^2 t^5
-q^3 t^5+2 a q^4 t^5-a q^5 t^5+q^2 t^6+2 q^4 t^6-a q^4 t^6
-2 q^3 t^7+a q^5 t^7+q^5 t^8-q^4 t^9+q^5 t^9+q^5 t^{10}\bigr) 
(1+a q) (1+a q^2)
+t^6 \bigl(1-t+q t-q t^2+q^2 t^2\bigr) (1+a q) (1+a q^2) (1+a q^3).
\)
}
\renewcommand{\baselinestretch}{1.2} 

The coincidence of $L_{6,4}$ with $\mathbf{H}_{6,4}$ 
was  verified only partially. Recall that $L(q,t,-1/q)$
is conjecturally the {\sf\em Z\'u{\oldt{n}}iga 
$L$-function}. 



\subsection{\bf Quasi-rho-invariants} \label{sec:rho}
The $\rho_{ab}-$ invariant
is the von Neumann invariant defined for the abelianization
representation $\pi_1(S^3\setminus K)\to \Z$. We will define 
a superpolynomial for a certain {\sf\em integer} variant
of $\rho_{ab}$ for algebraic knots $K$. The superduality 
$q^{\de}t^{2\de} \H_K(q,\frac{1}{qt},a)=\H_K(q,t,a)$
and the connection conjecture for $a\to -1/q$ will be used;
the dependence on the knot $K$ will be shown. 
We set:
\begin{align*}
&R_K(q,t,a)\equal
\bigl(\H_K(q,t,a)-t^{\de}\H_K(q,t\!=\!1,a)\bigr)/
\bigl((1-qt)(1-t)\bigr),\\
&\rho_K(q,t)\equal R_K(q,t,a\!=\!-\!1/q),\ \text{where}\,\ 
\H_K(q,1,a\!=\!-1/q)=q^{\de}. 
\end{align*} 
Switching to 
$L_{prncpl}$, our $\rho_K(q,t)$ is a sum of monic $q^i t^j$ 
(of multiplicity one); see the next section for the
exact formula.
The superduality holds
$q^{\de-1} t^{2\de-2}R_K(q,\frac{1}{qt},a)=R_K(q,t,a)$; 
the same holds for $\rho_K$. 

\vskip 0.2cm

Let us express $\rho_K(1,1)\!=\!\rho_K(q\!=\!1,t\!=\!1)$ 
in terms of 
$\Ga\!=\! \nu_z(\r\setminus\{0\})$. We set\, $G\equal\Z_+\setminus \Ga=
S$ for 
$S=\cup_{i=1}^{\varpi}\, [g_i,g'_i]$, a disjoint union of 
segments,   where
$g_i\le g_i'\in G\not\ni g'_i+1$. Then $\de\!=\!
\sum_{i=1}^\varpi m_i$ for $m_i\!\equal\!g'_i\!-\!g_i+1$.
This is actually for any Gorenstein $\r\subset \C[[z]]$.
 Setting $\varsigma(x)=x$ for $x\in S$ and $0$
otherwise,
\vskip 0.2cm

\centerline{
$
\rho_K(1,1)=\sum_{i=1}^{\varpi}m_i(g_i'+1-\frac{m_i}{2})-\frac{\de^2}{2}=
\int_0^\infty \varsigma(x)dx-\frac{\de^2}{2}.
$
}
\vskip 0.2cm

Geometrically,
 $\rho_{ab}=\int_0^1 \si_K(e^{2\pi i x})dx$ for
the {\sf\em Tristram-Levine signature} $\si_K$; see e.g.
{\sf\em ``Signatures of iterated torus knots"}
(Litherland, 1979). Our $\varsigma$ is some variant 
of $\si_K$; taking $\int_0^1$ in the formula for $\rho_{ab}$,
makes  $\rho_{ab}$
``additive" for iterated knots (see below).

\vskip 0.2cm

{\sf\em Quasi-rho for cables.} For $r,s>0$ 
such that gcd$(r,s)=1$,
one has:
$
\rho_{r,s}(1,1)=\frac{(r^2-1)(s^2-1)}{24}.
$
We note that this is the maximal size of the 
{\sf\em $(r,s)$-core partitions}.                
The classical $\rho_{ab}$ is 
$-\frac{1}{3}\frac{(r^2-1)(s^2-1)}{rs}$.
Thus, we basically obtain the same formula up to some
renormalization. The values of our  $\rho$ are always 
natural numbers; so it 
can be used for {\sf\em categorification}. 
This is not the case with $\rho_{ab}$ due to its 
geometric origin. Our  $\rho$ is expected to have some
natural geometric interpretation too.

\vskip 0.2cm

The classical $\rho_{ab}$ 
for cables is known to be 
additive. For instance let $K=C\!ab(m,n)C\!ab(s,r)$.
Then 
$\rho_{ab}=-\frac{1}{3}\bigl(\frac{(m^2-1)(n^2-1)}{mn}
+\frac{(r^2-1)(s^2-1)}{rs}\bigr)$. This deviates 
from our $\rho$, which is additive with some weights.

\vskip 0.2cm

Let $\r=\F[[z^{\upsilon r}, z^{\upsilon s} +z^{\upsilon s+p}]]$,
where, gcd$(r,s)=1$ as above, $\upsilon>1$ and gcd$(\upsilon,p)=1$ for
$p\ge 1$. Then $\Ga=\lan \upsilon r,\,
\upsilon s,\,\upsilon r s\!+\!p\ran$,\ 
$2\de=\upsilon^2 r s\!-\!\upsilon (r\!+\!s)+(\upsilon\!-\!1)p+1$ and 
$K=C\!ab(m\!=\!\upsilon r s\!+\!p,\,
n\!=\!\upsilon)\,
C\!ab(s,r)$. One obtains:  $\rho_K(1,1)=$
{\small $\frac{1}{24}\bigl((m^2-1)(n^2-1)
+\upsilon^2(r^2-1)(s^2-1)\bigr)$}. 
\vskip 0.1cm

More generally,
$\rho_K(1,1)=$
{\small $\frac{1}{24}\sum_{i=1}^k\upsilon_i^2 (a_i^2-1)(r_i^2-1)$}
for the cable $K=C\!ab(a_k,r_k)\cdots C\!ab(a_2,r_2)C\!ab(a_1,r_1)$,
where  $1\le i\le k$, $\upsilon_i=r_k\cdots r_{i+1}$ and
$\upsilon_{k}=1$. We will post the
details elsewhere. Here
$\upsilon_i\!=$\,gcd$(u_1,\cdots,u_{i+1})$ for 
$\Ga=\lan u_1,\cdots,u_{k+1}\ran$, where  $u_i<u_{i+1}$ and
$\upsilon_{i+1}|\upsilon_i$; it is known that
$\de=$
{\small $\frac{1}{2}\sum_{i=1}^k\upsilon_i (a_i-1)(r_i-1)$}.

\vskip 0.1cm

We note the following natural embedding for $K'=C\!ab(a,r)K$
(if both are algebraic knots).
If $\de$ is that for the ring $\r$ of $K$, then 
 $\rho_{K}(q,t)$ is the sum of monomials $q^it^j$ in
$\rho_{K'}(q,t)$ such that $j<2\de-1$.
\vskip 0.2cm

{\sf\em The case of C\!ab(13,2)C\!ab(2,3)}.
Here $\r\!=\!\C[[z^4,z^6\!+\!z^7]]$,
 $r\!=\!3,s\!=\!2,\upsilon\!=\!2$, $\de\!=\!8$. Then
$\rho(1,1)=25$ and  its refined
version is $\rho(q,t)=$ 
\renewcommand{\baselinestretch}{1.2} 
{\small
\(
1+q t+q^2 t^2+q^3 t^3+q^3 t^4+q^4 t^4+q^4 t^5+q^5 t^5+q^4 t^6
+q^5 t^6+q^6 t^6+q^5 t^7+q^6 t^7+q^7 t^7+q^5 t^8+q^6 t^8+q^7 t^8
+q^6 t^9+q^7 t^9+q^6 t^{10}+q^7 t^{10}
+q^7 t^{11}+q^7 t^{12}+q^7 t^{13}+q^7 t^{14}.
\)
}
\renewcommand{\baselinestretch}{1.2}

\vskip 0.1cm
RH holds for 
$\rho(q,t)$ when  $q\!<\!q_{sup}\!\approx\! 0.802$. 
Presumably, $\lim_{p\to\infty}q_{sup}\!=\!1$ 
for $C\!ab(2p\!+\!13,2)C\!ab(2,3)$;
for instance, $q_{sup}\!\approx\! 0.996$ for $p\!=\!2000$.

\vskip 0.2cm

Let us provide now full 
 $R(q,t,a)$ for  $\C[[z^4,z^6\!+\!z^7]]$. It is:

\renewcommand{\baselinestretch}{1.} 
{\footnotesize
\( 
1+t+q t+t^2+2 q t^2+q^2 t^2+t^3+2 q t^3+3 q^2 t^3+q^3 t^3+t^4
+2 q t^4+4 q^2 t^4+4 q^3 t^4+q^4 t^4+t^5+2 q t^5+4 q^2 t^5
+6 q^3 t^5+4 q^4 t^5+q^5 t^5+t^6+2 q t^6+4 q^2 t^6+7 q^3 t^6
+8 q^4 t^6+4 q^5 t^6+q^6 t^6+t^7+2 q t^7+4 q^2 t^7+7 q^3 t^7
+10 q^4 t^7+8 q^5 t^7+4 q^6 t^7+q^7 t^7+q t^8+2 q^2 t^8+4 q^3 t^8
+7 q^4 t^8+8 q^5 t^8+4 q^6 t^8+q^7 t^8+q^2 t^9+2 q^3 t^9
+4 q^4 t^9+6 q^5 t^9+4 q^6 t^9+q^7 t^9+q^3 t^{10}+2 q^4 t^{10}
+4 q^5 t^{10}+4 q^6 t^{10}+q^7 t^{10}+q^4 t^{11}+2 q^5 t^{11}
+3 q^6 t^{11}+q^7 t^{11}+q^5 t^{12}+2 q^6 t^{12}+q^7 t^{12}
+q^6 t^{13}+q^7 t^{13}+q^7 t^{14}+a \Bigl(q t+q t^2
+2 q^2 t^2+q t^3+3 q^2 t^3+3 q^3 t^3+q t^4+3 q^2 t^4+6 q^3 t^4
+3 q^4 t^4+q t^5+3 q^2 t^5+7 q^3 t^5+9 q^4 t^5+3 q^5 t^5
+q t^6+3 q^2 t^6+7 q^3 t^6+12 q^4 t^6+10 q^5 t^6+3 q^6 t^6
+q t^7+3 q^2 t^7+7 q^3 t^7+13 q^4 t^7+17 q^5 t^7+10 q^6 t^7
+3 q^7 t^7+q^2 t^8+3 q^3 t^8+7 q^4 t^8+12 q^5 t^8+10 q^6 t^8
+3 q^7 t^8+q^3 t^9+3 q^4 t^9+7 q^5 t^9+9 q^6 t^9+3 q^7 t^9
+q^4 t^{10}+3 q^5 t^{10}+6 q^6 t^{10}+3 q^7 t^{10}+q^5 t^{11}
+3 q^6 t^{11}
+3 q^7 t^{11}+q^6 t^{12}+2 q^7 t^{12}+q^7 t^{13}\Bigr)
+a^2 \Bigl(q^3 t^3+q^3 t^4+2 q^4 t^4
+q^3 t^5+3 q^4 t^5+3 q^5 t^5+q^3 t^6+3 q^4 t^6+6 q^5 t^6+3 q^6 t^6
+q^3 t^7+3 q^4 t^7+7 q^5 t^7+8 q^6 t^7+3 q^7 t^7+q^4 t^8
+3 q^5 t^8+6 q^6 t^8+3 q^7 t^8+q^5 t^9+3 q^6 t^9+3 q^7 t^9
+q^6 t^{10}+2 q^7 t^{10}+q^7 t^{11}\Bigr)
+a^3 \Bigl(q^6 t^6+q^6 t^7
+q^7 t^7+q^7 t^8\Bigr).
\)
}
\renewcommand{\baselinestretch}{1.2}

\vskip 0.1cm
Recall that $R(q,t,a)\mapsto\rho(q,t)$ upon the substitution 
$a\mapsto -\frac{1}{q}$  in  the parameters of  $\H(q,t,a)$, which is
generally the passage to the Heegaard-Floer homology and 
Alexander polynomials (when $q\!=\!1,a\!=\!-1$).

\section[\sc \hspace{1em}On physics connections] 
{\sc On physics connections} 
Generally, a challenge is to associate 
the Riemann and Lindel\"of hypotheses 
with some 
physics phenomena in SCFT
or similar theories. SCFT is connected
with quite a few recent  mathematical 
developments. DAHA can be considered as its part;
their origin was in the  
Knizhnik-Zamolodchikov equations.
DAHA superpolynomials
can be interpreted as some physics partition functions,
those for knot operators. 
\vskip 0.2cm

The 
{\sf\em $p$-adic strings} due to Witten and others must be 
mentioned in this context.
The starting point of this theory was an adelic product formula for the 
{\sf\em Veneziano amplitude}. In mathematics, we have a long history
of understanding ``geometry" via $p$-adic constructions and those
over $\F_q$, including recent Peter Scholze theory (say, in his 
lectures at IHES, 2023-24). 
The figures below are about such and similar
connections. For instance, we expect that $L$-functions of plane
curve singularities over $\F_q$ can be related to some
Dirac operators; this can be connected with
$p$-adic strings.

The {\sf\em Lee-Yang circle theorems} provide 
a different perspective. The Ising model with an external 
magnetic field is the key example. 
  

\subsection{\bf Lee-Yang theorem}
For any lattice (of any dimension) with $N$ vertices and
the connected pairs 
of vertices denoted by
$\lan n,n'\ran$, let $\z=\lim_{N\to\infty}\frac{\log(Z_N)}{N}$
for the partition function
$Z_N\!=\!\sum_{\{\si_n\}}e^{-\be \h}$, where the Hamiltonian is
$\h=-\sum_{\lan n,n'\ran}J_{n,n'}\si_n\si_{n'}-H\sum_n\si_n$
and $\si=\pm 1$. This is the Ising model with an
external magnetic field $H$. 
Here $\be=(k_B T)^{-1}$ is the inverse temperature for
the Boltzmann constant $k_B$.
Assuming that $J_{n,n'}\ge 0$ (the ferromagnetic case)
and $\be>0$, Lee-Yang proved
that the zeros of $Z_N$ in terms of the ``complex fugacity"
$\mu=e^{-2\be H}$ 
belong to the unit circle $|\mu|=1$; the
corresponding symmetry of $Z_N$  is simply $\si\mapsto -\si$. 
 For the square lattice with 
$J=const>0$, $Z_N$ is a polynomial in 
terms of $\mu$ and 
$0<u\equal e^{-4\be J}<1$. There is a $q$-version of this theorem
and other physics-statistical variants. 
\vskip 0.1cm

The Lee-Yang-Fisher zeros are 
when $u$ is considered as a free parameter for  
complex $T$.
Numerical experiments 
showed that $|\mu|=1$ for
the $\mu$-zeros can hold for some  $u<0$. 
The physics calculations are mostly when
$\mu$ is fixed and the $u$-zeros are considered, but
they can be used for the $\mu$-zeros too.
This phenomenon resembles 
the behavior
of the  $t$-zeros of our $\H(q,t,a)$. 
Actually, DAHA is
directly related to the $XXZ$-model, which is somewhat similar 
to the Ising model with $H$ as above, though all attempts to
``integrate" the latter failed. 
\vskip 0.2cm

Only $\z$  is physical; its phase transitions 
are positive {\sf\em real} limits as $N\to\infty$  of (complex) 
$\mu$-zeros 
of $Z_N$. Thus, these zeros can result in a phase
transition  only at $\mu=1$ due
to  RH for  $Z_N$, which point is the intersection of the unit
$\mu$-circle with $\R_+$. 
The relation between the failure of
RH and ``unwanted" phase transitions seems sufficiently general. 
Given $u<0$, the $\mu$-zeros 
of $Z_N$ quickly
become  wild near the real line 
when $N$ goes beyond $N_u$, the 
last $N$ when RH still holds for $u<0$.  So do the corresponding
points
of phase transition.  
This can be clearly seen in the 1D Ising model.

\vskip 0.2cm

{\sf\em One-dimensional case.} The zeros
of $Z_N$ can be found explicitly in this case. One has:
$-\be\h=\be\sum_{n=1}^N J \si_n\si_{n+1}+\frac{\be}{2}
H(\si_n+\si_{n+1})$,
where the periodicity $\si_{N+1}\!=\!\si_1$ is assumed.
Then $Z_N\!=\!\sum_{\{\si_n\}}e^{-\be \h}$ can be
calculated using the eigenvalues $\la_1,\la_2$
of the {\sf\em transfer matrix}\, 
{\small $\t=\begin{pmatrix} e^{\be(J+H)} & e^{-\be J}\\
e^{-\be J} & e^{\be(J-H)}\end{pmatrix}$.} Namely,
 $Z_N=$\ tr\,$(\t^N)=\la_1^N+\la_2^N$ for
$\la_{1,2}=e^{\be J}\cosh(\be H)\pm \sqrt{e^{2\be J} \sinh^2(\be H)+
e^{-2\be J}\,}$. Upon some algebraic manipulations,
the $\mu$-zeros of $Z_N$ are for $H=\imath\, \theta_n/\be$, where
$\cos(\th_n)=\sqrt{1-u\,}\,\cos\frac{(2n-1)\pi}{2N}$
\,for  $n=1,\ldots, N$ and $u=e^{-4\be J}$ as above.
We obtain that RH, which is the condition $\th_n\in \R$ for
any $n$,  holds 
if $\sqrt{1-u\,}\,\cos\frac{\pi}{2N}\le 1$, i.e. for $u\ge
-\tan^2\left( \frac{\pi}{2N}\right)$. Thus, RH can hold
for negative $u$, but the latter bound tends to $0$ when $N\to \infty$.

This example 
provides some physical insight.  A counterpart of $N$ is
our deg$_{\,t}\,\h^{\la}$, which is 
$\de\sum_i{m_i^2}$ for $\la=(m_i)$. For instance, 
 deg$_{\,t}\,\h^{\la}=\de n$ as $a=0$ for $n$-columns $\la.$ 
We note that the distribution of $\mu$-zeros above
resembles our one for $\H(q,t,a\!=\!0)$
in the case of uncolored $T(2,2p\!+\!1)$.

\subsection{\bf Landau-Ginzburg models}
Another physics approach to singularities 
is presented in  paper  
{\sf\em ``Catastrophes and the classification of conformal theories"} 
(Vafa-Warner, 1989). The authors consider LGSM,
Landau-Ginzburg Sigma Models,  for {\sf\em superpotentials}\,
$W(x,y)$\, corresponding to isolated singularities. They can be
with several variables, and more than one superpotential $W$ can
be considered. Let us mention here
two publications by Alexander Zamolodchikov in 1986. 
We mention that there are many other classes of 
superpotentials, for instance those for quiver 
varieties and KZ. 
The correspondence between SCFT and LGSM is one of the key 
in string theory.
\vskip 0.2cm

 A lot of information
can be obtained 
directly in terms of $W(x,y)$ and the corresponding
singularities. 
For instance, the {\sf\em Milnor number} $\mu=2\de$,
which is $(r\!-\!1)(s\!-\!1)$
for $W_{r,s}(x,y)\!=\!x^s\!-\!y^r$, 
coincides with the number 
of superfields. It is
is the {\sf\em Witten index} for plane curve singularities:\, the number
of zero energy bosonic vacuum states minus 
the number of zero energy fermionic vacuum states. 
The dimensions of superfields for $W_{r,s}(x,y)$ are proportional to
the corresponding quasi-homogeneity weights, which are $1/s$ for $x$,
$1/r$ for $y$ and so on.  

Another example is the central charge, which is  
$c\!=\!6\be$ for $\be=(\frac{1}{2}\!-\!\frac{1}{r})
(\frac{1}{2}\!-\!\frac{1}{s})$ for $W_{r,s}(x,y)$.
Generally, $\be$ is obtained from the asymptotic
formula $\int e^{\,\imath\,\la\, W(x_1,\ldots, x_m)}
\prod_{i=1}^m \la^{1/2}dx_i \sim
O(\la^{\be})$ for large $\la$. Also, the adjacency of singularities
plays an important role in this approach. 

\subsection{\bf Refined Witten index}
Refined  Witten and BPS indices were studied in the
literature; see e.g.
Gaiotto-Moore-Neitzke (Adv. Theor. Math. Phys. 2013). 
Generally, the challenge is
to ``split" the vacuum states counted by these indices
using some additional parameters.

Let us split $\mu=2\de$ using $\h^{mot}(q,t,a=-t/q)$, 
though we will actually use $L_{prncpl}$ in the following
calculation. We begin with  
$\de_{q,t}\!\equal\!$
{\large $\frac{\h^{mot}(q,t,a=-t/q)-(qt)^\de}{1-t}
=\frac{L_{pncpl}(\frac{q}{t},t)-(qt)^\de}{1-t}$}.
This formula is for any Gorenstein $\r\subset \C[[z]]$. The 
DAHA parameters $q,t$ from $\h$ are used (not $q_{new}$
from the definition of $\H$ and $L$).

\vskip 0.1cm

As above:
$G=\Z_+\!\!\setminus\! \Ga=$
$\cup_{i=1}^{\varpi}\, \{g_i\le x \le g'_i\}$, where
$g'_i+1\in \Ga$,
and $m_i\!=\!g'_i\!-\!g_i\!+\!1$. Then 
$
\de_{q,t}\!=\!$ {\large $\frac{1\!-\!t^{g_1}}{1\!-\!t}+
\sum_{i=1}^{\varpi-1}\frac{t^{g_i'+1}-t^{g_{i+1}}}{1-t}
(\frac{q}{t})^{m_1\!+\cdots+\!
m_i}$}, and $\de_{1,1}=\de$.  This formula was actually used
above in the definition of the refined quasi-rho invariant 
$\rho_K(q,t)$; also, see below.

\vskip 0.1cm

Let 
$\mu(q,t)\equal\de_{q,t}+(qt)^{\de-1} \de_{t^{-1},q^{-1}}=$
$\sum_{x=0}^{2\de-1} t^{v(x)-1}q^{g(x)}$, where
$v(x)=|\{\nu\!\in\! \Ga \mid 0\!\le\! \nu\!\le\! x\}|$ and 
$g(x)=|\{g\!\in\! G \mid 0\!\le\! g\! <\! x\}|$ for $G$ as above.
Then  $\mu(1,1)=\mu$ and
this definition ensures the superduality:
 $(qt)^{\de-1}\mu(1/t,1/q)=\mu(q,t)$.
Not all monomials are monic in $\mu(q,t)$: 
$\varpi$ of them are with coefficient $2$, which
correspond to $x\in \Ga \not\ni x+1$.

\vskip 0.1cm

For example,
$\mu(q,t)\!=${\small 
$2 + q + q^2 + 2 q^3 t + 2 q^4 t^2 + 2 q^5 t^3 + 2 q^6 t^4 + q^7 t^5 + 
 q^7 t^6 + 2 q^7 t^7$}
for $\r=\C[[z^4,z^6\!+\!z^7]].$
Upon $q\mapsto qt$,
it satisfies RH for $0\!<\!q\!<\!0.919090$. 
For $\r=\C[[z^6,z^9+z^{460}]]$,
this range becomes $0\!<\!q\!<\!0.852561$.
For $\r=\C[[z^6,z^8\!+\!z^{649}]]$, it is
$0\!<\!q\!<\!0.846566$ and  $0.848063$ for $z^8\!+\!z^{3003}$.
\vskip 0.2cm

Compare with the formula for 
$\varrho(q,t)\equal \frac{\h(q,t,a=-t/q)-q^{\de}}{(1-t)(1-q)}=
\rho(\frac{q}{t}, t)$: 

\centerline{
$\varrho(q,t)=\sum_{x\in G}q^{g(x)}\frac{1\!-\!t^{v(x)}}{1-t}\!=\!
\!\sum_{G\ni x>y\in \Ga}q^{g(x)}t^{v(y)-1}=(q t)^{\de-1}
\varrho(t^{-1},q^{-1}).
$
}

\noindent
Recall that we use the DAHA parameters $q,t$ from $\h$ in this section.
\vskip 0.2cm
\vfil

{\sf\em Adding colors.} The substitution $a\mapsto -\frac{t}{q}$ 
has remarkable properties 
for $\h^{\la}(q,t,a)$ for partitions $\la$ more general
than $\yng(1)$ (the uncolored case). 
Let $n$ be the number of rows of $\la$ and $m$ the number of its
columns.

 For {\sf\em hooks} $\la$, we expect that 
 $\h^{\la}(q,t,a\mapsto -\frac{t}{q})=(q t^{n-1})^{\de(n-1)}r^{\la}(q,t)$
for $ r^{\la}(q,t)\!=\!r(q\!\mapsto\! t^{n-1} q^m,
t\!\mapsto\! q^{m-1} t^n)$. Here $r(q,t)\!=\!
\h(q,t,a\!\mapsto\! -\frac{t}{q})$ is for $\la=\yng(1)$\,
considered above; its constant term is $1$. For instance,
$1\!-\!r^{\la}(q,t)$ is divisible by
$(1-q^m)$ for $n\!=\!1$. \Yboxdim5pt
This was checked for pure
columns/rows and several hooks with $m\!=\!2$ or $n\!=\!2$.
 For example,
{\small $\bigl(1-r^{\yng(2)}(q,t)\bigr)/(1-q^2)=qt 
(1 + q^2 + q^4 + q^7 t + q^{10} t^2 + q^{13} t^3 + q^{16} t^4 + 
   q^{21} t^7)$} for $K=C\!ab(13,2)C\!ab(2,3)$.


\vskip 0.1cm
This is
far from being that simple beyond the hooks.
For instance, for $T(3,2)$ and $\la=2\times 2=\yng(2,2)$\,:
{\small $1-r_{3,2}^{2\times 2}(q,t)=
q (1 - q t) (1 + q - q^2 + t - q^2 t + q^4 t + q^3 t^2 - q^5 t^2 + 
   q^2 t^3 - q^4 t^3 + q^5 t^4 + q^6 t^4 + q^4 t^5 - q^6 t^5 - 
   q^6 t^6 + q^6 t^7)$}.

\vskip 0.2cm
The definition of the $\varrho^\la(q,t)$ for {\sf\em symmetric} $\la$
can follow the uncolored case, but this is preliminary.
\Yboxdim6pt
 For instance, one can set:
$\varrho_{K}^{2\times 2}(q,t)\equal$ 
{\small 
$\bigl(r_{K}^{2\times 2}(q,t)-q^{4\de} t^{4\de}\bigr)/(1-qt)^2$}.
For the example above: \,$\varrho_{3,2}^{2\times 2}(q,t)=$  {\small
$1 - q - q^2 + q^3 + q t - q^2 t + q^4 t - q^5 t + 
   2 q^2 t^2 - q^3 t^2 - q^4 t^2 + q^5 t^2 + 2 q^3 t^3 - q^4 t^3 + 
   q^6 t^3 + 2 q^4 t^4 - q^5 t^4 - q^6 t^4 + q^5 t^5 - q^6 t^5 + 
   q^6 t^6
$}\,, which satisfies the superduality $q\leftrightarrow t^{-1}$
with the multiplier $q^6 t^6$. Presumably, $\varrho_K^{2\times 2}(1,1)
=4\varrho_K(1,1)$ for algebraic knots $K$, but there are negative
terms in $\varrho_{3,2}^{2\times 2}(q,t).$ Recall that
we omit $\la$ in uncolored 
$\h_K,\varrho_K$ and so on, i.e. for $\la=\yng(1)$\,.


\subsection{\bf {\em S}-duality}
The relation  SCFT $\leftrightsquigarrow$ LGSM suggests that 
the {\sf\em $S$-duality} in the former
 can be seen via the superpotential $W(x,y)$.
The superduality of physics 
superpolynomials can be connected with that in {\sf\em $M$-theory} and
the symmetry $\ep_1\leftrightarrow \ep_2$ in {\sf\em Nekrasov's
instanton sums}. The general physics  superduality (with $\la$) 
for superpolynomials was considered by
Gukov-Stosic (2012) (and in some prior works).
 For us, this correspondence
is between DAHA, a ``representative" of SCFT, 
and $L$-functions of singularities, which presumably
``represent"  LGSM.

More specifically, we focus on the DAHA 
superpolynomials, which are certain partition functions 
of knot operators and are related to the BPS states.  They can be 
expected to coincide with motivic superpolynomials (a solid mathematical
conjecture),  which are
presumably some partition functions of properly defined
{\sf\em motivic} LGSM for plane curve
singularities. Then, switching to 
the $L$-functions of the latter, the 
$S$-duality from SCFT becomes the Hasse-Weil functional equation,
with some potential toward various generalizations.

\vskip 0.2cm

This link  may be not too much surprising.
The $S$-duality and mirror symmetry (CPT) are very universal in physics.
The functional equation is certainly   
of the same calibre in mathematics. 
There are
more than $20$ different zeta-theories; RH does not always holds. 
By analogy with the Lee-Yang theorem,
one can speculate that (topological) LGSM associated with plane
curve singularities are 
``stable" when the ``coupling constant" $q$ is small
enough to ensure RH for $\H(q,t,a)$. More generally, can 
the failure of RH be somehow connected with the presence
of unwanted phase transitions?  
Something like this:  ``elementary particles" for
$W_{2,2p+1}=x^{2p+1}-y^2$ can be ``observed" for any 
$q\!>\!0$, but this requires $q<1/2$ or so for arbitrary
superpotentials $W(x,y)$.

\vskip 0.2cm

We took the  bound $1/2$ from {\sf\em strong RH}
for  $a\!=\!0$; generally, we can
make $a$ a constant or any quantity invariant under
the superduality.  For instance, let 
$a= -t/q$, which is super-invariant.
Then  we arrive
at 
$\varrho_{2,2p+1}(q,t)=\frac{q^p -t q^{p+1}-t (q t)^p+
(q t)^{p+1}+t-1}{(q-1) (t-1) (q t-1)}=
\sum_{0\le j\le i<p} q^it^j$, and  RH holds for $\rho_{2,2p+1}(q,t)=
\varrho_{2,2p+1}(qt,t)$ only for sufficiently small $q$ (not all).

\vskip 0.2cm

Let me finish this section with little something on Manin's
``{\sf\em Mathematics as
metaphor".} Yu.I. obviously expected number theory to play 
a major role in the alliance of physics and mathematics. 
If RH has something to
do with the absence of unwanted phase transitions in physics 
theories or their stability of any kind, then number theory
will not be just a ``metaphor".  Technically, DAHA
accumulated quite a few integrable models and the fact that
it appeared  very  ``motivic" can be meaningful physically. 
I thank my friends-physicists for various
talks on these matters (though they are not responsible for
what I wrote).

\section[\sc\hspace{1em}Zeta-functions  as invariants] 
{\sc Zeta-functions  as invariants} 
\subsection{\bf The first figure}
Modern mathematics and physics very much rely (more than ever)
on the progress in geometry. Any new geometric approaches to the classical
zeta and $L$-functions can open new avenues toward
the justification of ``Grand Conjectures". 
This of course does not diminish the role of analytic methods,
including classical Fourier analysis and (more recently)
the $p$-adic methods. The Dwork proof of the 
rationality of zeta-functions is an example of the latter.

\begin{figure*}[htbp]
\hskip -0.in
\includegraphics[scale=0.55]{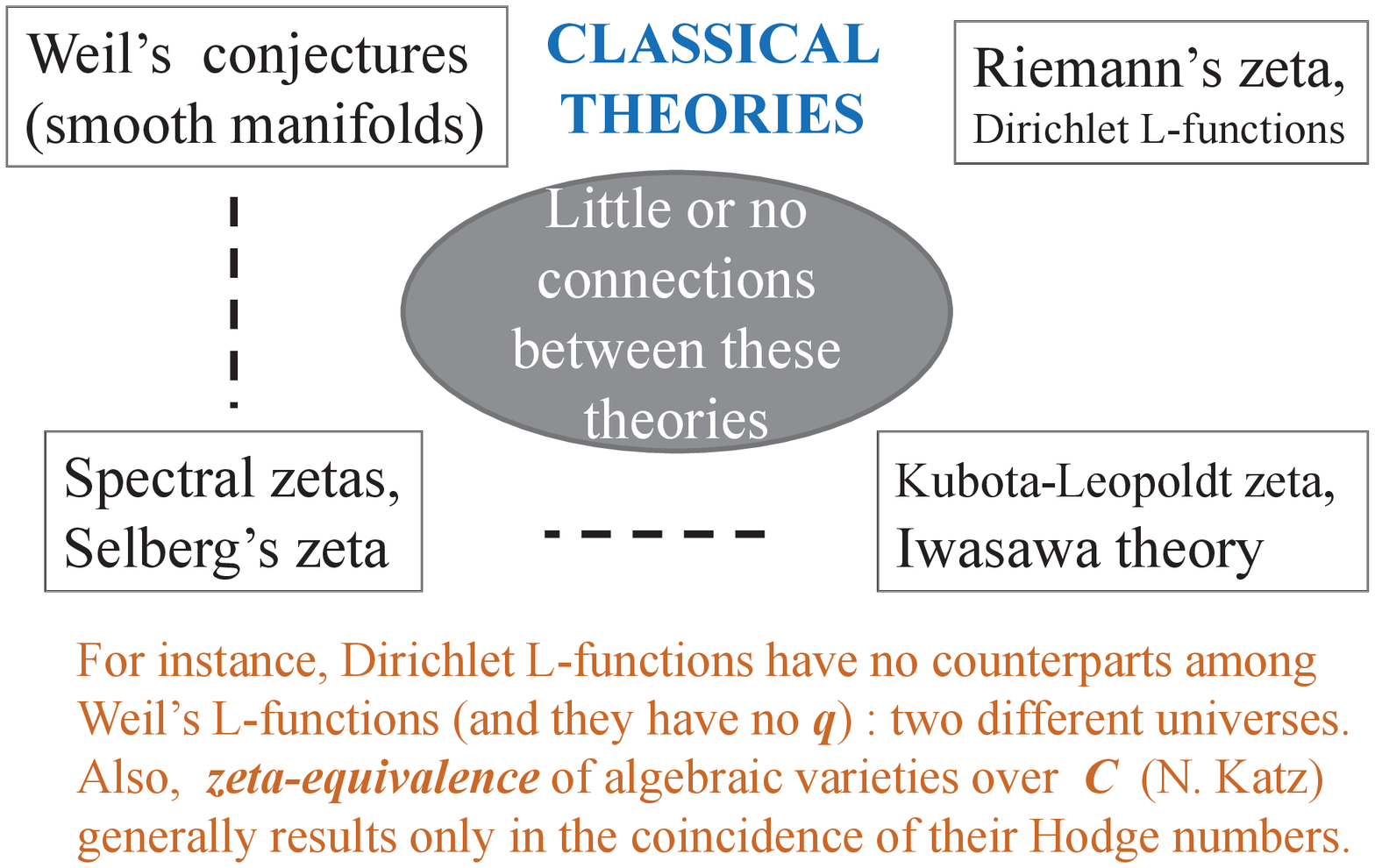}
\end{figure*}

In contrast to the Weil conjectures (proved by Deligne),
the  distribution of the zeros  
of Riemann's $\ze(s)$ in the critical strip 
does not seem  to reflect any ``geometry".  Riemann's zeta
does occur
in some geometric and physics considerations, but the
interpretation of the Grand Conjectures geometrically or
physically is (totally?) missing. 
We note here that the zeros of 
Selberg's zeta functions are much more ``geometric"
and have many analytic applications. The zeros of
Hasse-Weil zetas are very geometric too, which was the
key for the justification of the corresponding RH
in special cases and in general (Deligne).    
\vskip 0.2cm
 
Nicolas Katz proved that the zeta-equivalence
of algebraic varieties $X$ results in the coincidence of their
virtual Hodge numbers. One needs to add
the coefficients of the equations of $X$ to the corresponding
rings of functions and then consider the zeta function $\ze_X$
of the resulting scheme over $\Z.$ If $\ze_X=\ze_Y,$ then
we call $X$ and $Y$ zeta-equivalent. The coincidence
of Hodge numbers is of course very far from the existence 
of any kind of isomorphism between $X$ and $Y$. 

\vskip 0.2cm

\comment{
Generally, {\sf\em Kapranov's (motivic) zeta} 
and ``true motivic superpolynomials"
are with the coefficients in the 
Grothendieck ring $K_0(V\!ar/\F)$
of varieties over $\F$. Recall that we conjecture that
the Piontkowski cells are configurations of affine spaces
for plane curve singularities.
The map  $X\mapsto |X(\F_q)|$
is one of the {\sf\em motivic measures}. 
Let me mention here {\sf\em motives} and
Manin's {\sf\em ``Lectures on zeta functions 
and motives"} (1995).

\vskip 0.2cm

In the first figure, the following 4 theories
are presented as disconnected blocks: 
(a) the Hasse-Weil zetas (over $\F_q$), (b) the
Selberg's zetas (via the Laplace operators),
(c) the Kubota-Leopoldt zetas, and
(d) Riemann's $\ze(s)$ and the Dirichlet
$L$-functions. They are really different,
but not totally disconnected; there are some
deep relations. 

For instance, let us mention 
the connection (not reflected in this figure) of the 
Selberg's {\sf\em ``1/4 conjecture"} for arithmetic subgroups 
$\Ga\subset SL(2,\R)$ to the Riemann's zeta. 

The upper-left block is connected with the upper-right one:
global fields and the product formulas are here and there.
However this
does not help much to understand the zeros of 
the classical $\zeta(s),L(s)$. The main difference is obvious:
Riemann's $\zeta(s)$ 
does not contain $q$. One needs, at least, its $q$-deformation.
In the theory over $\F_q$, making $q$ a ``continuous parameter"
is necessary for any links. This can be hopefully 
achieved via the passage to zeta-functions of certain families of
isolated surface singularities. 
Let us try to outline this.
}

Generally, {\sf\em Kapranov's (motivic) zetas} 
are with the coefficients in the
Grothendieck ring $K_0(V\!ar/\mathbb{F})$
of varieties over $\mathbb{F}$. The map  
$X\mapsto |X(\mathbb{F}_q)|$
is a {\sf\em motivic measures}.
Conjecturally, our motivic superpolynomials
can be lifted to $K_0(V\!ar/\mathbb{F})$. 
Recall that we conjecture that 
the Piontkowski cells are configurations of affine spaces.
Let me mention 
Manin's {\sf\em "Lectures on zeta functions
and motives"} (1995).
\medskip

In the first figure, the following 4 theories
are presented as disconnected blocks:
(a) the Hasse-Weil zetas (over $\mathbb{F}_q$), (b) the
Selberg's zetas (via the Laplace operators),
(c) the Kubota-Leopoldt zetas, and
(d) Riemann's $\zeta(s)$ and the Dirichlet
$L$-functions. They are really different,
but not totally disconnected; there are some
deep relations.
For instance, 
the connection (not in this figure) of the
Selberg's {\sf\em "1/4 conjecture"} for arithmetic subgroups
$\Gamma\subset SL(2,\R)$ to the Riemann's zeta is such.

The upper-left block is connected with the upper-right one:
global fields and the product formulas are here and there.
However this
does not help much to understand the zeros of 
the classical $\zeta(s),L(s)$. The main difference is obvious:
Riemann's  $\zeta(s)$ 
does not contain $q$. One needs, at least, its $q$-deformation.
In the theory over $\mathbb{F}_q$, making $q$ a ``continuous parameter"
is necessary for any links. This can be hopefully 
achieved via the passage to zeta-functions of certain
isolated singularities. 

\begin{figure*}[htbp]
\hskip -0.in
\includegraphics[scale=0.55]{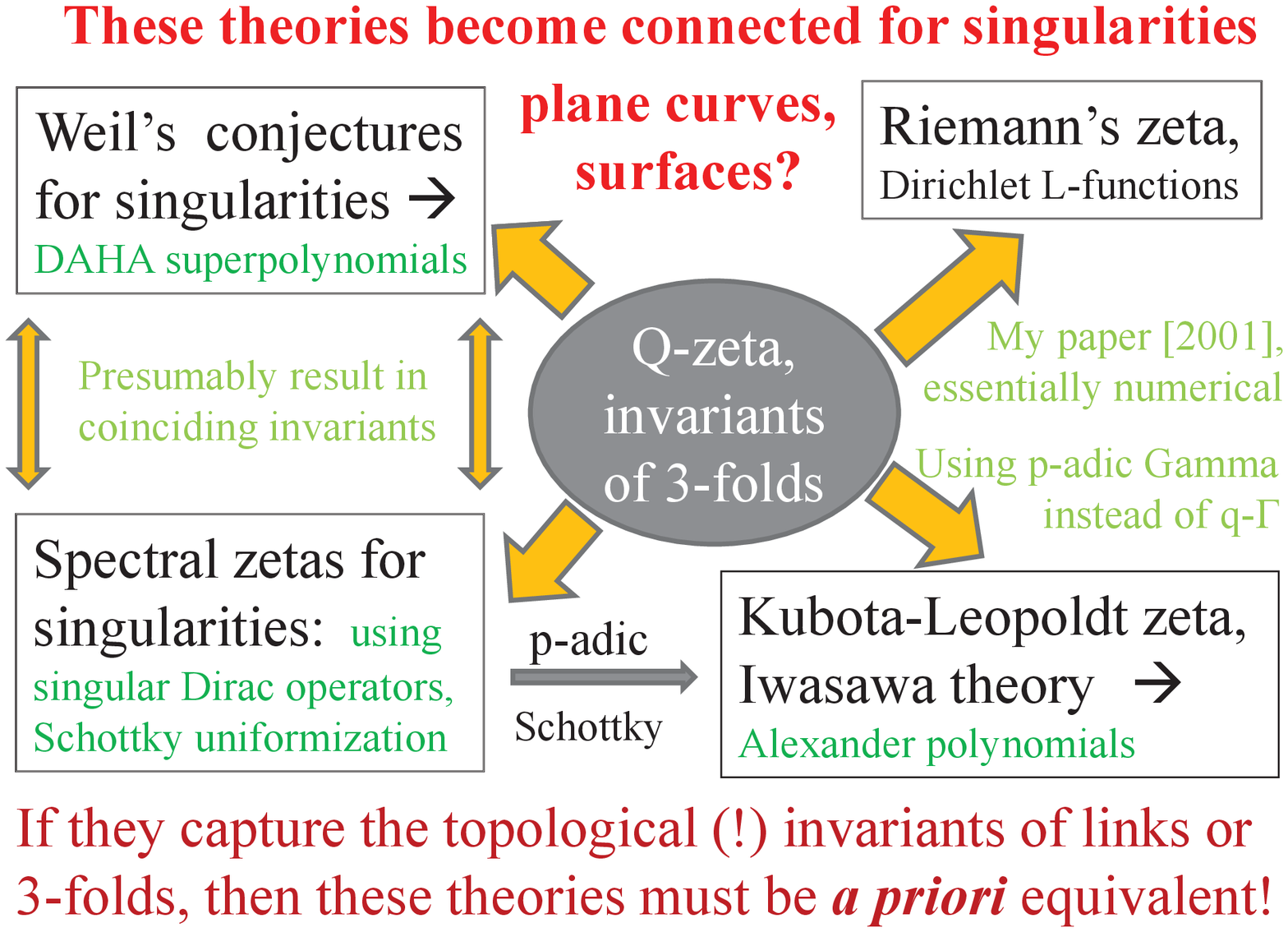}
\end{figure*}

\subsection{\bf The second figure}\label{sec:2ndfig}
Upon the switch to isolated singularities $\x$, 
the corresponding zetas can be expected to capture
the topological type for some ``good" $\x$, which
conjecturally holds for {\sf\em any} plane curve singularities. 
The ``Hasse-Weil block" of the  
$2${\small nd} figure, the upper-left one, 
is discussed in this note 
for plane curve singularities.  Let us comment a bit on the
``spectral block".

Algebraically, we study
the covers of $\mathbb P^2$ ramified at $\x$.
This is closely related to the {\sf\em Schottky uniformization} 
in the variant due to Tate-Mumford. 
Generally, it is for ``relative
curves" over
$\C[[z]]$
or $\Z_p$ ($2$-dimensional schemes), with a singular rational curve 
as the closed fiber.  Basically, 
the monoidal transformations are performed making
the closed fiber a divisor with normal crossings. We need this for
plane curve singularities, as closed fibers,
and their deformations.

More analytically, one can consider the 
{\sf\em versal deformations} of $\x$ and the invariants
of the zeta-function of the {\sf\em Dirac operator}
upon  the corresponding {\sf\em Schottky uniformization}.
Let us provide one reference:\,
{\sf\em ``Zeta functions that hear the shape
of a Riemann surface"} (Cornelissen-Marcolli, 2008).
This is in the smooth case, but Schottky uniformization
is compatible with the passage to singularities. 
For plane curve singularities $\x$, the corresponding 
Selberg-type zetas are expected to be related to
the superpolynomials of $\x$. 
\vfil

Let me mention
the $p$-adic uniformization of modular curves associated
with division algebras over totally real
number fields (my PhD thesis). 
It was used in  the {\sf\em $\vep$-conjecture} and
Ribet's theorem (in the proof of Last Fermat Theorem).
See, for instance, recent  
{\sf\em ``On the p-adic uniformization of quaternionic Shimura
curves"} (Boutot-Zink, 2022). 
This is an example of the
role the $p$-adic block can play, which also shows the potential of
the Schottky uniformization. 


\vskip 0.2cm
{\sf\em The middle oval}. The DAHA approach provides
a natural way of adding $q$ to the 
Riemann's  $\ze(s)$ and Dirichlet $L$-functions, though without the
symmetry $s\leftrightarrow 1-s$. 
The corresponding $q$-deformed
zeros seem more regular at $0\!<\!q\!<\!1$. If true, 
then the stochastic behavior of the classical ones can
be due to $q\!\to\! 1$.
\vskip 0.2cm

I defined such $q$-analogs  in 
{\sf\em ``On q-analogues of Riemann's zeta function"} (2001)
for $A_1$. The starting point was the
$q$-Mehta-Macdonald formula, which is a product formula
in terms of $q$-Gamma functions for  
$\int_{\imath \R^n} \ga(x)\, \mu(q^x)dx$ for the Gaussian
$\ga(x)=q^{-x^2/2}$ and the measure-function
$\mu(q^x)\!=${\large $\prod_{\al,j\ge 0}
\frac{\, (1-q^{(x,\al)+\nu_\al j})\, (1-q^{-(x,\al)+\nu_\al(j+1)})}
{(1-t_{\al} q^{-(x,\al)+\nu_\al j})
(1-t_\al q^{(x,\al)+\nu_\al(j+1)})},$}  
making the $E$-polynomials pairwise orthogonal. 
Here  
$\al$ are positive roots of a given reduced
irreducible root system $R\subset \R^n$,
normalized by the condition $(\al_{sht},\al_{sht})=2$; we set
$\nu_\al=1$ for short $\al$, and $\nu=2,3$ for long roots. Let
$t_{sht}=q^{k_{sht}}$ and 
$t_{lng}=q^{k_{lng}}$; $t_{\al}$ depends only on $|\al|$.
Also, let
$\rho_k=\frac{1}{2}\sum_{\al>0}k_\al \al$ and 
$\mu_1=\mu/\text{CT}(\mu)$ for the constant term functional
CT for Laurent series in terms of $X_b$ ($\mu_1$ is such).
 From now on, we assume that  $0\!<\!q\!<\!1$ and  $t_\al>1$.
\vskip 0.2cm

The link to the DAHA superpolynomials is due to the
theorem that $\int_{\imath \R^n} f(x)\ga(x)\, \mu(q^x)dx$ 
is proportional
to the {\sf\em coinvariant} of $f(x)$, which is
$f(x=-\rho_k)$ in suitable spaces of functions. Here the
integration can be replaced by taking CT;
see below. We note that the integration can be
over $\R^n$ in this formula:\, $\ga$ must be 
replaced by $\ga^{-1}=
q^{x^2/2}$ and the proportionality factor changes (significantly). 

Let $Z^+_n(q,t)=$ {\Large
$\frac{\int_{\vep+\imath R^n}\ga(x)/\left(1+ \ga(x)\right)\,\mu(q^x)dx}
{\int_{\vep+\imath R^n} \ga(x)\,\mu(q^x)dx}$ .}
Due to the {\sf\em Stirling-Moak formula} (Moak, 1984):\,
lim${}_{q\to 1_-}Z^+_n(q,t)=
\eta(s)\equal (1-2^{1-s})\ze(s)$ for the Riemann's $\ze(s)$, where
$s=k_{sht}|R_+^{sht}|+k_{lng}|R_+^{lng}|+\frac{n}{2}$. The usage
of $\vep$ improves the range of $k_\nu$ where $Z^+_n$ is analytic.  
The basic range is $\Re k_\nu>0$, which is for $\vep=0$.
If $k_{lng}=k=k_{sht}$ and $\vep=\rho/h$ for the Coxeter
number,  then $Z^+_n(q,t)$ is analytic for $\Re k>-1/h$, which
corresponds to 
$\Re s>0$, i.e. we cover the critical strip $0<s<1$ in the limit. 
By the analytic continuation, the convergence to $\eta(s)$
holds for any $s\in \C$. The analytic continuation
is essentially by the procedure of ``picking up 
residues" due to Weyl-Arthur-Heckman-Opdam.  
\vskip 0.2cm

{\sf\em The case of  $A_n$.} Setting
 $\up=n+1$
and  $\up^\circ=-k\up$, we obtain:
{\small
$$
s=k\frac{\up(\up-1)}{2}+\frac{\up-1}{2}=
\frac{1}{2}(\up-1)
(k\up^\circ+1)=-\frac{1}{2}(\up-1)(\up^\circ -1).
$$ 
}
\ The integral
$\i^+_n=\int_{\vep+\imath R^n}\frac{\ga(x)}{1+ \ga(x)}\,\mu(q^x)dx$
for $\vep=\rho/\up$ is an analytic function for $\Re k>-1/\up$
and, accordingly,  for  $s>-\frac{1}{\up}\frac{\up(\up-1)}{2}+
\frac{\up-1}{2}= 0$. Generally,
 $\i^+_R$ defined by the same formula
for a root system $R$ is analytic if 
$(a)$\, $k_{sht}>\max\{-(\vep,\al_i), (\vep,\th_{sht})-1\}$, where
$\al_i$ are short simple roots and $\th_{sht}$ is the maximal
short root in $R_+$, 
and $(b)$\, $k_{lng}$ satisfies the analogous inequalities for
long instead of short. 

We conjecture in type $A$, that there exists 
a meromorphic function  $\z(q,t,\aa)$ in terms of $q,t,\aa$
such that
$\z(t^{-1}, q^{-1}, \aa)=\z(q,t,\aa)$ and
$\eta(s)$ is the limit $q\to 1_-$ 
of $\z(q,t=q^k,\aa=t^\up)$. Note the usage of $\aa=-a$
instead of $a$ in superpolynomials.

The superduality becomes  $k\mapsto 1/k, \up\to -k\up=
\up^\circ$ in terms of $k,\up$.
The corresponding $s$ remains fixed under this
symmetry (it must!). However, we have
a nontrivial connection between the values of $\z$ at 
$k$ and $1/k$ in the $q$-theory. For instance, the $\aa$-coefficients
of $\z$ are
(conditionally) bounded as $|k|\!\to\! \infty$ and $\Re s\!>\!0$,
or its values for super-invariant $\,\aa$, which is a 
variant of the {\sf\em Lindel\"of hypothesis}. 

\vskip 0.2cm

One can expect similar features for classical root systems.
For $C_n$ in the case  $t_{\lng}=t_{sht}$,
briefly discussed above, the {\sf\em hyperpolynomials} 
are conjectured to depend polynomially on 
$q,t^{\pm 1},\aa$, where the passage to $C_n$ is $\aa=t^{2n}$;
the superduality is $q\leftrightarrow t^{-1}, \aa\mapsto \aa$, 
as that for the $A$-series. 

In this case:\, $k_{sht}=k=2 k_{lng}$, $\up=n$, and 
$s=k\bigl(\up(\up-1)+\up/2\bigr)+\up/2$. The superduality
becomes $k\mapsto 1/k, \up\mapsto -k\up$, i.e. the same as 
for $A$; however, $s$
is different for $C$. As for $A$,
this $s$ is fixed under the superduality:
$s\mapsto \frac{1}{k}(k\up)(\frac{1}{2}+k\up)-k\frac{\up}{2}=
k(\up^2-\frac{\up}{2})+\frac{\up}{2}=s.$
\vskip 0.2cm

{\sf\em DAHA vertex.}
The rationale for the existence of $\z(q,t,\aa)$
is the following theorem. Let ${\mathbb J}_m=
\text{CT}\bigl(\Th^m(q^x)\,\mu_1\bigr)$ for 
$\Th(q^x)\equal\sum_{b\in P}q^{b^2/2+ (x,b)}$.
The latter is  $\ga(x)$ ``presented"
as a Laurent series in terms of  $X_b=q^{(x,b)}$; its defining
property is the $P$-periodicity of  
$\ga(x)^{-1}\Th(q^x)$. This series naturally occurs
when we switch from $\int_{\vep+\imath \R^n} \{\cdots\} \mu(q^x)dx$ 
to the integration
over the periods of $\mu(q^x)$. For instance, 
$\int_{\vep+\imath R^n}\ga(x)\,\mu_1(q^x)dx$ 
coincides with $\mathbb J_1=
\text{CT}\bigl(\Th \mu_1(q^x)\bigr)$ up to a simple factor.
To see this, 
replace $\ga(x)$ by $\sum_{b\in Q}\ga(x+2\pi \imath\,\log(q)\,b)$
for the {\sf\em root lattice} $Q$, use the 
functional equation for $\Th$, and then switch to CT.

Given $m\ge 1$, the claim is that 
${\mathbb H}_m(q,t,\aa)$ exists 
such that 
${\mathbb H}_m(q,t,\aa=t^{n+1})={\mathbb J}_m/{\mathbb J}_1^m$ 
for any $A_n$, and the superduality holds:
${\mathbb H}_m(t^{-1},q^{-1},\aa)={\mathbb H}_m(q,t,\aa)$
(without any $q,t$-factors).

Technically, this theorem follows from the 
$\aa$-stabilization and super-invariance of 
{\large
$\frac{q^{b^2/2+k(\rho,b)}}{\lan \p_b,\p_b\ran_1}
$
} 
and those for $\p_b(q^{-c-k\rho})$ for any $b,c\in P_+$, 
where we set: $\p_b(X)\equal P_b(X)/P_b(t^{\rho})$ for
the Macdonald polynomials $P_b$, and  
$\lan f,g\ran_1\equal\text{CT}
\bigl(f(X)g(X^{-1})\mu_1(q^x)\bigr)$.
\Yboxdim5pt
For instance, let $b=\yng(2,1)\,= \om_1+\om_2,\,
c=\yng(1)\,=\om_1$. Then $\p_b(q^{-c-k\rho})=
1+\frac{(q^{-1}-1)(1-t)}{1-\aa}\bigl(1\!+\!q^{-1}\!+\!t\bigr)$,
which is, indeed, super-invariant, as well as:
{\large $\frac{q^{b^2/2+(b,\rho_k)}}{\lan \p_b,\p_b\ran_1}=$}
$$
\left(\frac{q}{t}\right)^{5/2}\aa^{-3/2}(1\!-\!\aa)
\frac{(1\!-\!t^{-1}\aa)(1\!-\!q\aa)(1\!-\!t^{-1}q\aa)
(1\!-\!t^{-1}q^2\aa)(1\!-\!t^2q\aa)}
{(1-q)^2(1-t^{-1})^2(1-qt^{-2})(1-q^2t^{-1})}.
$$
\Yboxdim6pt

These quantities are the
key in the 
theory of {\sf DAHA vertex} due to the author and Danilenko.
Namely, the series 
$\text{CT}\bigl(\Th^m(q^x)\mu_1(q^x)\bigr)$ can be expressed
in terms of products of them using the expansion $\Th(q^x)=$
{\large $\sum_{b\in P_+}\frac{q^{b^2/2+k(b,\rho)}}
{\lan \p_b,\p_b\ran_1}$}$\p_b(X)$
$\text{CT}\bigl(\Th\mu_1\bigr).$ Here the formula for
the proportionality coefficient is the $q$-Mehta-Macdonald identity:
$\text{CT}\bigl(\Th\mu_1  \bigr)=$  
\noindent
{\large $\prod_{\al\in R_+}\prod_{ j=1}^{\infty}
\Bigl(\frac{1-t_\al^{-1} q_\al^{(\rho_k,\al^\vee)+j}}
{1-q_\al^{(\rho_k,\al^\vee)+j}}\Bigr)$}.
It follows from
the theorem that the symmetric form 
$\lan f,\,g\,\Th(q^x)\,\ran_1$ 
corresponds to a {\sf\em Shapovalov-type} anti-involution of
$\HH$ and
is unique such up to proportionality.

When $t=0$, our ${\mathbb J}_m$ become
generalized {\sf\em Rogers-Ramanujan} series.
Actually, $\Th$-functions ``with characteristics" 
are needed for them, where $X_b\mapsto \xi(b)X_b$ for
characters $\xi: P \to \C^*$, which are via
$P/Q=\Z_{n+1}$ for $A_n$. 
See our paper with Boris Feigin (2012).
The case $t=0$ is incompatible with the super-invariance, but 
we obtain ``instead" the modularity of the resulting $q$-series.
\vfil

Without going int detail,  
${\mathbb H}_m$ can be interpreted as 
 invariants of the {\sf\em Lens
spaces} $L(m,1)$; they are sums of 
colored superpolynomials for the $m$-chains
of consecutive unknots with
linking numbers $-1$ between the neighboring ones. 
Not much is published in this direction. The necessary
DAHA theory is still in progress. Let me mention
{\sf\em ``BPS spectra and 3-manifold invariants"}
 (Gukov-Du Pei-Putrov-Vafa, 
2017).

The next level is the passage to $q$-zeta: 
when we consider special generating 
functions for the family  $\{{\mathbb H}_m\}$ and
more general ones.  Any ``geometric" interpretation of this 
passage can be valuable. There are quite a few 
challenges, including the following one.
\vskip 0.2cm

{\sf\em Higher theta-functions.}
Given $m$, the passage
from the imaginary integration in
$\int_{\vep+\imath R^n}\frac{\ga(x)}{1+ \ga(x)}\,\mu(q^x)dx=
\int_{\vep+\imath R^n}\sum_{m=1}^\infty
(-1)^{m-1}\ga^m(x)\,\mu(q^x)dx$ to CT results in 
$\ga^m(x)\mapsto \Th(q^{mx})$. The latter function
has the same multiplicator as $\Th^m(q^x)$ 
upon the action of $P$ and behaves as 
$\sim q^{-m x^2/2}$ when $q\to 1_-$ with a proper proportionality
factor. It is quite different from $\Th^m(q^x)$ for $m\ge 2$.  
To employ the theorem on super-invariance
of the quantities above, one needs either a counterpart of 
$\Z(q,t,a)$ defined 
in terms of $\{\Th^m(q^{x})\}$, which seems doable, or
the theory of {\sf\em DAHA vertex} based on 
$\Th(q^{mx})$ instead of  $\Th^m(q^{x})$. 

More generally, theta-functions of level $m$, those
with the same multiplicator as for $\Th^m$, are in 1-1
correspondence with the {\sf\em DAHA coinvariants} of level $m$. 
The latter are defined
algebraically using the action of $\tau_+^m$ in DAHA. 
We note that finding relations between theta-functions of arbitrary
levels $m$ is a subtle algebraic problem even for $A_1$. 

\vskip 0.2cm 
{\sf\em Further comments.}
The convergence of $Z^+_n(q,t)$  above will become
to $\ze(s)$ instead of $\eta(s)$ when we switch to $Z_n^-(q,t)$
with $(1-\ga(x))$
instead of $(1+\ga(x))$. However, this will hold only for $\Re s>n$.
It will diverge otherwise, which can be fixed upon multiplication
by proper $\om^\bullet$ for $q\equal e^{-1/\om}$;
the limit becomes then some
$\Ga$-type function (no zeta!). 
\vfil

For instance, let 
$\i^-(k)=\frac{1}{i}\int_{1/2+\imath \R} \bigl(q^{x^2}-1\bigr)^{-1} 
\mu(q^x)dx$ for $n=1$.
Then $\bigl(\frac{\om}{4}\bigr)^{k-1/2}\i^-(k)$ converges to 
$\Ga(s)\ze(s)$ for $\Re k\!>\!\frac{1}{2}$, where 
$s\!=\!k\!+\!\frac{1}{2}$.
However,  the limit of 
$\om^{2k-1}\i^-(k)$ becomes
$\tan(\pi k)\Ga(k)^2$ when $\Re k<\frac{1}{2}$. Generally,
$\om^{s-n}\i^{\pm}(k)$ converges to $\ze(s)$ or $\eta(s)$ times
proper products of $\Ga$-factors
 when $\Re s>n$ in the case of  $\i^{-}$ 
and $\Re s>0$ for $\i^{+}$. Using the analytic continuation,
$s$ can be arbitrary in the second case. 

Interestingly, 
the $q$-zeros of $\i^-(k)$ with $\Re k\sim 0$ (not in the
range $\Re k>\frac{1}{2}$)
approach the classical ones for $q$ sufficiently
close to $1_-$ (but not too close!), and then $\i$ 
slowly begin to ``switch" to the Gamma-limit as $q$ continues 
to approach $1_-$. This can be clearly seen numerically.
\vfil
\vskip 0.2cm

This can be potentially related to the {\sf\em Gram law}.
We  change $\i^-(k)$ to $\tilde{\i}^-(k)$ with the limit 
to the classical $\tilde{\ze}(s)=
\pi^{-s/2}\Gamma(\frac{s}{2})\ze(s)$, invariant under $s\mapsto 1-s$;
this limit is only  for $\Re s>1/2$ and upon $\om^\bullet$. 
Then the zeros of proper linear combination of $\i^-(k)$ and its
complex conjugation approach those 
of $\ze(k+1/2)$ for $\Re k\sim 0$ and large $\om$ (but not
too large), and then eventually tend to properly adjusted 
{\sf\em Gram points} as $\om$ increases.     
 This argument is 
of course qualitative. So is the Gram law: the zeros 
of $\ze(1/2+k)$ ``mostly"
alternate with the Gram points.  
\vskip 0.2cm

We note (again) that the invariance of $\tilde{\ze}(s)$ 
under  $s\mapsto 1-s$ 
fails for the $q$-deformations. However, a variant of 
``$q$-RH" can hold for $A_1$. 
This was discussed in my paper mostly
for the {\sf\em sharp $q$-zeta}
introduced below. A qualitative version is that
the tendency is strong for 
the {\sf\em sharp $q$-deformations} of the classical 
zeros 
to stay in the half-plane  $\Re s>1/2$
(in certain horizontal strips depending on $q<1$). The 
corresponding half-plane
becomes $\Re s<1/2$ for the imaginary integration. 
 
A direct counterpart of RH can hold too.
This requires the consideration of a proper
linear combination of $\i^+(k)$ and $\i^+(-k)$.  
The $q$-zeros of the latter we were able to find 
satisfy $\Re k=0$. The case of the simplest 
{\sf\em sharp} $q$-$L$-function is touched upon below. 

\vskip 0.2cm

Peter Sarnak noted once that many applications
are based  on the
absence of zeta-zeros with $\frac{1}{2}< \Re s <1$.  This
is, basically, what we see for ``small" $q$-zeros.
However, there was (and there is) 
uncertainty when the corresponding
neighboring  zeros of the classical
$\ze(s)$
are getting ``too close". Namely, the linear approximations
of sharp $q$-deformations of such ``unusual" zeros of $\ze(s)$
can be with $\Re s>1/2$ (in the opposite half-plane).
We think that the linear approximations
can be irrelevant for such $q$-zeros; theoretically,
the convergence
of the corresponding Taylor expansions is not known.
 Let me quote Harold Edwards: 
``the existence of nearly coincident zeros 
must give pause to even the most convinced believer"
(his {\sf\em ``Riemann's Zeta Function"}).


\subsection{\bf Sharp {\em q}-zeta}
As above,\, $q=\exp(-1/\om)$ for $\om>0$;
\def\shrpi
{\circ\kern-13.9pt\sqsubset\joinrel\mathrel
{\raise3.2pt\hbox{$\leftarrow$}}\kern-17.5pt
{\raise-3.2pt\hbox{$\rightarrow$}}}
let  $\si\equal\sqrt{\pi \om/2}$. The
integration path will be now
$\shrpi
^{\infty+\si i}_{\infty-\si i }$ 
around zero. 
For $A_1$ and 
$\de_k(x;q)\equal\prod_{j=0}^\infty \frac{(1-q^{j+2x})(1-q^{j-2x})}{
(1-q^{j+k+2x})(1-q^{j+k-2x})}$, the symmetric variant of $\mu$,  
the function 
$\mathfrak{Z}^{\sqsubset}_{q}(k) \equal 
\frac{1}{2i}\,{\raise1pt\hbox{\tiny$\sqsubset$}}
\kern-9pt\int^{\infty+\si i}_{\infty-\si i } 
\frac{\de_k(x;q)}{1+q^{-x^2}}\, dx
$
is analytic in the horizontal strip 
$K^{\sharp}\!=\!\{-2\si\!<\!\Im k\!<\!+2\si\}$ as  $\Re k>-1/2.$
Its meromorphic continuation to {\sf\em all} 
$k\in\C$ via Cauchy's theorem, the {\sf\em sharp $q$-zeta}, is:
{\footnotesize
\begin{align*}
&\mathfrak{Z}_{q}^{\sharp}(k) = 
-\frac{\om\pi}{ 2}\prod_{j=0}^\infty
\frac{(1-q^{j+k})(1-q^{j-k})}{
 (1-q^{j+2k})(1-q^{j+1})}\times\\
&\sum_{j=0}^\infty \frac{(1-q^{j+k})q^{-kj}}
{(1-q^{k})(q^{-\frac{(k+j)^2}{4}}+1)}
\prod_{l=1}^j
\frac{1-q^{l+2k-1}}{
 1-q^{l}}.
\end{align*}
}

\!\!\!\!It has poles at $\{-\frac{1}{2}\!-\!\Z_+\}$ in  $K^{\sharp}$.
This strip is between the first zeros of 
$1+q^{-\frac{k^2}{4}}.$
For all $k$ apart from the poles,
 $\lim_{\om\to \infty} (\frac{\om}{4})^{k-1/2}
\mathfrak{Z}_{q}^{\sharp}(k)$
$=\sin(\pi k)(1-2^{\frac{1}{2}-k})
\Gamma (k+\frac{1}{2})\zeta (k+\frac{1}{2}).$

Given a classical zero $k\!=\!z$ of $\ze(1/2+k)$, let us assume that 
its $\sharp$-deformation $z^\sharp(\om)$ 
exists and is differentiable with respect
to $\varpi=1/\om$. Then the formula for its linear approximation 
$\tilde{z}^\sharp(\om)$ is as follows:
\noindent
$\tilde{z}^\sharp(\om)=
z\bigl(1 -\frac{4(z+\frac{1}{2})\zeta_+(z+\frac{3}{2})-(z-1)
\zeta_+(z-\frac{1}{2})
}{ 12 \om\zeta'(z+\frac{1}{2})(1-2^{\frac{1}{2}-z}) }\bigr).
$
Thus, the linear $\varpi$-approximations of all 
classical zeros  $z$ exist if and
only all of them are simple, 
an interesting interpretation of the classical 
conjecture. 

A similar approach can be employed  for 
$\mathfrak{Z}_{q}^{\sharp}(k;d)$ with 
for $q^{-dx^2}$ instead of $q^{-x^2}$ in
$\mathfrak{Z}_{q}^{\sharp}$
and for the   
{\sf\em sharp $L$-functions} $\mathfrak{L}_q^\sharp(k;d)$ that are 
for 
$\frac{q^{x^2/2}-q^{-x^2/2}}{q^{(d+1) x^2/2}-q^{-(d+1)x^2/2}}.$
The $q$-deformations of
$L$-functions are somewhat simpler to analyze numerically; the
$1${\small st} $L$-zeros are smaller than those of $\ze(s)$. 
The usage of $d$ has some practical advantages too. 

\vskip 0.1cm

Taking the classical $z=14.1347i$ and $\om=750$ for  
$\mathfrak{Z}_{q}^{\sharp}(k;2)$:
$$
 z^\sharp=0.1304 + 14.1450i,\ 
\tz^\sharp=0.1302 + 14.1465i.
$$
Other zeros in $K^\sharp$ for $\om=750, d=2$ are:
{\small
\begin{align*}
& \ \ zeta\ \ \ \ \ \ \ \ \ \, sharp-zeta\ \ \ \ \ \ \ \, 
linear\ approx.\\
& 21.0220i\ \ \ 0.3514 + 21.0702i\ \ \ 
0.3504 + 21.0771i \\
& 25.0109i\ \ \ 0.5641 + 24.9586i\ \ \
0.5745 + 24.9643i \\
& 30.4249i\ \ \ 0.9046 + 30.4014i\ \ \
0.9134 + 30.4077i \\
& 32.9351i\ \ \ 1.1051 + 33.0341i\ \ \
1.0998 + 33.0854i \\
& 37.5862i\ \ \ 1.6449 + 37.9660i\ \ \
1.7675 + 38.1895i \\
& 40.9187i\ \ \ 1.9080 + 40.8119i\ \ \
1.9141 + 40.7816i \\
& 43.3271i\ \ \ 2.2860 + 43.2485i\ \ \
2.4497 + 43.3138i \\
& 48.0052i\ \ \ 2.9259 + 47.8424i\ \ \
3.1103 + 47.5578i.
\end{align*}
}
\!\!There is a clear tendency for $z^\sharp$ to move to the right.
If true, this would give the classical RH.
A direct $q$-counterpart of RH was the observation at the end of
paper ``On $q$-analogs ..." that
several (not too many) ``small" zeros  of 
$\mathfrak{L}_q^\sharp(k;d)-
\mathfrak{L}_q^\sharp(-k;d)$ were all with $\Re k=0$.
This is within the corresponding strips. The convergence is 
very good, including large $\om$, but the 
calculations become involved  for large $\Im(k)$; this restricted
our simulations.  
The passage from $A_1$ to $A_n$ and the stabilization add
{\sf\em superduality} to this approach. 


Our calculations indicate that 
the zeros of the $q$-deformed zeta-functions and $L$-functions
become  more ``regular" for $q<1$ than
the corresponding classical zeros.
It is expected that the zeros of Riemann zeta function 
in the critical strip are totally random subject to
the distribution for the 
eigenvalues of random Hermitian
matrices (Dyson, Montgomery, Odlyzko and others).
This can be due to the limit $q\to 1_-$ of ``relatively regular"
$q$-zeros. Similar quasi-classical limits, in physics and
mathematics, are known to create chaotic  
behavior. 
\vskip 0.2cm


{\sf\em Two bottom blocks.}
Concerning the ``$p$-adic block", it is expected that
there is a $p$-adic DAHA theory, where $q$-Gamma
functions are replaced by their $p$-adic counterparts.
This theory is important and doable. Technically, the
cyclotomic $q$-Gauss-Selberg sums in the DAHA theory,
where $q$ are roots of unity, will be replaced
by the classical {\sf\em modular} Gauss-Selberg sums,
those over $\F_q$. The 
$q$-Mehta-Macdonald formulas is expected
then to become in terms of 
the $p$-adic  Gamma. 

\vskip 0.2cm
The ``spectral block" has been partially discussed. 
The $\rho_{ab}$-invariants provide important links.
Generally, $\rho_\xi(M)$ is the value of at $s=0$ of
$\tilde{\eta}_M(s)-n\,\eta_M(s)$ 
for odd-dimensional closed oriented
Riemannian manifolds $M$, where $\eta_M(s)$ is  
due to  Atiyah-Patodi-Singer and ``twisted" $\tilde{\eta}_M(s)$
is defined for the flat bundle over $M$ associated with 
the representation $\xi: \pi_1(M)\to U_n$.
We take $M=S^3\setminus K$ for a knot $K$ and  
the abelianization of $\pi_1(M)$ as $\xi$. Thus, our $\rho_{q,t}=
R_K(q,t,a\!=\!-\!1/q)$ from Section \ref{sec:rho}
can be considered as a $q$-version
of $\tilde{\eta}_M(s)$. 
One can expect a ``triply-graded" categorification of 
$R(q,t,a)$. 
Recall that they are  defined 
in terms of $\h^{mot}$, which conjecturally coincide with
superpolynomials in any other theories and, also, with
the generalized $L$-functions of plane curve singularities. 
Thus, 
{\sf\em quasi-rho invariants} connect 
the ``Weil block" with
the ``spectral block" and then, potentially, with
the ``$p$-adic block"  via the $p$-adic Schottky uniformization.

\vskip 0.2cm
\vfil

{\sf\em Further perspectives.}
Assuming that zeta functions of singularities
$\x$ and their corresponding 
$a,q,t$-versions are topological
invariants/moduli of some sort, one can expect 
{\sf\em a priori} links between the Hasse-Weil zetas, 
Selberg's zetas and
$p$-adic zetas for such $\x$. If this is true, then these zetas are 
different invariants of the same singularity and 
must be connected as such. 
This is expected in the case of plane 
curve singularities.
\vfil

A program is to switch 
 from plane curve singularities to
{\sf\em surface} singularities 
serving Seifert 3-folds and more general
{\sf\em plumbed manifolds}. An example is our $q,t,a$-deformation
of the Dirichlet $\eta(s)$, which is, 
essentially, 
$\sum_{m=1}^\infty (-1)^{m-1} {\mathbb H}_m(q,t,\aa)$, an alternating
sum of certain invariants of Lens spaces $L(m,1)$. Thus, 
the invariants
of the latter 
``replace" $\frac{1}{m^s}$ in the classical  
$\eta(s)=\sum_{m=1}^\infty (-1)^{m-1}\frac{1}{m^s}$ in our approach. 
\vfil

One can try to interpret such deformed sums geometrically: by looking for
a manifold  $M$ such that its ``triply-graded"
homology reduces to the sums of those over its
special submanifolds, Lens or Seifert spaces. 
In geometry, this can be due to the {\sf\em localization} in
certain (co)homology of $M$ or via the  count of {\sf\em closed} 
totally-geodesic submanifolds in $M$ (with proper weights).
Also, it is not rare when some 
invariants of fibered spaces $M$ reduce to sums over proper
{\sf\em
special} fibers. However, here we have {\sf\em infinite} sums
and quite involved homology theories. 


\vskip 0.2cm
We note that the superduality  
holds for any (convergent) series in the 
form $\sum_{m=1}^\infty  c_m(q,t)
{\mathbb H}_m(q,t,\aa)$ provided
the super-invariance of $c_m(q,t)$.  The superduality alone does not
require specific $\lim_{q\to 1_-}c_m(q,t)=\pm 1$ 
from $\eta(s)$ 
and those for the $L$-functions.  However, the fact that we
``lift" the symmetry $k\mapsto  1/k, \upsilon\mapsto -\upsilon k$ 
of $s$
in $\eta(s)$ to the fundamental superduality 
$q\leftrightarrow t^{-1}$ in the theory
of superpolynomials is of interest. 


\subsection{\bf Strongly-polynomial count}
There are significant  restrictions for the types of singularities
and their usage in the 2{\small nd} figure. 

First, the motivic zetas are supposed to be {\sf\em 
topological} invariants. This is not granted in their
definition, which is in terms of the corresponding
singularity ring. They are of course  
of ``discrete nature", which makes them 
potentially topological; however, the topological invariance
is a conjecture even for general plane curve singularities.

Second, we need to
check that $\x$ can be defined over $\Z$ 
 within its topological type and has 
good reductions for almost all prime $p$. 

Third, the varieties of ``standard modules" 
and ideals of finite colength in 
the corresponding local rings 
must be of 
{\sf\em strongly-polynomial count\,}: 
the number of their points over $\F_q$
must depend  polynomially on $q$. 

The latter condition is very restrictive.
It holds if a variety is paved by 
{\sf\em configurations of affine spaces}, their unions
and differences in a bigger affine space. This is conjectured
for the {\sf\em Piontkowski cells} in our varieties
of standard modules $\j_{\ell}$
and their multi-rank generalizations.   
We mention that 
affine Springer fibers can be {\sf\em not}
of strong polynomial count (unless in type $A$).
There is an example of Bernstein-Kazhdan where the
counting their $\F_q$-points involves  
zeta-functions of elliptic curves over $\F_q$, 
certainly not $q$-polynomial.
See Appendix to {\sf\em ``Fixed point varieties on affine 
flag manifolds"} by Kazhdan-Lusztig (1988).
\vskip 0.2cm

Under these conditions,
the corresponding 
$\zeta_\x(q,t,a)$ or $L_\x(q,t,a)$  can be expected  
powerful topological invariants of $\x$. Presumably, they
can capture the topological types of $\x$ in some cases.
They certainly do this for 
plane curve singularities $\x$; however, the topological
invariance of the motivic superpolynomials $\h^{mot}$ is 
a conjecture. 

Plane curve singularities provide 
the main example by now.  There is a bunch of constructions, theorems
and conjectures in this case. One of the purposes of this work
is to show that the theory of their superpolynomials
can be developed in various directions, 
which, presumably, includes  
isolated surface singularities of toric type and the corresponding
Seifert-type manifolds. 



Needless to say that isolated singularities are (and
always were) among the key objects
of algebraic geometry.  Smooth projective manifolds 
proved to be very helpful in their study, but they are not really 
necessary
for many aspects of singularity theory.
We try to do as much  as we can directly
in terms of the singularity rings.

\vskip 0.2cm
 
{\sf\em Kn\"orrer's periodicity.}
In topology, there is a fundamental connection
between knots/links and 3-folds. Its certain 
algebraic counterpart is 
the Kn\"orrer's periodicity for singularities: 
a connection between
the plane curve singularities $W(x,y)\!=\!0$ and the
ones given by the equations $u^2\!=\!W(x,y)$ is its important 
part. Actually,  space singularities and
5-folds fit this picture too, but only ``good" ones.
For instance, the singularities in the form
$uv\!=\!W(x,y)$ naturally occur here;
such Calabi-Yau threefolds were considered by Vafa-Dijkgraaf.  

The expected connection between the {\sf\em superpolynomials} of 
algebraic links 
and the {\sf\em superseries} 
of the corresponding Seifert spaces  resembles
that between zeta-polynomials  and
the corresponding $L$-functions $L_\Phi(s,\chi)$, we began with.
This is very far from being exact.  Our superpolynomials
are much simpler than the zeta-polynomials, and 
our superduality is that from the Hasse-Weil 
functional equation, very different from
that for $L$-functions.  The passage from the {\sf\em superseries} 
of Lens spaces to $q,t,a$-deformations of Dirichlet $L$-functions 
can be viewed as an attempt to unify these two theories.
\vskip 0.2cm

Let me mention (again) that this note is
very incomplete concerning the names and
contributions; only very few papers are mentioned.
It is focused mostly on superpolynomials
and some perspectives of their theory. We tried to outline 
some number theoretical perspectives of this direction and
possible physics connections. The exposition is sketchy and
speculative in several places.  There are various 
omissions; for instance, 
we do not discuss much recent developments, even 
those directly related to the topics we touched upon. 
\vskip 0.2cm

To conclude, Manin's works and his vision of the role
of number theory greatly influenced a lot of people, 
certainly all his 
students. We thank very much Yuri Tschinkel, Michael Finkelberg
and the referee.

\end{document}